\RequirePackage{rotating}
\documentclass{mcom-l}
\usepackage{array}
\usepackage{amsmath,subfigure}
\usepackage{fix-cm}
\usepackage{amsfonts}
\usepackage{mathrsfs}
\usepackage{multirow}
\usepackage{anyfontsize}
\usepackage[dvips]{epsfig,graphicx}
\usepackage{stmaryrd}
\usepackage{listings,url,verbatim,color,bm}
\usepackage{lipsum}
\usepackage{algorithm}
\usepackage{algorithmicx}
\usepackage{algcompatible}

\newtheorem{problem}{Problem}
\newtheorem{theorem}{Theorem}[section]
\newtheorem{lemma}[theorem]{Lemma}

\theoremstyle{definition}

\theoremstyle{remark}
\newtheorem{remark}[theorem]{Remark}

\numberwithin{equation}{section}

\begin{document}

\title[Fixed-Rank Approximations in Tucker Format with Adaptive Shifts]{Fixed-Rank Approximations in Tucker Format with Adaptive Shifts based on Fast Randomized Algorithms}

\author{Maolin Che}
\address{School of Mathematics and Statistics and State Key Laboratory of Public Big Data, Guizhou University, Guiyang, Guizhou, P. R. of China 550025}
\email{mlche@gzu.edu.cn and chncml@outlook.com}
\thanks{Maolin Che was supported by the National Natural Science Foundation of China (No. 12561095), the Special Posts of Guizhou University (No. [2025]06), and Guizhou Provincial Major Project of Basic Research Program (Qiankehe Foundation VZD[2026]001).}

\author{Yimin Wei}
\address{School of Mathematical Sciences and Key Laboratory of Mathematics for Nonlinear Sciences, Fudan University, Shanghai, P. R. of China}
\email{ymwei@fudan.edu.cn}
\thanks{(Corresponding author) Yimin Wei was supported by the National Natural Science Foundation of China under grants U24A2001 and 12271108, and the Science and Technology Commission of Shanghai Municipality under grant 23JC1400501.}

\author{Chong Wu}
\address{Department of Electrical Engineering, City University of Hong Kong, Kowloon, Hong Kong SAR, P. R. of China}
\email{chongwu2-c@my.cityu.edu.hk}

\author{Hong Yan}
\address{Department of Electrical Engineering, City University of Hong Kong, Kowloon, Hong Kong SAR, P. R. of China}
\email{h.yan@cityu.edu.hk}
\thanks{Hong Yan was supported by the Hong Kong Research Grants Council (Project 11201825), and City University of Hong Kong (Projects 9610034 and 9610460).}

\date{}



\begin{abstract}
Randomized numerical linear algebra provides scalable approximations for tensor decomposition by connecting theoretical advances with practical computation. This paper introduces fast randomized algorithms for solving {\color{red}fixed-rank approximations in Tucker format}, through the integration of adaptive shifted power iterations. The proposed algorithms enhance randomized variants of truncated high-order singular value decomposition (T-HOSVD) and sequentially T-HOSVD (ST-HOSVD) by incorporating dynamic shift strategies. These accelerate convergence by refining the singular value gap and reduce the number of required power iterations while maintaining accuracy. Theoretical analyses provide probabilistic error bounds, demonstrating that the proposed methods achieve comparable or superior accuracy compared to deterministic approaches. Numerical experiments on synthetic and real-world datasets validate the efficiency and robustness of the proposed algorithms, showing a significant decline in runtime and approximation error over state-of-the-art techniques.
\end{abstract}

\maketitle

\section{Introduction}

Tensor decomposition refers to a family of mathematical methods that factorize {\color{red}high-dimensional arrays} (e.g., tensors, hypermatrices and multi-way data) into simpler and interpretable components. It brings matrix factorizations (e.g., singular value decomposition (SVD) and principal component analysis (PCA)) to the tensor case, capturing interactions across multiple dimensions. {\color{red}Tensor decomposition has} a wide range of applications, such as, data compression, quantum computing, signal processing, natural language processing (NLP) and machine learning.

{\color{red}A tensor} is denoted by $\mathcal{A}\in\mathbb{R}^{n_1\times n_2\times \dots \times n_d}$ with entries given by
$a_{i_1i_2 \dots i_d}\in\mathbb{R}$, where $i_k=1,2,\dots,n_k$ and $k=1,2,\dots,d$. In particular, when $d=1$, a first-order tensor $\mathcal{A}$ is a vector of size $n_1$, and when $d=2$, a second-order tensor $\mathcal{A}$ is a matrix of size $n_1\times n_2$. Common types of tensor decomposition include CANDECOMP/PARAFAC (CP) decomposition, Tucker decomposition, Hierarchical Tucker decomposition, tensor-train decomposition, tensor ring decomposition, fully-connected tensor network decomposition, tensor wheel decomposition, and others. Interested readers can refer to \cite{cichocki2016tensor,cichocki2017tensor,grasedyck2013literature,kolda2009tensor,wang2024randomized,zhao2016tensor,zheng2021fully} and references therein for details about these types of tensor decomposition. In this paper, we focus on the fixed-rank approximation in Tucker format, which is a special case of finding the approximation to Tucker decomposition. We now overview the existing algorithms for {\color{red} the fixed-rank approximation in Tucker format} and consider their disadvantages.

Existing algorithms for solving the fixed-rank approximation in Tucker format can be classified into iterative algorithms and direct algorithms. Iterative algorithms include the high-order orthogonal iteration (HOOI) (see \cite[Algorithm 4.2]{de2000best}), several optimization algorithms (such as the Newton, trust-region, quasi-Newton and preconditioned coordinate descent methods) on the product of matrix manifolds (see \cite{elden2009a,hamed2024riemannian,ishteva2011best,savas2010quasi}), and tensor Krylov-type methods (see \cite{goreinov2012wedderburn,elden2022krylov,savas2013krylov}). Note that HOOI is bottlenecked by the operation called the tensor times matrix-chain (TTMc). By using \textsc{TensorSketch} (see \cite{pagh2013compressed}), Malik and Becker \cite{malik2018low} proposed two randomized versions for HOOI, which avoid the expensive cost of TTMc. Ma and Solomonik \cite{ma2021fast} introduced another sketched version for HOOI with \textsc{TensorSketch} and leveraged score sampling (see \cite{drineas2012fast}), which contains a sequence of sketched rank-constrained linear least squares subproblems. The truncated high-order singular value decomposition (T-HOSVD) and the sequentially T-HOSVD (ST-HOSVD) are two common algorithms for solving the fixed-rank approximation in Tucker format (see \cite{de2000best,vannieuwenhoven2012new}). 

Recently, several randomized algorithms based on random projection and power scheme were developed to accelerate T-HOSVD and ST-HOSVD (called randomized variants for T-HOSVD and ST-HOSVD), for example, see (\cite{ahmadi2021randomized,che2019randomized,che2023randomized,che2020computation,che2021efficient,che2025efficient-acom,che2021randomized,minster2020randomized,minster2024parallel,sun2020low,zhou2014decomposition}). {\color{red}One key of these algorithms is that power iteration on subspaces (also called subspace iteration or orthogonal iteration) with a random start (cf. \cite[Sections 7.3.2 and 8.2.4]{Golub2013matrix}) is used to improve the accuracy of approximate Tucker decompositions with a given Tucker rank.  Rokhlin {\it et al}. \cite{rokhlin2010randomized} introduced the idea of using randomized subspace iteration to obtain low rank approximations of matrices by selecting a large oversampling parameter and a small power parameter. Halko {\it et al.} \cite{halko2011finding} refactored and simplified the algorithm, and presented a complete theoretical justification for the approach. Subsequent analyses can be also found in  \cite{gu2015subspace,martinsson2020randomized}.} Another key of these algorithms lies in the selection of the random projection matrix, such as standard Gaussian matrices in \cite{ahmadi2021randomized,che2025efficient-acom,minster2020randomized,sun2020low,zhou2014decomposition}, {\color{red}the Khatri-Rao product of standard Gaussian matrices} in \cite{che2019randomized,sun2020low}, the Kronecker product of standard Gaussian matrices in \cite{che2020computation,che2021randomized}, and the Kronecker product of structured random projection matrices in \cite{che2023randomized,che2021efficient,minster2024parallel}. {\color{red}For the sake of completeness, we provide two forms of randomized variants for T-HOSVD and ST-HOSVD in Appendix \ref{dash-rthosvd:appsec1}.}

A related problem for approximating Tucker decomposition with a given tolerance is called the fixed-precision problem in Tucker format (see \cite[Problem 1.1]{che2025efficient}). Several scholars have proposed efficient algorithms for this problem. 

Ehrlacher {\it et al.} \cite{ehrlacher2021adaptive} proposed a greedy strategy to estimate the Tucker-rank of a tensor with a given tolerance. A rank-adaptive (RA) variant of HOOI (see \cite{xiao2024rank}) is proposed, which proved that the method is locally optimal and monotonically convergent. The key for adaptive variants of T-HOSVD and ST-HOSVD is to calculate the truncated SVD with a given tolerance $0<\epsilon<1$ of each mode unfolding of the corresponding tensor. When replacing the standard Gaussian vectors used in the adaptive randomized range finder (see \cite[Algorithm 4.2]{halko2011finding}) with the Khatri-Rao product of standard Gaussian vectors, Che and Wei \cite{che2019randomized} proposed adaptive randomized variants of T-HOSVD and ST-HOSVD. Minster {\it et al.} \cite{minster2020randomized} presented adaptive randomized algorithms for the fixed-precision problem of approximation of Tucker decomposition by using a version of the adaptive randomized range finder algorithm in \cite{martinsson2016randomized,yu2018efficient}. Hashemi and Nakatsukasa \cite{hashemi2023rtsms} proposed an adaptive randomized algorithm for approximating Tucker decomposition, denoted by RTSMS (Randomized Tucker via Single-Mode-Sketching). Che {\it et al.} \cite{che2025efficient} modified the adaptive randomized range finder algorithm in \cite{martinsson2016randomized,yu2018efficient} by using the Khatri-Rao product of standard Gaussian matrices (or uniform random matrices) and proposed block versions for the adaptive randomized variants of T-HOSVD and ST-HOSVD.

Simultaneous iteration (also called orthogonal iteration) is the key part in (adaptive) randomized variants of T-HOSVD and ST-HOSVD. As shown in \cite{Golub2013matrix}, shifted simultaneous iteration, combined the shift strategy with simultaneous iteration, was proposed to improve the convergence rate for computing dominant subspaces of any matrix $\mathbf{A}\in\mathbb{R}^{m\times n}$. {\color{red}Shift strategies accelerate the convergence of simultaneous iterations by enlarging its singular value gap. Reducing the ratio \(\frac{\sigma_i(\mathbf{A})}{\sigma_j(\mathbf{A})}\) for \(1\leq j<i\leq \min\{m,n\}\) speeds up subspace identification, with \(\sigma_i(\mathbf{A})\) being the i-th singular value of \(\mathbf{A}\).} Typical shift strategies comprise single-shift, double-shift, double-implicit-shift, and Rayleigh shift schemes. {\color{red}The Rayleigh shift strategy constitutes a classical acceleration technique for eigenproblems: it adopts the Rayleigh quotient as a data-driven shift parameter to rapidly refine eigenvalue and subspace estimates. This strategy exemplifies the general principle of adaptive shifting and serves as the motivation for our dynamic shift design (cf. \cite[Section 8.2.3]{Golub2013matrix})}. 

{\bf Our works.}  In this paper, we introduce the adaptive shifted strategy to resolve the core speed/accuracy trade-off in randomized variants of T-HOSVD and ST-HOSVD. Standard power iteration converges slowly when singular values are closely spaced, requiring many iterations to maintain accuracy, while shifting widens the singular value gap near the target rank to accelerate convergence greatly, enabling high accuracy with fewer power steps or better accuracy for the same computational cost. This step is central to our contribution because it makes randomized variants of T-HOSVD and ST-HOSVD both faster and more accurate than non-shifted variants, matching deterministic performance without extra overhead.  Our adaptive midpoint shift follows this spirit while remaining simple, stable, and tractable for tensor randomized algorithms.

To ensure that randomized variants of T-HOSVD and ST-HOSVD are similar to the deterministic algorithms in terms of accuracy, power iteration should be introduced to improve their accuracy. Properly setting the power parameter is a difficult issue. Using the per vector error (PVE) criterion (see \cite{musco2015randomized}) can resolve the difficulty of setting a suitable power parameter and enables automatic termination of the power iteration.
\subsection{{\color{red}Organization}}

The remainder of this paper is organized as follows. Some basic definitions and results for Tucker decomposition are introduced in Section \ref{dash-rthosvd:sec2}. In Section \ref{dash-rthosvd:sec3}, we deduce two efficient randomized algorithms for the fixed-rank approximation in Tucker format, analyze their computational complexities, and consider the choice of the power parameters in Algorithms \ref{dash-rthosvd:alg2} and \ref{dash-rthosvd:alg2:v2} based on PVE accuracy control . In Section \ref{dash-rthosvd:sec4}, we discuss the error bounds for approximating Tucker decomposition with a given Tucker-rank, which is obtained from the proposed algorithms. In Section \ref{dash-rthosvd:sec5}, several examples are used to illustrate the accuracy and efficiency of our algorithms. We conclude this paper in Section \ref{dash-rthosvd:sec6}. In the Appendix sections, we give two models for randomized variants of T-HOSVD and ST-HOSVD (Section \ref{dash-rthosvd:appsec1}), and deduce another type for Algorithm \ref{dash-rthosvd:alg2:v2} according to Algorithm \ref{dash-rthosvd:alg1:app} (Section \ref{dash-rthosvd:appsec2}).

\section{Preliminary}
\label{dash-rthosvd:sec2}
We now introduce several notational forms used in this paper. Lower case letters (e.g., $x,u,v$), lower case bold letters (e.g., $\mathbf{x},\mathbf{u},
\mathbf{v}$), bold capital letters (e.g., $\mathbf{A},\mathbf{B},\mathbf{C}$) and calligraphic letters (e.g., $\mathcal{A},\mathcal{B},\mathcal{C}$) are used, respectively, for representing scalars, vectors, matrices and tensors. This notation is used consistently for the lower-order parts of a given structure. For example, the entry with the row index $i$ and the column index $j$ in a matrix ${\bf A}$, that is, {\color{red}$({\bf A})_{ij}$ or $\mathbf{A}(i,j)$}, is represented as $a_{ij}$ (also $(\mathbf{x})_i=x_i$ and $(\mathcal{A})_{i_1i_2\dots i_d}=a_{i_1i_2\dots i_d}$).

{\color{red}The symbols $\mathbf{A}^\top$, $\mathbf{A}^{\dag}$, $\|\mathbf{A}\|$ and $\|\mathbf{A}\|_2$ are denoted, respectively,  by the transpose, the Moore-Penrose pseudoinverse, the Frobenius norm, and the spectral norm of $\mathbf{A}\in\mathbb{R}^{n_1\times n_2}$.} We use $\mathbf{I}_n$ to denote the identity matrix in $\mathbb{R}^{n\times n}$. An orthogonal projection $\mathbf{P}\in\mathbb{R}^{n\times n}$ is defined as ${\bf P}^2={\bf P}$ and ${\bf P}^\top={\bf P}$. A matrix $\mathbf{Q}\in\mathbb{R}^{n\times r}$ with $r< n$ is orthonormal if $\mathbf{Q}^\top\mathbf{Q}=\mathbf{I}_r$. We use $\mathbb{S}_d$ to denote the $d$-th-order symmetric group on the set $(1,2,\dots,d)$. We use $O(\cdot)$ to refer to the class of functions whose growth is bounded above and below by a constant.

For any $k$, the mode-$k$ product of $\mathcal{A}\in \mathbb{R}^{n_1\times n_2\times \dots \times n_d}$ by ${\bf B}\in \mathbb{R}^{m_k\times n_k}$, denoted
by $\mathcal{A}\times_{k}{\bf B}$, is $\mathcal{C}\in\mathbb{R}^{n_1\times \dots\times n_{k-1}\times m_k\times n_{k+1}\times \dots \times n_{d}}$, where its entries are given by
$$c_{i_{1}\dots i_{k-1}ji_{k+1}\dots i_{d}}=\sum_{i_{k}=1}^{n_k}a_{i_{1}\dots i_{k-1}i_{k}i_{k+1}\dots i_{d}}b_{ji_{k}}.$$

For two tensors $\mathcal{A},\mathcal{B}\in \mathbb{R}^{n_1\times n_2\times \dots \times n_d}$, the {\it Frobenius norm} of tensor $\mathcal{A}$ is given by {\color{red}$\|\mathcal{A}\|=\sqrt{\langle\mathcal{A},\mathcal{A}\rangle}$} and the scalar product $\langle\mathcal{A},\mathcal{B}\rangle$ is defined as
$$\langle\mathcal{A},\mathcal{B}\rangle=\sum_{i_1,i_2,\dots,i_d=1}^{n_1,n_2,\dots,n_d}a_{i_{1}i_{2}\dots i_{d}}b_{i_{1}i_{2}\dots i_{d}}.$$

The Kronecker product (see \cite[Section 1.3.6]{Golub2013matrix}) of two matrices $\mathbf{A}\in\mathbb{R}^{m\times p}$ and $\mathbf{B}\in \mathbb{R}^{n\times q}$, denoted by  $\mathbf{A}\otimes \mathbf{B}\in\mathbb{R}^{mn\times pq}$, is given by
\begin{equation*}
    \mathbf{A}\otimes \mathbf{B}=
    \begin{pmatrix}
    a_{11}\mathbf{B}& a_{12}\mathbf{B}& \dots & a_{1p}\mathbf{B}\\
    a_{21}\mathbf{B}& a_{22}\mathbf{B}& \dots & a_{2p}\mathbf{B}\\
    \vdots & \vdots & \ddots & \vdots\\
    a_{m1}\mathbf{B}& a_{m2}\mathbf{B}& \dots & a_{mp}\mathbf{B}\\
    \end{pmatrix}.
\end{equation*}
{\color{red}The Khatri-Rao product (see \cite[Section 12.3.3]{Golub2013matrix}) of $\mathbf{A}\in\mathbb{R}^{m\times k}$ and $\mathbf{B}\in\mathbb{R}^{n\times k}$ is given by $\mathbf{A}\odot\mathbf{B}=[\mathbf{a}_1\otimes \mathbf{b}_1,\cdots,\mathbf{a}_k\otimes \mathbf{b}_k]$, where $\mathbf{a}_i$ and $\mathbf{b}_i$ are the $i$-th columns of $\mathbf{A}$ and $\mathbf{B}$, respectively.}

Notations $C_{{\rm mm}}$ and $C_{{\rm svd}}$ from \cite{martinsson2016randomized,yu2018efficient} are used to count the computational complexity of all algorithms in this paper: $C_{{\rm mm}}$ and $C_{{\rm svd}}$ represent constants so that the cost of multiplying two dense matrices of size $n_1\times n_2$ and $n_2\times n_3$ is $C_{{\rm mm}}n_1n_2n_3$, and the cost of computing the best rank-$r$ approximation of an $n_1\times n_2$ matrix with $r\leq \min\{n_1,n_2\}$ is $C_{{\rm svd}}n_1n_2r$. For a given tolerance $0<\delta<1$, the $\delta$-rank of any given matrix $\mathbf{A}\in\mathbb{R}^{n_1\times n_2}$, denoted by ${\rm rank}_{\delta}(\mathbf{A})$, is defined as the minimum rank of matrix $\mathbf{B}$ over all matrices $\mathbf{B}\in\mathbb{R}^{n_1\times n_2}$ satisfying $\|\mathbf{A}-\mathbf{B}\|\leq \delta$.

Any matrix $\mathbf{\Omega}\in\mathbb{R}^{m\times k}$ is a standard Gaussian matrix (cf. \cite{tropp2017practical}) if its entries form an independent family of standard normal random variables (i.e., Gaussian distribution with mean zero and variance one). For a matrix $\mathbf{A}\in\mathbb{R}^{m\times n}$ with $m\geq n$, the matrix $\mathbf{Q}={\rm orth}(\mathbf{A})$ is an orthonormal matrix whose columns form a basis for the range of $\mathbf{A}$. In practice, this matrix is typically achieved most efficiently by calling a packaged QR factorization (e.g., in MATLAB, we write $\mathbf{Q}={\rm qr}(\mathbf{A},0)$). 
\subsection{Overview of Tucker decomposition}
As shown in \cite{tucker1966some}, the Tucker decomposition is a form of high-order principal component analysis that factorizes a tensor into a core tensor and a set of factor matrices (one for each mode). In detail, for a given $d$-tuple $(r_1,r_2,\dots,r_d)$, Tucker decomposition of $\mathcal{A}\in\mathbb{R}^{n_1\times n_2\times \dots\times n_d}$ is given as
\begin{equation}
\label{dash-rthosvd:tucker-form}
    \mathcal{A}=\mathcal{G}\times_1\mathbf{U}_1\times_2\mathbf{U}_2\dots\times_d\mathbf{U}_d
    +\mathcal{E},
\end{equation}
where $\mathcal{E}\in\mathbb{R}^{n_1\times n_2\times \dots\times n_d}$ is any noise tensor, $\mathcal{G}\in\mathbb{R}^{r_1\times r_2\times \dots\times r_d}$ is the core tensor, and {\color{red}$\mathbf{U}_k\in\mathbb{R}^{n_k\times r_k}$ is the mode-$k$ factor matrix} with $k=1,2,\dots,n$. The Tucker-rank of $\mathcal{A}$ is denoted by the $d$-tuple $${\color{red}({\rm rank}(\mathbf{A}_{(1)}),{\rm rank}(\mathbf{A}_{(2)}),\dots,{\rm rank}(\mathbf{A}_{(d)})),}$$ where {\color{red}$\mathbf{A}_{(k)}\in\mathbb{R}^{n_k\times n_1\cdots n_{k-1}n_{k+1}\cdots n_d}$} is the mode-$k$ unfolding matrix of $\mathcal{A}$. When $\mathcal{E}$ is the zero tensor, it follows from (\ref{dash-rthosvd:tucker-form}) that
\begin{equation*}
    \mathbf{A}_{(k)}=\mathbf{U}_k\mathbf{G}_{(k)}(\mathbf{U}_1\otimes\dots\otimes\mathbf{U}_{k-1}\otimes\mathbf{U}_{k+1}\otimes\dots\otimes\mathbf{U}_d),
\end{equation*}
which implies that ${\rm rank}(\mathbf{A}_{(k)})\leq r_k$.

\subsection{ The proposed problem and two common algorithms}
We now introduce a general form for the fixed-rank approximation in Tucker format, which is given in the following problem.
\begin{problem}
	Suppose that $\mathcal{A}\in\mathbb{R}^{n_1\times n_2\times \dots \times n_d}$. For a given $d$-tuple of positive integers $(r_1,r_2,\dots,r_d)$ with $r_k<n_k$, the goal is to find $d$ orthonormal matrices ${\bf U}_{k}\in\mathbb{R}^{n_{k}\times r_{k}}$ such that
	\begin{equation*}
        {\color{red}\mathcal{A}\approx\mathcal{A}\times_1(\mathbf{U}_{1}\mathbf{U}_{1}^\top)\times_2(\mathbf{U}_{2}\mathbf{U}_{2}^\top)\cdots\times_d(\mathbf{U}_{d}\mathbf{U}_{d}^\top)}.
	\end{equation*}
    \label{dash-rthosvd:prob1}
\end{problem}

Suppose that $\{\mathbf{U}_1,\mathbf{U}_2,\dots,\mathbf{U}_d\}$ is any solution for Problem \ref{dash-rthosvd:prob1} with a given $d$-tuple $(r_1,r_2,\dots,r_d)$. We call the tensor $\widehat{\mathcal{A}}=\mathcal{A}\times_1(\mathbf{U}_1\mathbf{U}_1^\top)\times_2(\mathbf{U}_2\mathbf{U}_2^\top)\dots\times_d(\mathbf{U}_d\mathbf{U}_d^\top)$ as an approximation for the Tucker decomposition of $\mathcal{A}$ with the $d$-tuple $(r_1,r_2,\dots,r_d)$. 

The truncated high-order singular value decomposition (T-HOSVD) (see \cite[Fig. 4.3]{kolda2009tensor}) and the sequentially T-HOSVD (ST-HOSVD) (see \cite[Algorithm 1]{vannieuwenhoven2012new}) are two common algorithms for solving Problem \ref{dash-rthosvd:prob1}. The detailed processes for T-HOSVD and ST-HOSVD are summarized in Algorithm \ref{dash-rthosvd:alg:sthosvd}. For clarity, for the case of ST-HOSVD, the processing order $(1,2,\dots,d)$ can be replaced by any other element in the permutation group $\mathbb{S}_d$ generated by $(1,2,\dots,d)$.

\begin{algorithm}[htb]
     \caption{T-HOSVD and ST-HOSVD with a fixed Tucker-rank}
     \begin{algorithmic}[1]
        \STATEx {\bf Input}: A tensor $\mathcal{A}\in\mathbb{R}^{n_1\times n_2\times\dots\times n_d}$, and the desired $d$-tuple $\mathbf{r}=(r_1,r_2,\dots,r_d)$ with $r_k< n_k$ and $k=1,2,\dots,d$.
        \STATEx {\bf Output}: The core tensor $\mathcal{G}\in\mathbb{R}^{r_1\times r_2\times\dots\times r_d}$ and the mode-$k$ factor matrix $\mathbf{U}_k\in\mathbb{R}^{n_k\times r_k}$ such that $\mathcal{A}\approx\widehat{\mathcal{A}}:=\mathcal{G}\times_1\mathbf{U}_1\times_2\mathbf{U}_2\dots\times_d\mathbf{U}_d$.
        \STATE Set a temporary tensor $\mathcal{G}:=\mathcal{A}$.
        \FOR{$k=1,2,\dots,d$}
           \IF{in T-HOSVD case}
              \STATE Compute $[\mathbf{Q}_k,\sim,\sim]={\rm svd}(\mathbf{A}_{(k)},'{\rm econ}')$.
            \ELSIF{in ST-HOSVD case}
               \STATE Compute $[\mathbf{Q}_k,\sim,\sim]={\rm svd}(\mathbf{G}_{(k)},'{\rm econ}')$.
            \ENDIF
            \STATE Form $\mathbf{U}_k=\mathbf{Q}_k(:,1:r_k)$.
            \STATE Update $\mathcal{G}=\mathcal{G}\times_k\mathbf{U}_k^\top$.
        \ENDFOR
    \end{algorithmic}
    \label{dash-rthosvd:alg:sthosvd}
\end{algorithm}

From Algorithm \ref{dash-rthosvd:alg:sthosvd}, we see that $\mathcal{G}=\mathcal{A}\times_1\mathbf{U}_1^\top\times_2\mathbf{U}_2^\top\dots\times_d\mathbf{U}_d^\top$, which implies that $\widehat{\mathcal{A}}=\mathcal{A}\times_1(\mathbf{U}_1\mathbf{U}_1^\top)\times_2(\mathbf{U}_2\mathbf{U}_2^\top)\dots\times_d(\mathbf{U}_d\mathbf{U}_d^\top)$. Hence, according to \cite[Theorems 5.1 and 6.5]{vannieuwenhoven2012new}, in terms of the Frobenius norm of a tensor, the error in approximating $\mathcal{A}$ using Algorithm \ref{dash-rthosvd:alg:sthosvd} is shown in the following theorem.
{\color{red}
\begin{theorem}
    For a given multilinear rank $\mathbf{r}=(r_1,r_2,\dots,r_d)$ with $r_k< n_k$ and $k=1,2,\dots,d$, let $\widehat{\mathcal{A}}:=\mathcal{G}\times_1\mathbf{U}_1\times_2\mathbf{U}_2\dots\times_d\mathbf{U}_d$ be obtained by applying Algorithm \ref{dash-rthosvd:alg:sthosvd} to $\mathcal{A}\in\mathbb{R}^{n_1\times n_2\times \dots\times n_d}$. Then, for both the T-HOSVD and ST-HOSVD cases, one has
    \begin{equation*}
        \|\mathcal{A}-\widehat{\mathcal{A}}\|^2\leq\sum_{k=1}^d\sum_{i=r_{k}+1}^{\widehat{n}_k}\sigma_{i}(\mathbf{A}_{(k)})^2,
    \end{equation*}
    where for each $k$, we introduce $\widehat{n}_k=\min\{n_k,n_1\dots n_{k-1}n_{k+1}\dots n_d\}$, and $\sigma_{i}(\mathbf{A}_{(k)})$ is the $i$-th singular value of $\mathbf{A}_{(k)}$.
    \label{dash-rthosvd:general:errror}
\end{theorem}
}
We introduce the computational complexity of Algorithm \ref{dash-rthosvd:alg:sthosvd}. In detail, T-HOVSD and ST-HOSVD require
\begin{equation*}
    \begin{cases}
        \begin{aligned}
            &T_{{\rm T}}(n_1,\dots,n_d;r_1,\dots,r_d)
            =\sum_{k=1}^dC_{{\rm mm}}r_{1}r_2\dots r_kn_kn_{k+1}\dots n_d\\
            &\quad+\sum_{k=1}^d\left( O(n_1n_2\dots n_d)+C_{{\rm svd}}n_1n_2\dots n_dr_k\right),\\
            &T_{{\rm ST}}(n_1,\dots,n_d;r_1,\dots,r_d)
            =\sum_{k=1}^dC_{{\rm mm}}r_{1}r_2\dots r_kn_kn_{k+1}\dots n_d\\
            &\quad+\sum_{k=1}^d\left( O(r_1\dots r_{k-1}n_k\dots n_d)+C_{{\rm svd}}r_1\dots r_{k-1}n_{k}\dots n_dr_k\right)
        \end{aligned}
    \end{cases}
\end{equation*}
operations, respectively, to obtain the approximate Tucker decomposition $\widehat{\mathcal{A}}$ of $\mathcal{A}$ with a given Tucker-rank $\mathbf{r}=(r_1,r_2,\dots,r_d)$. Here, for each $k$, reshaping the tensor $\mathcal{A}$ to the mode-$k$ unfolding $\mathbf{A}_{(k)}$ requires $O(n_1n_2\dots n_d)$ operations.

\begin{remark}
    As shown in Algorithm \ref{dash-rthosvd:alg:sthosvd}, all the mode-$k$ factor matrices $\mathbf{U}_k$ are orthonormal. Meanwhile, there exist several types for Tucker decomposition, such as, nonnegative Tucker decomposition (cf. \cite{zhou2012fast}), high-order interpolatory decomposition/tensor CUR decomposition (cf. \cite{cai2021mode,che2022perturbations,drineas2007a,saibaba2016hoid}), and structure-preserving decomposition in Tucker format (cf. \cite{minster2020randomized}).
\end{remark}
\section{The proposed algorithms}
\label{dash-rthosvd:sec3}

{\color{red}By combining the random projection technique and shifted simultaneous iteration within Algorithm \ref{dash-rthosvd:alg:sthosvd}, we now present efficient algorithms for solving Problem \ref{dash-rthosvd:prob1} in the T-HOSVD and ST-HOSVD formats. We count the computational complexity of these algorithms and give a detailed analysis of the approximation of Tucker decomposition derived from these algorithms.}

For the purpose of further discussion, {\color{red}we transform Problem \ref{dash-rthosvd:prob1} into the following optimization problem (Problem \ref{dash-rthosvd:prob3}) by using the Frobenius norm of a tensor}.
\begin{problem}
	Suppose that $\mathcal{A}\in\mathbb{R}^{n_1\times n_2\times \dots \times n_d}$. For a given $d$-tuple of positive integers $(r_1,r_2,\dots,r_d)$ with $r_k<n_k$, the goal is to find $d$ orthonormal matrices ${\bf U}_{k}\in\mathbb{R}^{n_{k}\times r_{k}}$ such that the minimum value of the function $f(\mathbf{V}_1,\mathbf{V}_2,\dots,\mathbf{V}_d)$ is attained at $\{\mathbf{U}_1,\mathbf{U}_2,\dots,\mathbf{U}_d\}$, with
	\begin{equation*}
        f(\mathbf{V}_1,\dots,\mathbf{V}_d)
        =\|\mathcal{A}\times_1(\mathbf{V}_1\mathbf{V}_1^\top)
        \times_2(\mathbf{V}_2\mathbf{V}_2^\top)\dots
        \times_d(\mathbf{V}_d\mathbf{V}_d^\top)-\mathcal{A}\|
	\end{equation*}
	where for each $k$, $\mathbf{V}_k\in\mathbb{R}^{n_k\times r_k}$ is orthonormal.
\label{dash-rthosvd:prob3}
\end{problem}

As shown in \cite{saibaba2016hoid,vannieuwenhoven2012new}, we have
\begin{align}
\label{dash-rthosvd:thosvd-general}
    \left\|\mathcal{A}\times_1\left({\bf V}_1{\bf V}_1^\top\right)\dots\times_d\left({\bf V}_d{\bf V}_d^\top\right)-\mathcal{A}\right\|^2\leq\sum_{k=1}^d\left\|\mathcal{A}\times_k\left({\bf V}_k{\bf V}_k^\top\right)-\mathcal{A}\right\|^2,
\end{align}
and
\begin{align}
\label{dash-rthosvd:sthosvd-general}
    &\left\|\mathcal{A}\times_1\left({\bf V}_1{\bf V}_1^\top\right)\times_2\left({\bf V}_2{\bf V}_2^\top\right)\dots\times_d\left({\bf V}_d{\bf V}_d^\top\right)-\mathcal{A}\right\|^2\nonumber\\
    &\leq\left\|\mathcal{A}\times_{p_1}\left(\mathbf{I}_{n_{p_1}}-{\bf V}_{p_1}{\bf V}_{p_1}^\top\right)\right\|^2
    +\left\|\left(\mathcal{A}\times_{p_1}{\bf V}_{p_1}^\top\right)\times_{p_2}\left(\mathbf{I}_{n_{p_2}}-{\bf V}_{p_2}{\bf V}_{p_2}^\top\right)\right\|^2\nonumber\\
    &+\dots+\left\|\left(\mathcal{A}\times_{p_1}{\bf V}_{p_1}^\top\times_{p_2}\mathbf{V}_{p_2}^\top\dots\times_{p_{d-1}}\mathbf{V}_{p_{d-1}}^\top\right)\times_{p_d}\left(\mathbf{I}_{n_{p_d}}-{\bf V}_{p_d}{\bf V}_{p_d}^\top\right)\right\|^2,
\end{align}
{\color{red}with $(p_1,p_2,\dots,p_d)\in\mathbb{S}_d$}.
Hence, a quasi-optimal solution $(\mathbf{U}_1,\mathbf{U}_2,\cdots,\mathbf{U}_d)$ to Problem \ref{dash-rthosvd:prob3} can be obtained from the following two problems.
\begin{problem}
\label{dash-rthosvd:prob4}
	Suppose that $\mathcal{A}\in\mathbb{R}^{n_1\times n_2\times \dots \times n_d}$. For a given positive integer $r_k$ with $r_k<n_k$, the goal is to find an orthonormal matrix ${\bf U}_{k}\in\mathbb{R}^{n_{k}\times r_{k}}$ such that the minimum value of the function $f_k(\mathbf{V}_k)$ is attained at $\mathbf{U}_k$, with
	\begin{equation*}
        f_k(\mathbf{V}_k)
        =\|\mathcal{A}\times_k(\mathbf{V}_k\mathbf{V}_k^\top)-\mathcal{A}\|^2.
	\end{equation*}
\end{problem}
\begin{problem}
\label{dash-rthosvd:prob5}
	Suppose that $\mathcal{A}\in\mathbb{R}^{n_1\times n_2\times \dots \times n_d}$. Let $(p_1,p_2,\dots,p_d)\in\mathbb{S}_d$. For a given $k$, let $\mathbf{U}_{p_l}\in\mathbb{R}^{n_{p_l}\times r_{p_l}}$ be orthonormal with $l=1,2,\dots,k-1$. For a given positive integer $r_{p_k}$ with $r_{p_k}<n_{p_k}$, the goal is to find an orthonormal matrix ${\bf U}_{p_k}\in\mathbb{R}^{n_{p_k}\times r_{p_k}}$ such that the minimum value of the function $g_k(\mathbf{U}_{p_1},\dots,\mathbf{U}_{p_{k-1}},\mathbf{V}_{p_k})$ is attained at $\mathbf{U}_{p_k}$, with
	\begin{equation*}
        g_k(\mathbf{U}_{p_1},\dots,\mathbf{U}_{p_{k-1}},\mathbf{V}_{p_k})=\left\|\left(\mathcal{A}\times_{p_1}{\bf U}_{p_1}^\top\dots\times_{p_{k-1}}\mathbf{U}_{p_{k-1}}^\top\right)\times_{p_k}\left(\mathbf{I}_{n_{p_k}}-{\bf V}_{p_k}{\bf V}_{p_k}^\top\right)\right\|^2.
	\end{equation*}
\end{problem}
   \begin{remark} 
    Note that in the fixed Tucker-rank case, the framework for T-HOSVD is based on Problem \ref{dash-rthosvd:prob4} and the framework for ST-HOSVD is based on Problem \ref{dash-rthosvd:prob5}.
\end{remark}

\begin{remark}
    All solutions $\mathbf{U}_k$ with $k=1,2,\cdots,d$ to Problem \ref{dash-rthosvd:prob4} or \ref{dash-rthosvd:prob5} form a quasi-optimal solution $(\mathbf{U}_1,\mathbf{U}_2,\cdots,\mathbf{U}_d)$ to Problem \ref{dash-rthosvd:prob3}. Similar to (\ref{dash-rthosvd:thosvd-general}) and (\ref{dash-rthosvd:sthosvd-general}), one has
    \begin{align*}
    \left\|\mathcal{A}\times_1\left({\bf U}_1{\bf U}_1^\top\right)\cdots\times_d\left({\bf U}_d{\bf U}_d^\top\right)-\mathcal{A}\right\|^2\leq\sum_{k=1}^d\left\|\mathcal{A}\times_k\left({\bf U}_k{\bf U}_k^\top\right)-\mathcal{A}\right\|^2,
\end{align*}
and
\begin{align*}
    &\left\|\mathcal{A}\times_1\left({\bf U}_1{\bf U}_1^\top\right)\cdots\times_d\left({\bf U}_d{\bf U}_d^\top\right)-\mathcal{A}\right\|^2\nonumber\\
    &\leq\left\|\mathcal{A}\times_{1}\left(\mathbf{I}_{n_{1}}-{\bf U}_{1}{\bf U}_{1}^\top\right)\right\|^2
    +\left\|\left(\mathcal{A}\times_{1}{\bf U}_{1}^\top\right)\times_{2}\left(\mathbf{I}_{n_{2}}-{\bf U}_{2}{\bf U}_{2}^\top\right)\right\|^2\nonumber\\
    &+\dots+\left\|\left(\mathcal{A}\times_{1}{\bf U}_{1}^\top\times_{2}\mathbf{U}_{2}^\top\dots\times_{{d-1}}\mathbf{U}_{{d-1}}^\top\right)\times_{d}\left(\mathbf{I}_{n_{d}}-{\bf U}_{d}{\bf U}_{d}^\top\right)\right\|^2
\end{align*}
with $(p_1,p_2,\cdots,p_d)=(1,2,\cdots,d)$. Theorems \ref{dash-rthosvd:thm1} and \ref{dash-rthosvd:thm2} were deduced, respectively, according to the above two inequalities. 

However, there exists an unresolved issue: does there exist a constant $C>0$ such that
    \begin{align*}
    \sum_{k=1}^d\left\|\mathcal{A}\times_k\left({\bf U}_k{\bf U}_k^\top\right)-\mathcal{A}\right\|^2\leq C\cdot \left\|\mathcal{A}\times_1\left({\bf U}_1{\bf U}_1^\top\right)\cdots\times_d\left({\bf U}_d{\bf U}_d^\top\right)-\mathcal{A}\right\|^2,
\end{align*}
and
\begin{align*}
    &\left\|\mathcal{A}\times_{1}\left(\mathbf{I}_{n_{1}}-{\bf U}_{1}{\bf U}_{1}^\top\right)\right\|^2
    +\left\|\left(\mathcal{A}\times_{1}{\bf U}_{1}^\top\right)\times_{2}\left(\mathbf{I}_{n_{2}}-{\bf U}_{2}{\bf U}_{2}^\top\right)\right\|^2\\
    &+\dots+\left\|\left(\mathcal{A}\times_{1}{\bf U}_{1}^\top\times_{2}\mathbf{U}_{2}^\top\dots\times_{{d-1}}\mathbf{U}_{{d-1}}^\top\right)\times_{d}\left(\mathbf{I}_{n_{d}}-{\bf U}_{d}{\bf U}_{d}^\top\right)\right\|^2\\
    &\leq C\cdot \left\|\mathcal{A}\times_1\left({\bf U}_1{\bf U}_1^\top\right)\cdots\times_d\left({\bf U}_d{\bf U}_d^\top\right)-\mathcal{A}\right\|^2?
\end{align*}
\end{remark}

\subsection{Shifted versions for randomized T-HOSVD and ST-HOSVD}

Following Algorithm \ref{dash-rthosvd:alg1}, {\color{red}for each $k$, the computation $[\mathbf{Q}_k,\sim,\sim]={\rm svd}(\mathbf{A}_{(k)}\mathbf{A}_{(k)}^\top\mathbf{Q}_k,'{\rm econ}')$ with $t=1,2,\cdots,q$ matches the simultaneous iteration used to compute dominant subspaces of $\mathbf{A}_{(k)}\mathbf{A}_{(k)}^\top$ (cf. \cite[Section 8.2.4]{Golub2013matrix})}. As shown in \cite[Section 7.5]{Golub2013matrix}, {\color{red}the shifted method} can be used to accelerate the convergence of the power method because of the reduction of the ratio between the second largest and the largest eigenvalues, and proper selection of the shift is crucial to maintain the dominance of the target eigenvalue and avoid complications. Feng {\it et al.} \cite{feng2024algorithm} applied a dynamically shifted power iteration technique to improve the accuracy of the randomized SVD method (e.g., \cite{halko2011finding}).

\begin{algorithm}
     \caption{A shifted version for the randomized variant of T-HOSVD with a fixed Tucker-rank}
     \begin{algorithmic}[1]
        \STATEx {\bf Input}: A tensor $\mathcal{A}\in\mathbb{R}^{n_1\times n_2\times\dots\times n_d}$, the desired $d$-tuple $\mathbf{r}=(r_1,r_2,\dots,r_d)$ with $r_k<n_k$ and $k=1,2,\dots,d$, the $d$-tuple of oversampling parameters $\mathbf{s}=(s_1,s_2,\dots,s_d)$, and the power parameter $q\geq 1$.
        \STATEx {\bf Output}: The core tensor $\mathcal{G}\in\mathbb{R}^{r_1\times r_2\times\dots\times r_d}$ and the mode-$k$ factor matrix $\mathbf{U}_k\in\mathbb{R}^{n_k\times r_k}$ such that $\mathcal{A}\approx\widehat{\mathcal{A}}:=\mathcal{G}\times_1\mathbf{U}_1\times_2\mathbf{U}_2\dots\times_d\mathbf{U}_d$.
        \STATE Let $l_k=r_k+s_k$ with $k=1,2,\dots,d$ and set a temporary tensor $\mathcal{G}:=\mathcal{A}$.
        \FOR{$k=1,2,\dots,d$}
            \STATE Draw a standard Gaussian matrix $\mathbf{\Omega}_k\in\mathbb{R}^{n_1\dots n_{k-1}n_{k+1}\dots n_{d}\times l_k}$.
            \STATE Form the mode-$k$ unfolding $\mathbf{A}_{(k)}$ from $\mathcal{A}$.
            \STATE Compute $[\mathbf{Q}_k,\sim,\sim]={\rm svd}(\mathbf{A}_{(k)}\mathbf{\Omega}_k,'{\rm econ}')$ and set {\color{red}$\alpha_{0}^{(k)}=0$}.
            \FOR{$t=1,2,\dots,q$}
                \STATE Compute $[\mathbf{Q}_k,\bm{\Sigma}_k,\sim]={\rm svd}(\mathbf{A}_{(k)}(\mathbf{A}_{(k)}^\top\mathbf{Q}_k)-{\color{red}\alpha_{t-1}^{(k)}} \mathbf{Q}_k,'{\rm econ}')$.
                \IF{$\bm{\Sigma}_k(l_k,l_k)>{\color{red}\alpha_{t-1}^{(k)}}$}
                    \STATE Update {\color{red}$\alpha_{t}^{(k)}=(\mathbf{\Sigma}_k(l_k,l_k)+\alpha_{t-1}^{(k)})/2$}.
                \ENDIF
            \ENDFOR
            \STATE Set $\mathbf{U}_k=\mathbf{Q}_k(:,1:r_k)$ and update $\mathcal{G}=\mathcal{G}\times_k\mathbf{U}_{k}^\top$.
        \ENDFOR
    \end{algorithmic}
    \label{dash-rthosvd:alg2}
\end{algorithm}

For each $k$, let $l_k=r_k+s_k$ satisfy $l_k<\min\{n_k, n_1\dots n_{k-1}n_{k+1}\dots n_d\}$ and ${\rm rank}(\mathbf{A}_{(n)})\geq l_k$, which implies that each element in the Tucker-rank of $\mathcal{A}$ is larger than the corresponding one in $(r_1,r_2,\dots,r_d)$. By combining the shifted method with Algorithm \ref{dash-rthosvd:alg1}, we obtain a shifted version for randomized T-HOSVD with a fixed Tucker-rank, which is summarized as Algorithm \ref{dash-rthosvd:alg2} and can be viewed as the generalization of Algorithm 2 in \cite{feng2024algorithm} from matrices to tensors. We now discuss several key points in this algorithm.
\begin{enumerate}
    \item[(a)] Similar to the description in \cite{feng2024algorithm}, if $0<\alpha<\sigma_{l_k}(\mathbf{A}_{(k)}\mathbf{A}_{(k)}^\top)/2$, then for each $i_k\leq l_k$, $\sigma_{i_k}(\mathbf{A}_{(k)}\mathbf{A}_{(k)}^\top-\alpha \mathbf{I}_{n_k})=\sigma_{i_k}(\mathbf{A}_{(k)}\mathbf{A}_{(k)}^\top)-\alpha$. If $\sigma_{i_k}(\mathbf{A}_{(k)}\mathbf{A}_{(k)}^\top-\alpha \mathbf{I}_{n_k})\neq \sigma_{l_k+1}(\mathbf{A}_{(k)}\mathbf{A}_{(k)}^\top-\alpha \mathbf{I}_{n_k})$, then the left singular vector corresponding to the $i_k$-th singular value of $\mathbf{A}_{(k)}\mathbf{A}_{(k)}^\top-\alpha \mathbf{I}_{n_k}$ is that corresponding to the $i_k$-th singular value of $\mathbf{A}_{(k)}\mathbf{A}_{(k)}^\top$. Following the above descriptions, for each $k$, the computation $[\mathbf{Q}_k,\sim,\sim]={\rm svd}(\mathbf{A}_{(k)}\mathbf{A}_{(k)}^\top\mathbf{Q}_k,'{\rm econ}')$ in Algorithm \ref{dash-rthosvd:alg1} is changed to $[\mathbf{Q}_k,\sim,\sim]={\rm svd}(\mathbf{A}_{(k)}\mathbf{A}_{(k)}^\top\mathbf{Q}_k-\alpha\mathbf{Q}_k,'{\rm econ}')$ with $0<\alpha<\sigma_{l_k}(\mathbf{A}_{(k)}\mathbf{A}_{(k)}^\top)/2$, which can be used to obtain the same approximated dominant subspace.
    \item[(b)] Under the case of $0<\alpha<\sigma_{l_k}(\mathbf{A}_{(k)}\mathbf{A}_{(k)}^\top)/2$, for $j_k<i_k\leq l_k$, suppose that $\sigma_{j_k}(\mathbf{A}_{(k)}\mathbf{A}_{(k)}^\top)>\sigma_{i_k}(\mathbf{A}_{(k)}\mathbf{A}_{(k)}^\top)$, then one has
    \begin{align}
    \label{dash-rthosvd:appeqn1}
        f(\alpha):=\frac{\sigma_{i_k}(\mathbf{A}_{(k)}\mathbf{A}_{(k)}^\top-\alpha\mathbf{I}_{n_k})}{\sigma_{j_k}(\mathbf{A}_{(k)}\mathbf{A}_{(k)}^\top-\alpha\mathbf{I}_{n_k})}=\frac{\sigma_{i_k}(\mathbf{A}_{(k)}\mathbf{A}_{(k)}^\top)-\alpha}{\sigma_{j_k}(\mathbf{A}_{(k)}\mathbf{A}_{(k)}^\top)-\alpha},
    \end{align}
   which is strictly monotonically decreasing with respect to $\alpha$, implying that the singular values of $\mathbf{A}_{(k)}\mathbf{A}_{(k)}^\top-\alpha\mathbf{I}_{n_k}$ decay faster than those of $\mathbf{A}_{(k)}\mathbf{A}_{(k)}^\top$ and the larger the value of $\alpha$ is, the smaller the ratio $f(\alpha)$ becomes, reflecting faster decay of singular values. This makes the dominant subspace of $\mathbf{A}_{(k)}\mathbf{A}_{(k)}^\top-\alpha\mathbf{I}_{n_k}$ more distinguishable, accelerating power iteration. Hence, the shifted method avoids slow convergence for slowly decaying spectra.To maximize the effect of the shifted power iteration on improving the accuracy, we should choose $\alpha$ as large as possible while satisfying $0<\alpha<\sigma_{l_k}(\mathbf{A}_{(k)}\mathbf{A}_{(k)}^\top)/2$.
    \item[(c)] How to set $\alpha$ properly is a key issue. Note that for large $n_k$, it is impossible to directly calculate $\sigma_{l_k}(\mathbf{A}_{(k)}\mathbf{A}_{(k)}^\top)$. One idea is to use the singular value of $\mathbf{A}_{(k)}\mathbf{A}_{(k)}^\top\mathbf{Q}_k$ to approximate $\sigma_{l_k}(\mathbf{A}_{(k)}\mathbf{A}_{(k)}^\top)$. Since $\sigma_{i_k}(\mathbf{A}_{(k)}\mathbf{A}_{(k)}^\top\mathbf{Q}_k)\leq \sigma_{i_k}(\mathbf{A}_{(k)}\mathbf{A}_{(k)}^\top)$ with $i_k=1,2,\dots,l_k$, setting $\alpha=\sigma_{l_k}(\mathbf{A}_{(k)}\mathbf{A}_{(k)}^\top\mathbf{Q}_k)/2$ satisfies $0<\alpha<\sigma_{l_k}(\mathbf{A}_{(k)}\mathbf{A}_{(k)}^\top)/2$. When $0<\alpha<\sigma_{l_k}(\mathbf{A}_{(k)}\mathbf{A}_{(k)}^\top)/2$, we also have
    \begin{align*}
        \sigma_{i_k}((\mathbf{A}_{(k)}\mathbf{A}_{(k)}^\top-\alpha\mathbf{I}_{n_k})\mathbf{Q}_k)+\alpha&\leq \sigma_{i_k}(\mathbf{A}_{(k)}\mathbf{A}_{(k)}^\top-\alpha\mathbf{I}_{n_k})+\alpha\\
        &=\sigma_{i_k}((\mathbf{A}_{(k)}\mathbf{A}_{(k)}^\top)
    \end{align*}
    with $i_k=1,2,\dots,l_k$, which implies how to use the singular values of $(\mathbf{A}_{(k)}\mathbf{A}_{(k)}^\top-\alpha\mathbf{I}_{n_k})\mathbf{Q}_k$ to approximate those of $\mathbf{A}_{(k)}\mathbf{A}_{(k)}^\top$.
\end{enumerate}

{\color{red}We now prove that under certain assumptions, the function $ f(\alpha)$ in (\ref{dash-rthosvd:appeqn1}) is strictly monotonically decreasing with respect to $\alpha$.
\begin{proof}
Without loss of generality, let $a=\sigma_{i_k}(\mathbf{A}_{(k)}\mathbf{A}_{(k)}^\top)$ and $b=\sigma_{j_k}(\mathbf{A}_{(k)}\mathbf{A}_{(k)}^\top)$. It follows from $\sigma_{j_k}(\mathbf{A}_{(k)}\mathbf{A}_{(k)}^\top)>\sigma_{i_k}(\mathbf{A}_{(k)}\mathbf{A}_{(k)}^\top)$ with $j_k<i_k\leq l_k$ that $a<b$. Hence, the function $ f(\alpha)$ can be rewritten as
\begin{equation*}
 f(\alpha)=\frac{a-\alpha}{b-\alpha}
\end{equation*}
with $0<\alpha<\sigma_{l_k}(\mathbf{A}_{(k)}\mathbf{A}_{(k)}^\top)/2$. The differentiation of $ f(\alpha)$ is given by
\begin{equation*}
f'(\alpha)=\frac{-(b-\alpha)+(a-\alpha)}{(b-\alpha)^2}=\frac{a-b}{(b-\alpha)^2}<0,
\end{equation*}
where the inequality holds as $a<b$. Hence, the proof is complete.
\end{proof}
}

Therefore, for each $k$, we obtain an efficient way to update $\alpha_t^{(k)}$ in Algorithm \ref{dash-rthosvd:alg2} to ensure faster decay of singular values. In particular, for the same oversampling vector $\mathbf{s}=(s_1,s_2,\dots,s_d)$, Algorithm \ref{dash-rthosvd:alg2} has similar computational complexity to one version of the randomized variant of T-HOSVD, shown in Algorithm \ref{dash-rthosvd:alg1}; however, the accuracy of Algorithm \ref{dash-rthosvd:alg2} is better than that of Algorithm \ref{dash-rthosvd:alg1}.

We also propose a shifted version for the randomized variant of ST-HOSVD with a fixed Tucker-rank, which is summarized in Algorithm \ref{dash-rthosvd:alg2:v2}. Specifically, for each $k$, a quasi-optimal solution $\mathbf{U}_k$ for Problem \ref{dash-rthosvd:prob5} with $(p_1,p_2,\dots,p_d)=(1,2,\dots,d)$ is obtained from Algorithm \ref{dash-rthosvd:alg2:v2}. For the same oversampling vector $\mathbf{s}=(s_1,s_2,\dots,s_d)$, Algorithm \ref{dash-rthosvd:alg2:v2} has a similar computational complexity to the existing randomized ST-HOSVD, shown in Algorithm \ref{dash-rthosvd:alg1}; however, the accuracy of Algorithm \ref{dash-rthosvd:alg2:v2} is better than that of Algorithm \ref{dash-rthosvd:alg1}.

\begin{algorithm}[htb]
     \caption{A shifted version for the randomized variant of ST-HOSVD with a fixed Tucker-rank}
     \begin{algorithmic}[1]
        \STATEx {\bf Input}: A tensor $\mathcal{A}\in\mathbb{R}^{n_1\times n_2\times\dots\times n_d}$, the desired $d$-tuple $\mathbf{r}=(r_1,r_2,\dots,r_d)$ with $r_k<n_k$ and $k=1,2,\dots,d$, the $d$-tuple of oversampling parameters $\mathbf{s}=(s_1,s_2,\dots,s_d)$, and the power parameter $q\geq 1$.
        \STATEx {\bf Output}: The core tensor $\mathcal{G}\in\mathbb{R}^{r_1\times r_2\times\dots\times r_d}$ and the mode-$k$ factor matrix $\mathbf{U}_k\in\mathbb{R}^{n_k\times r_k}$ such that $\mathcal{A}\approx\widehat{\mathcal{A}}:=\mathcal{G}\times_1\mathbf{U}_1\times_2\mathbf{U}_2\dots\times_d\mathbf{U}_d$.
        \STATE Let $l_k=r_k+s_k$ with $k=1,2,\dots,d$ and set a temporary tensor $\mathcal{G}:=\mathcal{A}$.
        \FOR{$k=1,2,\dots,d$}
            \STATE Draw a standard Gaussian matrix $\mathbf{\Omega}_k\in\mathbb{R}^{r_1\dots r_{k-1}n_{k+1}\dots n_{d}\times l_k}$.
            \STATE Form the mode-$k$ unfolding $\mathbf{G}_{(k)}$ from $\mathcal{G}$.
            \STATE Compute $[\mathbf{Q}_k,\sim,\sim]={\rm svd}(\mathbf{G}_{(k)}\mathbf{\Omega}_k,'{\rm econ}')$ and set {\color{red}$\alpha_0^{(k)}=0$}.
            \FOR{$t=1,2,\dots,q$}
                \STATE Compute $[\mathbf{Q}_k,\bm{\Sigma}_k,\sim]={\rm svd}(\mathbf{G}_{(k)}(\mathbf{G}_{(k)}^\top\mathbf{Q}_k)-\alpha \mathbf{Q}_k,'{\rm econ}')$.
                \IF{$\bm{\Sigma}_k(l_k,l_k)>{\color{red}\alpha_{t-1}^{(k)}}$}
                    \STATE Update {\color{red}$\alpha_{t}^{(k)}=(\mathbf{\Sigma}_k(l_k,l_k)+\alpha_{t-1}^{(k)})/2$}.
                \ENDIF
            \ENDFOR
            \STATE Set $\mathbf{U}_k=\mathbf{Q}_k(:,1:r_k)$ and update $\mathcal{G}=\mathcal{G}\times_k\mathbf{U}_{k}^\top$.
        \ENDFOR
    \end{algorithmic}
    \label{dash-rthosvd:alg2:v2}
\end{algorithm}

\begin{remark}
For each $k$ in Algorithms \ref{dash-rthosvd:alg2} and \ref{dash-rthosvd:alg2:v2}, the standard Gaussian matrix $\mathbf{\Omega}_k$ may be large for high-dimensional tensors, which then require an expensive dense matrix–matrix product with this matrix. We can reduce the runtime by using a structured test matrix $\mathbf{\Omega}_k$ that supplies a faster algorithm for the matrix-matrix product. The most popular examples are sparse random matrices (cf. \cite{achlioptas2001database,clarkson2009numerical,kane2014sparser}), subsampled randomized trigonometric transforms (cf. \cite{ailon2009the,tropp2011improved,woolfe2008a}) and the Khatri-Rao product of random matrices (i.e., random Khatri–Rao products) (cf. \cite{al2023randomized,bamberger2022johnson,camano2025faster,saibaba2026improved}), which reduce computational cost while preserving approximation quality. 
\end{remark}
\begin{remark}
In Section \ref{dash-rthosvd:sec5:subsec5}, we compare the efficiency of Algorithms \ref{dash-rthosvd:alg2} and \ref{dash-rthosvd:alg2:v2} under the standard Gaussian matrix and two types of random Khatri–Rao products. For clarity, the theoretical analyses for Algorithms \ref{dash-rthosvd:alg2} and \ref{dash-rthosvd:alg2:v2} are based on the theory of standard Gaussian matrices, as shown in Section \ref{dash-rthosvd:sec4}. Theoretical analyses based on other types of random matrices can be considered as future work.
\end{remark}

\subsection{Computational complexity}

\begin{figure}[htb]
    \setlength{\tabcolsep}{4pt}
    \renewcommand\arraystretch{1}
    \centering
    \subfigure[Algorithm \ref{dash-rthosvd:alg2}]{\includegraphics[width=4.5in,height=1.6in]{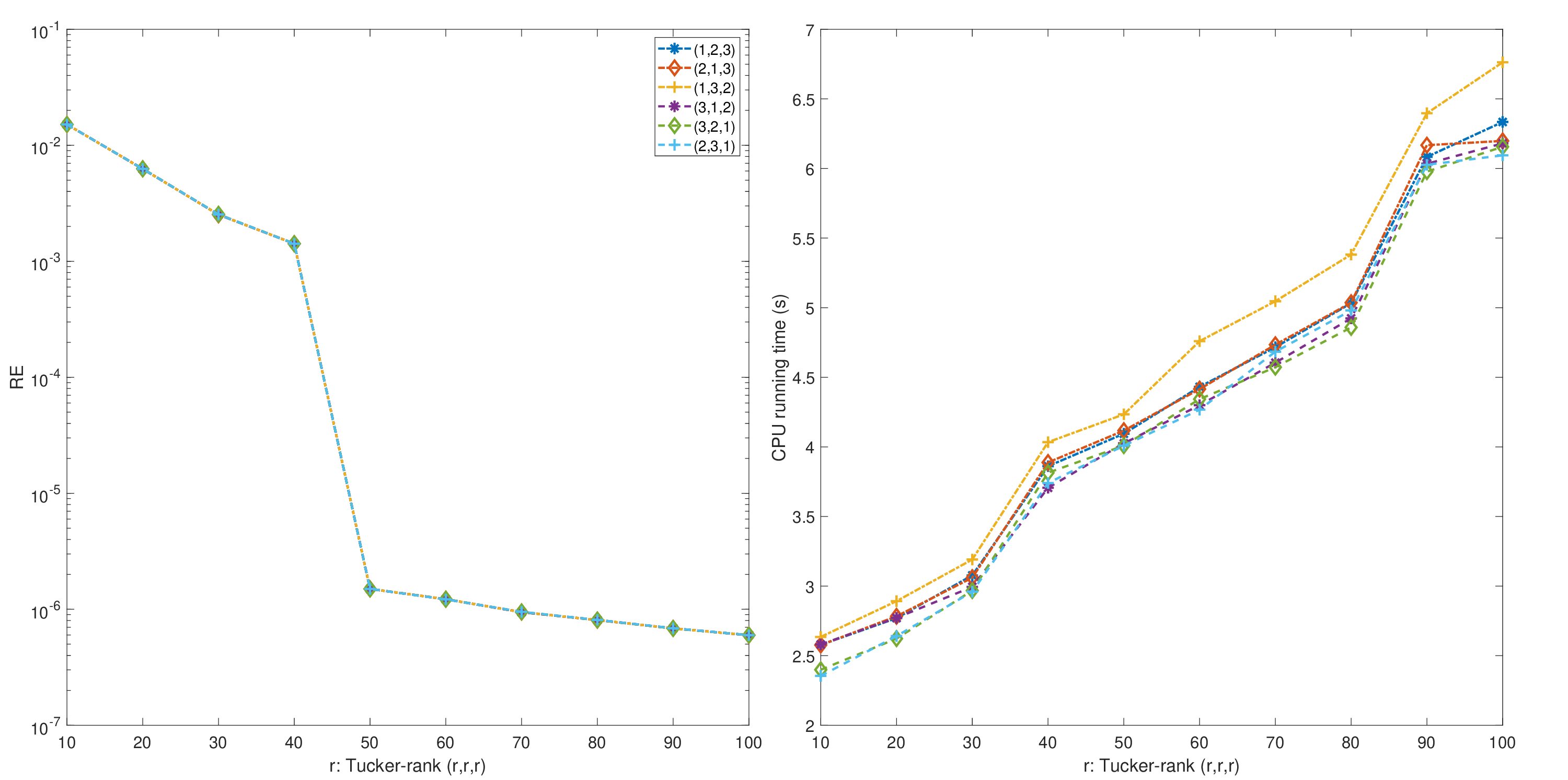}}\\
	\subfigure[Algorithm
\ref{dash-rthosvd:alg2:v2}]{\includegraphics[width=4.5in,height=1.6in]{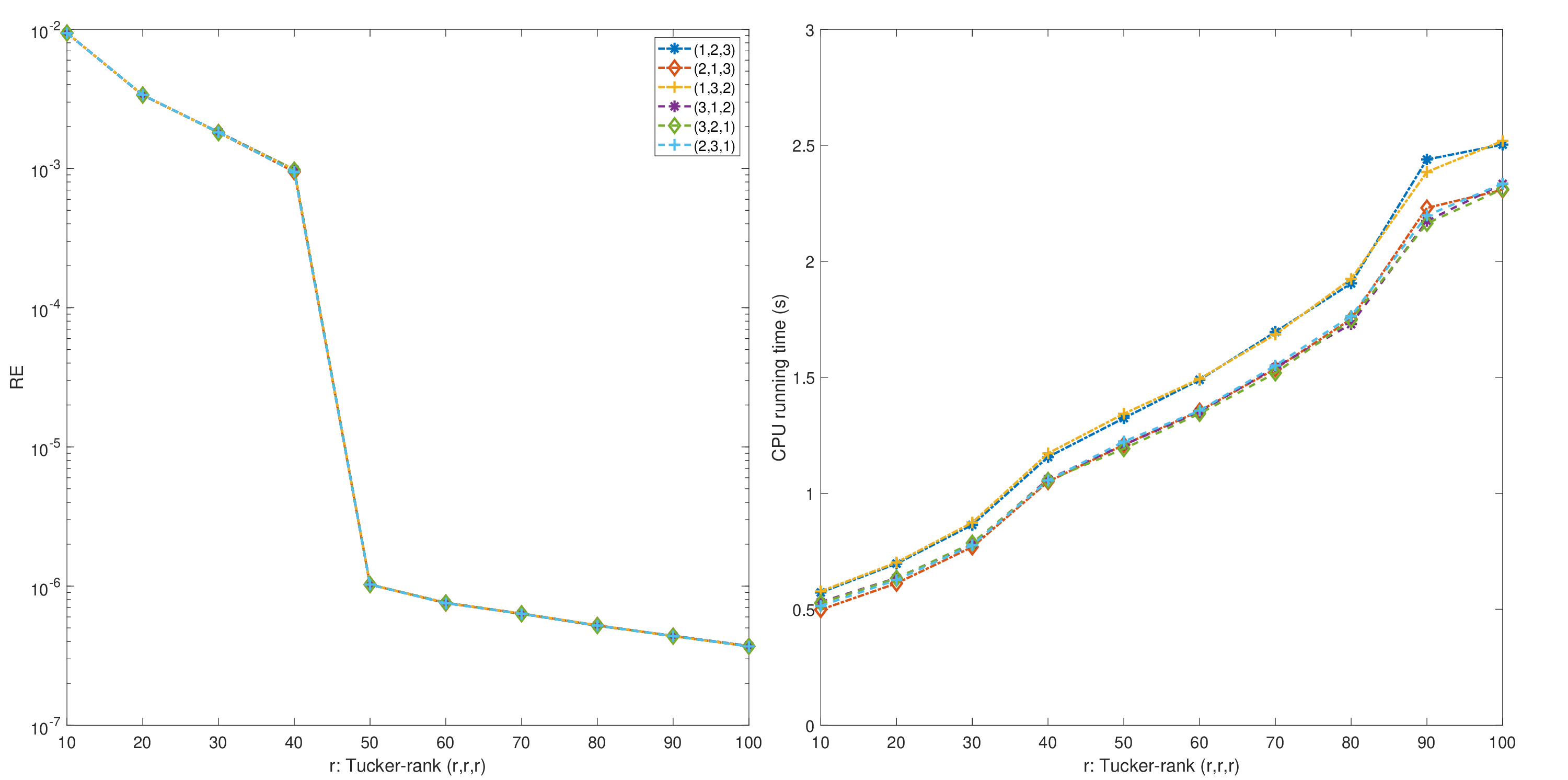}}\\
    \caption{For fixed $s=10$ and $q=1$, numerical simulation results of Algorithms \ref{dash-rthosvd:alg2} and \ref{dash-rthosvd:alg2:v2}
    with different Tucker-ranks $(r,r,r)$ and processing orders on the test tensor $\mathcal{A}$ are shown.}\label{dash-rthosvd:figure3:add}
\end{figure}

Note that the dimensions $(n_1,n_2,\cdots,n_d)$ of the tensors are not all the same. {\color{red}For any $(n_1,n_2,\cdots,n_d)$, there exists a processing order $(p_1,p_2,\dots,p_d)\in\mathbb{S}_d$ such that $n_{p_1}\geq n_{p_2}\geq \cdots\geq n_{p_d}$. When analyzing the computational complexity of Algorithms \ref{dash-rthosvd:alg2} and \ref{dash-rthosvd:alg2:v2}, we assume that $n_1\geq n_2\geq \dots\geq n_d$.  We now illustrate the plausibility of this hypothesis via a test tensor $\mathcal{C}\in\mathbb{R}^{600\times 700\times 800}$:
\begin{equation*}
    \mathcal{A}=\sum_{i=1}^{50}\frac{1000}{i}\mathbf{x}_i\circ\mathbf{y}_i\circ\mathbf{z}_i
    +\sum_{i=51}^{600}\frac{1}{i}\mathbf{x}_i\circ\mathbf{y}_i\circ\mathbf{z}_i
\end{equation*}
where $\mathbf{x}_i\in\mathbb{R}^{600}$, $\mathbf{y}_i\in\mathbb{R}^{700}$ and $\mathbf{z}_i\in\mathbb{R}^{800}$ are four sparse vectors with only 0.05 sparsity. By using Algorithms \ref{dash-rthosvd:alg2} and \ref{dash-rthosvd:alg2:v2} with different Tucker-ranks $(r,r,r)$ and processing orders to the test tensor $\mathcal{C}$, the results are shown in Figure \ref{dash-rthosvd:figure3:add}, which illustrates that for the test tensor $\mathcal{C}$, the processing order $(3,2,1)$ is a reasonable choice, and the processing orders $(1,2,3)$ and $(1,3,2)$ are two bad unreasonable choices in terms of RE (see (\ref{dash-rthosvd:eqn:re})) and CPU run time.}

We now calculate the computational complexities of Algorithms \ref{dash-rthosvd:alg2} and \ref{dash-rthosvd:alg2:v2}. In detail, for each $k$ in Algorithm \ref{dash-rthosvd:alg2}, one has
\begin{enumerate}
    \item[(a)] to draw a standard Gaussian matrix that requires $O(n_1\dots n_{k-1}n_{k+1}\dots n_dl_k)$ operations;
    \item[(b)] to form the mode-$k$ unfolding $\mathbf{A}_{(k)}$ that needs $O(n_1n_2\dots n_d)$ operations;
    \item[(c)] to obtain $\mathbf{Q}_k$ requires $(2q+1)C_{{\rm mm}}n_1n_2\dots n_d l_k+(q+1)C_{{\rm svd}}n_kl_k^2$ operations;
    \item[(d)] it costs $C_{{\rm mm}}r_1\dots r_kn_kn_{k+1}\dots n_d$ operations to  update $\mathcal{G}$.
\end{enumerate}
By summing up the complexity of all the steps above, one has that Algorithm \ref{dash-rthosvd:alg2} requires
\begin{align*}
    &\sum_{k=1}^d(O(n_1\dots n_{k-1}n_{k+1}\dots n_dl_k)+C_{{\rm mm}}r_1\dots r_kn_kn_{k+1}\dots n_d)\\
    &+\sum_{k=1}^d((2q+1)C_{{\rm mm}}n_1n_2\dots n_d l_k+(q+1)C_{{\rm svd}}n_kl_k^2)
\end{align*}
operations to find a quasi-optimal solution for Problem \ref{dash-rthosvd:prob1}. Similarly, for each $k$ in Algorithm \ref{dash-rthosvd:alg2:v2}, one has
\begin{enumerate}
    \item[(a)] to draw a standard Gaussian matrix that requires $O(r_1\dots r_{k-1}n_{k+1}\dots n_dl_k)$ operations;
    \item[(b)] to form the mode-$k$ unfolding $\mathbf{G}_{(k)}$ that needs $O(r_1\dots r_{k-1}n_{k}\dots n_d)$ operations;
    \item[(c)] to obtain $\mathbf{Q}_k$ requires $(2q+1)C_{{\rm mm}}r_1\dots r_{k-1}n_{k}\dots n_d l_k+(q+1)C_{{\rm svd}}n_kl_k^2$ operations;
    \item[(d)] it costs $C_{{\rm mm}}r_1\dots r_kn_kn_{k+1}\dots n_d$ operations to  update $\mathcal{G}$.
\end{enumerate}
Hence, Algorithm \ref{dash-rthosvd:alg2:v2} requires
\begin{align*}
    &\sum_{k=1}^d(O(r_1\dots r_{k-1}n_{k+1}\dots n_dl_k)+C_{{\rm mm}}r_1\dots r_kn_kn_{k+1}\dots n_d)\\
    &+\sum_{k=1}^d((2q+1)C_{{\rm mm}}r_1\dots r_{k-1}n_{k}\dots n_d l_k+(q+1)C_{{\rm svd}}n_kl_k^2)
\end{align*}
operations to find a quasi-optimal solution for Problem \ref{dash-rthosvd:prob1}.
\begin{remark}
    For the same random projection matrices $\{\mathbf{\Omega}_1,\mathbf{\Omega}_2,\dots,\mathbf{\Omega}_d\}$ and the power parameter $q$, the complexity of Algorithms \ref{dash-rthosvd:alg2} and \ref{dash-rthosvd:alg2:v2} are similar, respectively, to the randomized variants, shown in Algorithm \ref{dash-rthosvd:alg1}, of T-HOSVD and ST-HOSVD, respectively. However, numerical examples illustrate that under the same settings, Algorithm \ref{dash-rthosvd:alg2:v2} is faster than the existing randomized variants of T-HOSVD and ST-HOSVD, and Algorithm \ref{dash-rthosvd:alg2} is faster than the existing randomized variants of T-HOSVD, while maintaining similar approximation errors.
\end{remark}
\subsection{Choosing the power parameters in Algorithm \ref{dash-rthosvd:alg2:v2}}
\label{dash-rthosvd:appsec3}

Without loss of generality, we consider the choice of the power parameter $q$ in Algorithm \ref{dash-rthosvd:alg2:v2} to enable automatic termination of the power iteration according to the PVE bound-based accuracy criterion. The cases of other randomized T-HOSVD and ST-HOSVD can be considered similarly.

When we apply Algorithm \ref{dash-rthosvd:alg2:v2} to find a quasi-optimal solution $(\mathbf{U}_1,\mathbf{U}_2,\cdots,\mathbf{U}_d)$ to Problem \ref{dash-rthosvd:prob3}, we suppose that any power parameter $q$ is suitable for each $k$. In general, we assume that any tuple of power parameters $(q_1,q_2,\dots,q_d)$ with $q_k\geq 1$ is needed in Algorithm \ref{dash-rthosvd:alg2:v2}. For each $k$ in Algorithm \ref{dash-rthosvd:alg2:v2}, the power parameter $q_k$ is assumed known as a prior. Without loss of generality, the processing order $(p_1,p_2,\cdots,p_d)$ in Problem \ref{dash-rthosvd:prob5} is set to $(1,2,\cdots,d)$. Note that each $\mathbf{U}_k$ is also a quasi-optimal solution to Problem \ref{dash-rthosvd:prob5} and satisfies that for a given $0<\epsilon<1$, the expectation of inequality $\|\mathbf{U}_k\mathbf{U}_k^\top\mathbf{B}_k-\mathbf{B}_k\|$ is less than or equal to $(1+\epsilon)\|(\mathbf{B}_k)_{r_k}-\mathbf{B}_k\|$ (see \cite[Theorem 2.4]{minster2020randomized}), where $\mathbf{B}_k$ is the mode-$k$ unfolding of $\mathcal{A}\times_1\mathbf{U}_1^\top\dots\times_{k-1}\mathbf{U}_{k-1}^\top$ and $(\mathbf{B}_k)_{r_k}$ is the best rank-$r_k$ approximation to $\mathbf{B}_k$. Hence, for each $k$, the smaller $\epsilon$ is, the larger $q_k$ in Algorithm \ref{dash-rthosvd:alg2:v2} needs to be to achieve the associated accuracy guarantee.

For many practical situations, multiple examinations may be performed to find a suitable and not too large (cf. \cite{che2025efficient,che2025efficient-acom}). This brings a lot of computational cost. Therefore, adaptively and efficiently determining the tuple of power parameters $(q_1,q_2,\dots,q_d)$ should be executed to ensure the required accuracy of the approximate Tucker decomposition with a fixed Tucker-rank, derived from Algorithm \ref{dash-rthosvd:alg2:v2}, is of concern. Algorithm \ref{dash-rthosvd:alg3:app} is proposed by combining Algorithm \ref{dash-rthosvd:alg2:v2} with the PVE bound criterion (cf. \cite{musco2015randomized}). We use this criterion to find a suitable power parameter $q_k$ to ensure the accuracy of the approximate Tucker decomposition with a fixed Tucker-rank. The parameter $\textit{tol}$ in Algorithm \ref{dash-rthosvd:alg3:app} is supposed to be specified by the user and $q_{\max}$ is a parameter denoting the upper bound of the tuple of power parameters $(q_1,\dots,q_d)$.

\begin{algorithm}
     \caption{Modification of Algorithm \ref{dash-rthosvd:alg2:v2} based on the PVE accuracy control}
     \begin{algorithmic}[1]
        \STATEx {\bf Input}: A tensor $\mathcal{A}\in\mathbb{R}^{n_1\times n_2\times\dots\times n_d}$, the desired $d$-tuple $\mathbf{r}=(r_1,r_2,\dots,r_d)$ with $r_k\leq n_k$ and $k=1,2,\dots,d$, the $d$-tuple of oversampling parameters $\mathbf{s}=(s_1,s_2,\dots,s_d)$, the maximum of the power parameter $q_
        {\max}$, and an error tolerance $\textit{tol}$.
        \STATEx {\bf Output}: The power parameters $(q_1,q_2,\dots,q_d)$, and the core tensor $\mathcal{G}\in\mathbb{R}^{r_1\times r_2\times\dots\times r_d}$ and the mode-$k$ factor matrix $\mathbf{U}_k\in\mathbb{R}^{n_k\times r_k}$ such that $\mathcal{A}\approx\widehat{\mathcal{A}}:=\mathcal{G}\times_1\mathbf{U}_1\times_2\mathbf{U}_2\dots\times_d\mathbf{U}_d$.
        \STATE Let $l_k=r_k+s_k$ and $q_k=0$ with $k=1,2,\dots,d$ and set a temporary tensor $\mathcal{G}:=\mathcal{A}$.
        \FOR{$k=1,2,\dots,d$}
            \STATE Draw a standard Gaussian matrix $\mathbf{\Omega}_k\in\mathbb{R}^{r_1\dots r_{k-1}n_{k+1}\dots n_{d}\times l_k}$.
            \STATE Form the mode-$k$ unfolding $\mathbf{G}_{(k)}$ from $\mathcal{G}$.
            \STATE Compute $[\mathbf{Q}_k,\sim,\sim]={\rm svd}(\mathbf{G}_{(k)}\mathbf{\Omega}_k,'{\rm econ}')$ and set ${\color{red}\alpha_{0}^{(k)}}=0$ and ${\color{red}\hat{\bm{\sigma}}_{1}^{(k)}}=\mathbf{0}_{l_k}$.
            \FOR{$t=1,2,\dots,q_{max}$}
                \STATE Update $q_k=q_k+1$.
                \STATE Compute $[\mathbf{Q}_k,\bm{\Sigma}_k,\sim]={\rm svd}(\mathbf{G}_{(k)}(\mathbf{G}_{(k)}^\top\mathbf{Q}_k)-{\color{red}\alpha_{t-1}^{(k)}} \mathbf{Q}_k,'{\rm econ}')$.
                \STATE Compute {\color{red}$\bm{\sigma}_{t}^{(k)}={\rm diag}(\bm{\Sigma}_k)+\alpha_{t-1}^{(k)}$}.
                \IF{$\max\{|{\color{red}\bm{\sigma}_{t}^{(k)}(i)-\hat{\bm{\sigma}}_{t}^{(k)}}|:i=1,2,\dots,r_k\}\leq \textit{tol}\cdot {\color{red}\bm{\sigma}_{t}^{(k)}(i)}(r_k+1)$}
                    \STATE Break.
                \ENDIF
                \IF{$\bm{\Sigma}_k(l_k,l_k)>{\color{red}\alpha_{t-1}^{(k)}}$}
                    \STATE Update ${\color{red}\alpha_{t}^{(k)}}=(\bm{\Sigma}_k(l_k,l_k)+{\color{red}\alpha_{t-1}^{(k)}})/2$.
                \ENDIF
                \STATE Update ${\color{red}\hat{\bm{\sigma}}_{t+1}^{(k)}=\bm{\sigma}_{t}^{(k)}}$.
            \ENDFOR
            \STATE Set $\mathbf{U}_k=\mathbf{Q}_k(:,1:r_k)$ and update $\mathcal{G}=\mathcal{G}\times_k\mathbf{U}_{k}^\top$.
        \ENDFOR
    \end{algorithmic}
    \label{dash-rthosvd:alg3:app}
\end{algorithm}

\begin{remark}
    In Algorithm \ref{dash-rthosvd:alg3:app}, we can obtain another version for Randomized T-HOSVD with the PVE accuracy control (details omitted). The related algorithm is denoted as Algorithm \ref{dash-rthosvd:alg2}-v3 and the version Algorithm \ref{dash-rthosvd:alg3:app} for ST-HOSVDis denoted by Algorithm \ref{dash-rthosvd:alg2:v2}-v3. Numerical comparisons for Algorithm \ref{dash-rthosvd:alg2}-v1/v3 and Algorithm \ref{dash-rthosvd:alg2:v2}-v1/v3 are given in Section \ref{dash-rthosvd:sec5:subsec6}.
\end{remark}

\section{Theoretical analysis}
\label{dash-rthosvd:sec4}
Without loss of generality, for each $k$, let $\widehat{n}_k=\min\{n_k,n_1\cdots n_{k-1}n_{k+1}\cdots n_d\}$ and $\widetilde{n}_k=\min\{n_k,r_1\cdots r_{k-1}n_{k+1}\cdots n_d\}$, where two $d$-tuples $(n_1,n_2,\cdots,n_d)$ and $(r_1,r_2,\cdots,r_d)$ satisfy $r_k<n_k$. It is easy to see that for each $k$, $\widetilde{n}_k\leq \widehat{n}_k$.

Let $\{\mathcal{G},\mathbf{U}_1,\mathbf{U}_2,\dots,\mathbf{U}_d\}$ be obtained from Algorithm \ref{dash-rthosvd:alg2}, then $\mathcal{G}=\mathcal{A}\times_1\mathbf{U}_1^\top\dots\times_d\mathbf{U}_d^\top$. For the case of all the $\mathbf{\Omega}_k$ being standard Gaussian matrices, we obtain an upper bound for $\|\mathcal{G}\times_1\mathbf{U}_1\dots\times_d\mathbf{U}_d-\mathcal{A}\|_F$, which reflects how close the computed low Tucker-rank approximation is to that derived by T-HOSVD.

Similar to (\ref{dash-rthosvd:thosvd-general}), we have
\begin{align}
    \|\mathcal{G}\times_1\mathbf{U}_1\dots\times_d\mathbf{U}_d-\mathcal{A}\|&\leq\sum_{k=1}^d\|\mathcal{A}\times_k(\mathbf{U}_k\mathbf{U}_k^\top)-\mathcal{A}\|.
    \label{dash-rthosvd:thosvd-general-two}
\end{align}
Hence, {\color{red}by obtaining the upper bound for each term $\|\mathcal{A}\times_k(\mathbf{U}_k\mathbf{U}_k^\top)-\mathcal{A}\|_F^2$, we can deduce the upper bound for $\|\mathcal{G}\times_1\mathbf{U}_1\dots\times_d\mathbf{U}_d-\mathcal{A}\|_F^2$, which is summarized in the following theorem}.
\begin{theorem}
    \label{dash-rthosvd:thm1}
   (Proof follows in Section \ref{dash-rthosvd:sec4:sub2}) Let $r_k$, $s_k$, and $q$ be the parameters in Algorithm \ref{dash-rthosvd:alg2} such that {\color{red}$r_k<l_k=r_k+s_k\leq \widehat{n}_k$}. Suppose that the $(d+1)$-tuple $\{\mathcal{G},\mathbf{U}_1,\mathbf{U}_2,\dots,\mathbf{U}_d\}$ is obtained by applying Algorithm \ref{dash-rthosvd:alg2} to $\mathcal{A}\in\mathbb{R}^{n_1\times n_2\times\dots\times n_d}$. If for each $k$, there exist two real numbers $\beta_k>0$ and $\gamma_k>1$ such that $0<\sum_{k=1}^d\Phi_k<1$, where
    \begin{align}
    \label{dash-rthosvd:main-probability}
        \Phi_k=&\frac{1}{\sqrt{2\pi(l_k-r_k+1)}}\left(\frac{e}{(l_k-r_k+1)\beta_k}\right)^{l_k-r_k+1}+\frac{1}{4(\gamma_k^2-
        1)\sqrt{\pi\widehat{n}_k\gamma_k^2}}\left(\frac{2\gamma_k^2}{e^{\gamma_k^2-1}}\right)^{\widehat{n}_k}\nonumber\\
        &+\frac{1}{4(\gamma_k^2-1)\sqrt{\pi n_1\dots n_{k-1}n_{k+1}\dots n_d\gamma_k^2}}\left(\frac{2\gamma_k^2}{e^{\gamma_k^2}-1}\right)^{n_1\dots n_{k-1}n_{k+1}\dots n_d}.
    \end{align}
    Then
    \begin{align*}
        \|\mathcal{G}\times_1\mathbf{U}_1\dots\times_d\mathbf{U}_d-\mathcal{A}\|\leq \sum_{k=1}^d(f_k+g_k)\cdot\Delta_{r_k}(\mathbf{A}_{(k)}),
    \end{align*}
    holds with a probability at least $1-\sum_{k=1}^d\Phi_k$, where
    \begin{align*}
        \Delta_{r_k}(\mathbf{A}_{(k)})=&\sqrt{\sum_{i_k=r_k+1}^{\widehat{n}_k}\sigma_{i_k}(\mathbf{A}_{(k)})^2},\\ f_k=&2\cdot\left(\sqrt{2\widehat{n}_k}\gamma_k\frac{\prod_{t=1}^q\left(\sigma_{r_k+1}(\mathbf{A}_{(k)})^2-\alpha_t^{(k)}\right)}{\prod_{t=1}^q\left(\sigma_{r_k}(\mathbf{A}_{(k)})^2-\alpha_t^{(k)}\right)}+1\right),\\
        g_k=&2\cdot\sqrt{2n_1\dots n_{k-1}n_{k+1}\dots n_dl_k}\beta_k\gamma_k\frac{\prod_{t=1}^q(\sigma_{r_k+1}(\mathbf{A}_{(k)})^2-\alpha_t^{(k)})}{\prod_{t=1}^q\left(\sigma_{r_k}(\mathbf{A}_{(k)})^2-\alpha_{t}^{(k)}\right)}.
    \end{align*}
    Here, for each $t$, $\alpha_t^{(k)}$ denotes the shift value in the $t$-th shifted power iteration.
\end{theorem}
{\color{red}
\begin{remark}
    For each $k$, suppose that the singular values of $\mathbf{A}_{(k)}$ satisfy $\sigma_{r_{k+1}}(\mathbf{A}_{(k)})<\sigma_{r_{k}}(\mathbf{A}_{(k)})$. Then, for each $t$, one has
    \begin{align*}
        \frac{\left(\sigma_{r_{k+1}}(\mathbf{A}_{(k)})^2-\alpha_t^{(k)}\right)}{\left(\sigma_{r_k}(\mathbf{A}_{(k)})^2-\alpha_{t}^{(k)}\right)}<1,
    \end{align*}
which implies that
    \begin{align*}
        \frac{\prod_{t=1}^q\left(\sigma_{r_{k+1}}(\mathbf{A}_{(k)})^2-\alpha_t^{(k)}\right)}{\prod_{t=1}^q\left(\sigma_{r_k}(\mathbf{A}_{(k)})^2-\alpha_{t}^{(k)}\right)}\rightarrow0,\quad \text{as}\quad q\rightarrow \infty.
    \end{align*}
    Hence, when the singular values of $\mathbf{A}_{(k)}$ satisfy $\sigma_{r_{k+1}}(\mathbf{A}_{(k)})<\sigma_{r_{k}}(\mathbf{A}_{(k)})$, then, as $q\rightarrow\infty$, one has $\|\mathcal{G}\times_1\mathbf{U}_1\dots\times_d\mathbf{U}_d-\mathcal{A}\|_F\leq 4\cdot\Delta_{r_k}(\mathbf{A}_{(k)})$ with high probability.
    
    It follows from $f(\alpha)$ being strictly monotonically decreasing with respect to $\alpha$ that
    \begin{equation*}
    \frac{\prod_{t=1}^q\left(\sigma_{r_{k+1}}(\mathbf{A}_{(k)})^2-\alpha_t^{(k)}\right)}{\prod_{t=1}^q\left(\sigma_{r_k}(\mathbf{A}_{(k)})^2-\alpha_{t}^{(k)}\right)}\leq \frac{\prod_{t=1}^q\left(\sigma_{r_{k+1}}(\mathbf{A}_{(k)})^2\right)}{\prod_{t=1}^q\left(\sigma_{r_k}(\mathbf{A}_{(k)})^2\right)}
    \end{equation*}
    which implies that the error factor $$\frac{\prod_{t=1}^q\left(\sigma_{r_{k+1}}(\mathbf{A}_{(k)})^2-\alpha_t^{(k)}\right)}{\prod_{t=1}^q\left(\sigma_{r_k}(\mathbf{A}_{(k)})^2-\alpha_{t}^{(k)}\right)}$$ decays faster with shifting than in the unshifted case. As $q\rightarrow \infty$, this factor converges to 0 at a faster rate, yielding better accuracy with fewer iterations.
    
In Section \ref{dash-rthosvd:sec5:subsec4}, we show that Algorithm \ref{dash-rthosvd:alg2} matches and outperforms HOOI and its randomized variants in terms of speed and approximation errors. Algorithm \ref{dash-rthosvd:alg2} is comparable to T-HOSVD, ST-HOSVD and their existing randomized variants with $q=1$, and better than their randomized algorithms with $q=0$.
\end{remark}
}

Under the case of standard Gaussian matrices, we now consider the upper bound for $\|\mathcal{G}\times_1\mathbf{U}_1\dots\times_d\mathbf{U}_d-\mathcal{A}\|_F$, where $\{\mathcal{G},\mathbf{U}_1,\mathbf{U}_2,\dots,\mathbf{U}_d\}$ is obtained from Algorithm \ref{dash-rthosvd:alg2:v2}. Similar to (\ref{dash-rthosvd:sthosvd-general}), one has
\begin{align}
\label{dash-rthosvd:sthosvd-general-two}
    &\|\mathcal{G}\times_1\mathbf{U}_1\dots\times_d\mathbf{U}_d-\mathcal{A}\|=\left\|\mathcal{A}\times_1\left({\bf U}_1{\bf U}_1^\top\right)\dots\times_d\left({\bf U}_d{\bf U}_d^\top\right)-\mathcal{A}\right\|\nonumber\\
    &\quad\quad\leq\left\|\mathcal{A}\times_{1}\left(\mathbf{I}_{n_{1}}-{\bf U}_{1}{\bf U}_{1}^\top\right)\right\|
    +\left\|\left(\mathcal{A}\times_{1}{\bf U}_{1}^\top\right)\times_{2}\left(\mathbf{I}_{n_{2}}-{\bf U}_{2}{\bf U}_{2}^\top\right)\right\|\nonumber\\
    &\quad\quad\quad+\dots+\left\|\left(\mathcal{A}\times_{1}{\bf U}_{1}^\top\times_{2}\mathbf{U}_{2}^\top\dots\times_{d-1}\mathbf{U}_{d-1}^\top\right)\times_{d}\left(\mathbf{I}_{n_{d}}-{\bf U}_{d}{\bf U}_{d}^\top\right)\right\|.
\end{align}

For clarity, let $\mathcal{B}_1=\mathcal{A}$ and {\color{red}$\mathcal{B}_k=\mathcal{A}\times_{1}\mathbf{U}_{1}^\top\dots\times_{k-1}\mathbf{U}_{k-1}^\top$} with $k=2,3,\dots,d$. For each $k$, the mode-$k$ unfolding of $\mathcal{B}_k$ is denoted by $\mathbf{B}_k$. Hence, we rewrite (\ref{dash-rthosvd:sthosvd-general-two}) as follows:
\begin{align*}
    \|\mathcal{G}\times_1\mathbf{U}_1\dots\times_d\mathbf{U}_d-\mathcal{A}\|\leq\sum_{k=1}^d\left\|{\bf U}_{k}{\bf U}_{k}^\top\mathbf{B}_k-\mathbf{B}_k\right\|.
\end{align*}

We now give a rough upper bound for Algorithm \ref{dash-rthosvd:alg2:v2}, which is illustrated in the following theorem.
\begin{theorem}
\label{dash-rthosvd:thm2}
    (Proof follows in Section \ref{dash-rthosvd:sec4:sub3}) Suppose that $r_k$, $s_k$, and $q$ are the parameters in Algorithm \ref{dash-rthosvd:alg2:v2} such that {\color{red}$r_k<l_k=r_k+s_k\leq \widehat{n}_k$}. Let $\{\mathcal{G},\mathbf{U}_1,\mathbf{U}_2,\dots,\mathbf{U}_d\}$ be obtained by applying Algorithm \ref{dash-rthosvd:alg2:v2} to $\mathcal{A}\in\mathbb{R}^{n_1\times n_2\times\dots\times n_d}$. 
    
    If for each $k$, there exists a positive integer $j_k$ with $j_k<r_k$ and two real numbers $\beta_k>0$ and $\gamma_k>1$ such that $0<\sum_{k=1}^d\Phi_k<1$, where for each $k$, $\Phi_k$ is given in (\ref{dash-rthosvd:main-probability}). Then, one has
    \begin{align*}
        &\|\mathcal{G}\times_1\mathbf{U}_1\dots\times_d\mathbf{U}_d-\mathcal{A}\|_F\leq2\cdot\sum_{k=1}^d\left(\sqrt{2\widehat{n}_k}\gamma_k
        +1\right)\sqrt{\sum_{i_k=r_k+1}^{\widehat{n}_k}\sigma_{i_k}(\mathbf{A}_{(k)})^2}\\
        &+2\cdot\sum_{k=1}^d\sqrt{2n_1\dots n_{k-1}n_{k+1}\dots n_dl_k}\beta_k\gamma_k\sqrt{\sum_{i_k=r_k+1}^{\widehat{n}_k}\sigma_{i_k}(\mathbf{A}_{(k)})^2}
    \end{align*}
    holds with a probability at least $1-\sum_{k=1}^d\Phi_k$.
\end{theorem}
\begin{remark}
    Following Theorem \ref{dash-rthosvd:thm2}, we cannot ensure that one has $\|\mathcal{G}\times_1\mathbf{U}_1\dots\times_d\mathbf{U}_d-\mathcal{A}\|_F\leq 4\cdot\Delta_{r_k}(\mathbf{A}_{(k)})$ with high probability when the singular values of $\mathbf{A}_{(k)}$ satisfy $\sigma_{r_{k+1}}(\mathbf{A}_{(k)})<\sigma_{r_{k}}(\mathbf{A}_{(k)})$, then, as $q\rightarrow\infty$.
\end{remark}
\begin{remark}
In this paper, we only obtain high-probability Frobenius-norm bounds (see Theorems \ref{dash-rthosvd:thm1} and \ref{dash-rthosvd:thm2}), which are standard for randomized tensor methods. Expected-error bounds can be derived using similar techniques but are omitted for brevity and are left as future work.
\end{remark}
\subsection{Basic lemmas}
We first introduce two lemmas to illustrate the upper bound on the largest singular value and the lower bound on the least singular value of any standard Gaussian matrix.
\begin{lemma}{{\bf (see \cite[Lemma 2.1]{rokhlin2010randomized})}}
\label{dash-rthosvd:lem1}
    Let $\mathbf{\Omega}\in\mathbb{R}^{l\times n}$ be a standard Gaussian matrix with $l<n$. For $\gamma>1$, suppose that
    \begin{equation}\label{dash-rthosvd:app:eqn2}
	1-\frac{1}{4(\gamma^2-
        1)\sqrt{\pi
        n\gamma^2}}\left(\frac{2\gamma^2}{e^{\gamma^2-1}}\right)^n\geq 0,
    \end{equation}
    then, the largest singular value of $\mathbf{\Omega}$ is at most $\sqrt{2n}\gamma$ with a probability not less than the amount in {\rm (\ref{dash-rthosvd:app:eqn2})}.
\end{lemma}
\begin{lemma}{{\bf (see \cite[Lemma 2.2]{rokhlin2010randomized})}}
\label{dash-rthosvd:lem2}
    For $\beta>0$, let $n$ and $l$ be positive integers such that $l\leq n$ and
    \begin{equation}\label{dash-rthosvd:app:eqn3}
	    1-\frac{1}{\sqrt{2\pi(n-
        l+1)}}\left(
        \frac{e}{(n-
        l+1)\beta}\right)^{n-l+1}\geq 0.
    \end{equation}
    Suppose that $\mathbf{\Omega}\in\mathbb{R}^{l\times n}$ is a standard Gaussian matrix. Then, the smallest singular value of $\mathbf{\Omega}$ is at least $1/(\sqrt{n}\beta)$ with a probability not less than the amount in {\rm (\ref{dash-rthosvd:app:eqn3})}.
\end{lemma}
Two other lemmas, which will be used in the proof of Theorem \ref{dash-rthosvd:thm1}, are also introduced.
\begin{lemma}{{\bf (see \cite[Lemma 3]{che2025efficient-acom})}}
\label{dash-rthosvd:lem3}
    Let $k$, $m$ and $n$ be positive integers with $k<n\leq m$. For a given matrix $\mathbf{A}\in\mathbb{R}^{m\times n}$, there exists an orthonormal matrix $\mathbf{Q}\in\mathbb{R}^{m\times k}$ and $\mathbf{S}\in\mathbb{R}^{k\times n}$ such that
    \begin{equation*}
       \|\mathbf{Q}\mathbf{S}-\mathbf{A}\|=\left(\sum_{i=k+1}^n\sigma_i(\mathbf{A})^2\right)^{1/2},
    \end{equation*}
    where $\sigma_i(\mathbf{A})$ is the $i$-th singular value of $\mathbf{A}$ and $i=1,2,\dots,n$.
\end{lemma}
\begin{lemma}{{\bf (see \cite[Corollary 5.1]{che2021randomized})}}
\label{dash-rthosvd:lem4}
    Let $\mathbf{A}\in\mathbb{R}^{m\times n}$ and $\mathbf{B}\in\mathbb{R}^{n\times p}$ with $p\leq\min\{m,n\}$. Then for all $k=1,2,\dots,p$, we have
    \begin{equation*}
        \sum_{i=k}^{p}\sigma_i
        (\mathbf{A}\mathbf{B})^2\leq
        \|\mathbf{B}\|_2^2\sum_{i=k}^{\min
        \{m,n\}}\sigma_i(\mathbf{A})^2.
    \end{equation*}
\end{lemma}
\begin{remark}
    Let $\mathbf{B}\in\mathbb{R}^{n\times p}$ with $p\leq n$ be any standard Gaussian matrix. Martinsson {\it et al.} \cite{martinsson2011randomized} considered a tighter upper bound for $\sigma_k
    (\mathbf{A}\mathbf{B})$ with $k=1,2,\dots,p$ (see Lemma 4.6 in \cite{martinsson2011randomized}). In this paper, we consider a rougher upper bound for  $\sigma_k(\mathbf{A}\mathbf{B})$ by using Lemmas \ref{dash-rthosvd:lem1} and \ref{dash-rthosvd:lem4}. In detail, one has
    \begin{align*}
        \sqrt{\sum_{i=k}^{p}\sigma_i
        (\mathbf{A}\mathbf{B})^2}\leq
        \sqrt{2n}\gamma\sqrt{\sum_{i=k}^{\min
        \{m,n\}}\sigma_i(\mathbf{A})^2}
    \end{align*}
    with a probability at least
    \begin{equation*}
	1-\frac{1}{4(\gamma^2-
        1)\sqrt{\pi
        n\gamma^2}}\left(\frac{2\gamma^2}{e^{\gamma^2-1}}\right)^n,
    \end{equation*}
    where $\gamma>1$ satisfies
    \begin{align*}
        0<\frac{1}{4(\gamma^2-
        1)\sqrt{\pi
        n\gamma^2}}\left(\frac{2\gamma^2}{e^{\gamma^2-1}}\right)^n<1.
    \end{align*}
\end{remark}
\subsection{Proof of Theorem \ref{dash-rthosvd:thm1}}
\label{dash-rthosvd:sec4:sub2}
We now consider the term $\|\mathcal{A}\times_k(\mathbf{U}_k\mathbf{U}_k^\top)-\mathcal{A}\|_F^2$ with $k=1,2,\cdots,d$ in the right-hand side of (\ref{dash-rthosvd:thosvd-general-two}). For a given $k$, the matrix $\mathbf{U}_k\in\mathbb{R}^{n_k\times r_k}$ is obtained by Algorithm \ref{dash-rthosvd:alg2}, so the matrix $\mathbf{U}_k\mathbf{U}_k^\top\mathbf{A}_{(k)}$ is a good approximation to the matrix $\mathbf{A}_{(k)}$, provided that there exists matrices $\mathbf{\Omega}_k\in\mathbb{R}^{n_1\dots n_{k-1}n_{k+1}\dots n_d\times l_k}$ and $\mathbf{R}_k\in\mathbb{R}^{r_k\times l_k}$ such that:
\begin{enumerate}
    \item[(a)] $\mathbf{U}_k$ is orthonormal;
    \item[(b)] $\mathbf{U}_k\mathbf{R}_k$ is a good approximation of $\prod_{t=1}^q(\mathbf{A}_{(k)}\mathbf{A}_{(k)}^\top-\alpha_t^{(k)}\mathbf{I}_{n_k})\mathbf{A}_{(k)}\mathbf{\Omega}_k$;
    \item[(c)] there exists a matrix $\mathbf{F}_k\in\mathbb{R}^{l_k\times n_1\dots n_{k-1}n_{k+1}\dots n_d}$ such that $\prod_{t=1}^q(\mathbf{A}_{(k)}\mathbf{A}_{(k)}^\top-\alpha_t^{(k)}\mathbf{I}_{n_k})\mathbf{A}_{(k)}\mathbf{\Omega}_k\mathbf{F}_k$ is a good approximation of $\mathbf{A}_{(k)}$ and  $\|\mathbf{F}_k\|_2$ is not too large,
\end{enumerate}
where for each $t$, $\alpha_t^{(k)}$ denotes the shift value in the $t$-th shifted power iteration.
\begin{lemma}
\label{dash-rthosvd:lem5}
    For a given $k$, suppose that $\mathbf{A}_{(k)}\in\mathbb{R}^{n_k\times n_1\dots n_{k-1}n_{k+1}\dots n_d}$,
    $\mathbf{U}_k\in\mathbb{R}^{n_k\times r_k}$ is orthonormal, $\mathbf{F}_k\in\mathbb{R}^{l_k\times n_1\dots n_{k-1}n_{k+1}\dots n_d}$, $\mathbf{R}_k\in\mathbb{R}^{r_k\times l_k}$, and $\mathbf{\Omega}_k\in\mathbb{R}^{n_1\dots n_{k-1}n_{k+1}\dots n_d\times l_k}$ with $r_k<l_k=r_k+s_k\leq \widehat{n}_k$. For a given positive integer $q>0$, we have
    \begin{align}\label{dash-rthosvd:app:eqn4}
        \|\mathcal{A}\times_k(\mathbf{U}_k\mathbf{U}_k^\top)-\mathcal{A}\|^2
	      &\leq 2\left\|\prod_{t=1}^q\left(\mathbf{A}_{(k)}\mathbf{A}_{(k)}^\top-\alpha_t^{(k)}
        \mathbf{I}_{n_k}\right)\mathbf{A}_{(k)}\mathbf{\Omega}_k\mathbf{F}_k-\mathbf{A}_{(k)}\right\|^2\nonumber\\
          &+2\|\mathbf{F}_k\|_2^2\left\|\mathbf{U}_k\mathbf{R}_k-\prod_{t=1}^q\left(\mathbf{A}_{(k)}\mathbf{A}_{(k)}^\top
          -\alpha_t^{(k)}\mathbf{I}_{n_k}\right)\mathbf{A}_{(k)}\mathbf{\Omega}_k\right\|^2,
    \end{align}
    where for each $t$, $\alpha_t^{(k)}$ denotes the shift value in the $t$-th shifted power iteration.
\end{lemma}

\begin{proof}
	The proof is straightforward, but tedious, as follows. By using the triangular inequality, we have
	\begin{equation}\label{dash-rthosvd:eqn11}
	\begin{aligned}
	&\|\mathbf{A}_{(k)}-\mathbf{U}_k\mathbf{U}_k^\top\mathbf{A}_{(k)}\|^2
 \leq\left\|\mathbf{A}_{(k)}-\prod_{t=1}^q\left(\mathbf{A}_{(k)}\mathbf{A}_{(k)}^\top-\alpha_t^{(k)}
        \mathbf{I}_{n_k}\right)\mathbf{A}_{(k)}\mathbf{\Omega}_k\mathbf{F}_k\right\|^2\\
	&+\left\|\prod_{t=1}^q\left(\mathbf{A}_{(k)}\mathbf{A}_{(k)}^\top-\alpha_t^{(k)}
        \mathbf{I}_{n_k}\right)\mathbf{A}_{(k)}\mathbf{\Omega}_k\mathbf{F}_k
	-\mathbf{U}_k\mathbf{U}_k^\top\prod_{t=1}^q\left(\mathbf{A}_{(k)}\mathbf{A}_{(k)}^\top-\alpha_t^{(k)}
        \mathbf{I}_{n_k}\right)\mathbf{A}_{(k)}\mathbf{\Omega}_k\mathbf{F}_k\right\|^2\\
	&+\left\|\mathbf{U}_k\mathbf{U}_k^\top\prod_{t=1}^q\left(\mathbf{A}_{(k)}\mathbf{A}_{(k)}^\top-\alpha_t^{(k)}
        \mathbf{I}_{n_k}\right)\mathbf{A}_{(k)}\mathbf{\Omega}_k\mathbf{F}_k
	-\mathbf{U}_k\mathbf{U}_k^\top\mathbf{A}_{(k)}\right\|^2.
	\end{aligned}
	\end{equation}
	Note that we have
	\begin{equation*}
	\begin{aligned}
&\left\|\mathbf{U}_k\mathbf{U}_k^\top\prod_{t=1}^q\left(\mathbf{A}_{(k)}\mathbf{A}_{(k)}^\top-\alpha_t^{(k)}
        \mathbf{I}_{n_k}\right)\mathbf{A}_{(k)}\mathbf{\Omega}_k\mathbf{F}_k
	-\mathbf{U}_k\mathbf{U}_k^\top\mathbf{A}_{(k)}\right\|^2\\
 &\leq
	\left\|\prod_{t=1}^q\left(\mathbf{A}_{(k)}\mathbf{A}_{(k)}^\top-\alpha_t^{(k)}
        \mathbf{I}_{n_k}\right)\mathbf{A}_{(k)}\mathbf{\Omega}_k\mathbf{F}_k
	-\mathbf{A}_{(k)}\right\|^2\|\mathbf{U}_k\mathbf{U}_k^\top\|_2^2,
	\end{aligned}
	\end{equation*}
	which implies that
	\begin{equation}\label{dash-rthosvd:eqn12}
	\begin{aligned}
	&\left\|\mathbf{U}_k\mathbf{U}_k^\top\prod_{t=1}^q\left(\mathbf{A}_{(k)}\mathbf{A}_{(k)}^\top-\alpha_t^{(k)}
        \mathbf{I}_{n_k}\right)\mathbf{A}_{(k)}\mathbf{\Omega}_k\mathbf{F}_k
	-\mathbf{U}_k\mathbf{U}_k^\top\mathbf{A}_{(k)}\right\|^2\\
	&\leq
	\left\|\prod_{t=1}^q\left(\mathbf{A}_{(k)}\mathbf{A}_{(k)}^\top-\alpha_t^{(k)}
        \mathbf{I}_{n_k}\right)\mathbf{A}_{(k)}\mathbf{\Omega}_k\mathbf{F}_k
	-\mathbf{A}_{(k)}\right\|^2.
	\end{aligned}
	\end{equation}
	For the second term in the right-hand side of (\ref{dash-rthosvd:eqn11}), we have
	\begin{equation*}
	\begin{aligned}
	&\left\|\prod_{t=1}^q\left(\mathbf{A}_{(k)}\mathbf{A}_{(k)}^\top-\alpha_t^{(k)}
        \mathbf{I}_{n_k}\right)\mathbf{A}_{(k)}\mathbf{\Omega}_k\mathbf{F}_k
	-\mathbf{U}_k\mathbf{U}_k^\top\prod_{t=1}^q\left(\mathbf{A}_{(k)}\mathbf{A}_{(k)}^\top-\alpha_t^{(k)}
        \mathbf{I}_{n_k}\right)\mathbf{A}_{(k)}\mathbf{\Omega}_k\mathbf{F}_k\right\|^2\\
	&\leq\left\|\prod_{t=1}^q\left(\mathbf{A}_{(k)}\mathbf{A}_{(k)}^\top-\alpha_t^{(k)}
        \mathbf{I}_{n_k}\right)\mathbf{A}_{(k)}\mathbf{\Omega}_k
	-\mathbf{U}_k\mathbf{U}_k^\top\prod_{t=1}^q\left(\mathbf{A}_{(k)}\mathbf{A}_{(k)}^\top-\alpha_t^{(k)}
        \mathbf{I}_{n_k}\right)\mathbf{A}_{(k)}\mathbf{\Omega}_k\right\|^2\|\mathbf{F}_k\|_2^2.
	\end{aligned}
	\end{equation*}
	From the triangular inequality, we have
	\begin{equation*}
	\begin{aligned}
	&\left\|\prod_{t=1}^q\left(\mathbf{A}_{(k)}\mathbf{A}_{(k)}^\top-\alpha_t^{(k)}
        \mathbf{I}_{n_k}\right)\mathbf{A}_{(k)}\mathbf{\Omega}_k
	-\mathbf{U}_k\mathbf{U}_k^\top\prod_{t=1}^q\left(\mathbf{A}_{(k)}\mathbf{A}_{(k)}^\top-\alpha_t^{(k)}
        \mathbf{I}_{n_k}\right)\mathbf{A}_{(k)}\mathbf{\Omega}_k\right\|^2\\
	&\quad\leq\left\|\prod_{t=1}^q\left(\mathbf{A}_{(k)}\mathbf{A}_{(k)}^\top-\alpha_t^{(k)}
        \mathbf{I}_{n_k}\right)\mathbf{A}_{(k)}\mathbf{\Omega}_k-\mathbf{U}_k\mathbf{R}_k\right\|^2+\|\mathbf{U}_k\mathbf{R}_k-\mathbf{U}_k\mathbf{U}_k^\top\mathbf{U}_k\mathbf{R}_k\|^2\\
        &\quad\quad+\left\|\mathbf{U}_k\mathbf{U}_k^\top\mathbf{U}_k\mathbf{R}_k
	-\mathbf{U}_k\mathbf{U}_k^\top\prod_{t=1}^q\left(\mathbf{A}_{(k)}\mathbf{A}_{(k)}^\top-\alpha_t^{(k)}
        \mathbf{I}_{n_k}\right)\mathbf{A}_{(k)}\mathbf{\Omega}_k\right\|^2.
	\end{aligned}
	\end{equation*}
	By the fact that $\mathbf{U}_k^\top\mathbf{U}_k=\mathbf{I}_{r_k}$, we have $\|\mathbf{U}_k\mathbf{R}_k-\mathbf{U}_k\mathbf{U}_k^\top\mathbf{U}_k\mathbf{R}_k\|^2=0$ and
	\begin{equation*}
	\begin{aligned}
	&\left\|\mathbf{U}_k\mathbf{U}_k^\top\mathbf{U}_k\mathbf{R}_k
	-\mathbf{U}_k\mathbf{U}_k^\top\prod_{t=1}^q\left(\mathbf{A}_{(k)}\mathbf{A}_{(k)}^\top-\alpha_t^{(k)}
        \mathbf{I}_{n_k}\right)\mathbf{A}_{(k)}\mathbf{\Omega}_k\right\|^2\\
	&\quad\quad\leq\left\|\mathbf{U}_k\mathbf{R}_k
	-\prod_{t=1}^q\left(\mathbf{A}_{(k)}\mathbf{A}_{(k)}^\top-\alpha_t^{(k)}
        \mathbf{I}_{n_k}\right)\mathbf{A}_{(k)}\mathbf{\Omega}_k\right\|^2.
	\end{aligned}
	\end{equation*}
	Therefore,
	\begin{equation}\label{dash-rthosvd:eqn13}
	\begin{aligned}
	&\left\|(\mathbf{I}_{n_k}-\mathbf{U}_k\mathbf{U}_k^\top)\prod_{t=1}^q\left(\mathbf{A}_{(k)}\mathbf{A}_{(k)}^\top-\alpha_t^{(k)}
        \mathbf{I}_{n_k}\right)\mathbf{A}_{(k)}\mathbf{\Omega}_k\right\|^2\\
        &\quad\leq 2\left\|\mathbf{U}_k\mathbf{R}_k
	-\prod_{t=1}^q\left(\mathbf{A}_{(k)}\mathbf{A}_{(k)}^\top-\alpha_t^{(k)}
        \mathbf{I}_{n_k}\right)\mathbf{A}_{(k)}\mathbf{\Omega}_k\right\|^2.
	\end{aligned}
	\end{equation}
	Combining (\ref{dash-rthosvd:eqn11}), (\ref{dash-rthosvd:eqn12}) and (\ref{dash-rthosvd:eqn13}) yields (\ref{dash-rthosvd:app:eqn4}).
\end{proof}

{\color{red}
\begin{remark}
When $\alpha_t^{(k)}=0$ and the Frobenius norm is replaced by the spectral norm, the inequality (\ref{dash-rthosvd:app:eqn4}) is rewritten as
    \begin{align*}
      \|\mathbf{A}_{(k)}-\mathbf{U}_k\mathbf{U}_k^\top\mathbf{A}_{(k)}\|^2
       &\leq 2\left\|(\mathbf{A}_{(k)}\mathbf{A}_{(k)}^\top)^q\mathbf{A}_{(k)}\mathbf{\Omega}_k\mathbf{F}_k-\mathbf{A}_{(k)}\right\|^2\nonumber\\
       &+2\|\mathbf{F}_k\|_2^2\left\|\mathbf{U}_k\mathbf{R}_k-(\mathbf{A}_{(k)}\mathbf{A}_{(k)}^\top)^q\mathbf{A}_{(k)}\mathbf{\Omega}_k\right\|^2,
    \end{align*}
    which is the same as (A.1) in \cite[Lemma A.1]{rokhlin2010randomized}. When $\alpha_t^{(k)}=0$, the inequality (\ref{dash-rthosvd:app:eqn4}) is rewritten as
    \begin{align*}
      \|\mathbf{A}_{(k)}-\mathbf{U}_k\mathbf{U}_k^\top\mathbf{A}_{(k)}\|^2
       &\leq 2\left\|(\mathbf{A}_{(k)}\mathbf{A}_{(k)}^\top)^q\mathbf{A}_{(k)}\mathbf{\Omega}_k\mathbf{F}_k-\mathbf{A}_{(k)}\right\|^2\nonumber\\
       &+2\|\mathbf{F}_k\|_2^2\left\|\mathbf{U}_k\mathbf{R}_k-(\mathbf{A}_{(k)}\mathbf{A}_{(k)}^\top)^q\mathbf{A}_{(k)}\mathbf{\Omega}_k\right\|^2,
    \end{align*}
    which is the same as (4.1) in \cite[Lemma 4.5]{che2020computation}.
\end{remark}
}

By combining Lemmas \ref{dash-rthosvd:lem1}, \ref{dash-rthosvd:lem3} and \ref{dash-rthosvd:lem4}, the upper bound for the second part in the right-hand side of (\ref{dash-rthosvd:app:eqn4}) is given in the following theorem.
\begin{theorem}
\label{dash-rthosvd:app:thm1}
    For a given $k$, suppose that
    $\mathbf{\Omega}_k\in\mathbb{R}^{n_1\dots n_{k-1}n_{k+1}\dots n_d\times (r_k+p_k)}$ is any standard Gaussian matrix with $l_k=r_k+p_k<\widehat{n}_k$ and $\mathbf{A}_{(k)}\in\mathbb{R}^{n_k\times n_1\dots n_{k-1}n_{k+1}\dots n_d}$ is the mode-$k$ unfolding of $\mathcal{A}\in\mathbb{R}^{n_1\times n_2\times \cdots\times n_d}$. For $\gamma_k>1$, let
    \begin{align}\label{dash-rthosvd:app:eqn5}
        1-\frac{1}{4(\gamma_k^2-1)\sqrt{\pi n_1\dots n_{k-1}n_{k+1}\dots n_d\gamma_k^2}}\left(\frac{2\gamma_k^2}{e^{\gamma_k^2}-1}\right)^{n_1\dots n_{k-1}n_{k+1}\dots n_d}>0.
    \end{align}
    For a given positive integer $q>0$, there exists an orthonormal matrix $\mathbf{U}_k\in\mathbb{R}^{n_k\times r_k}$ and $\mathbf{R}_k\in\mathbb{R}^{r_k\times l_k}$ such that
    \begin{align*}
        &\left\|\mathbf{U}_k\mathbf{R}_k-\prod_{t=1}^q\left(\mathbf{A}_{(k)}\mathbf{A}_{(k)}^\top-\alpha_t^{(k)}\mathbf{I}_{n_k}\right)\mathbf{A}_{(k)}\mathbf{\Omega}_k\right\|\\
        &\leq \sqrt{2n_1\dots n_{k-1}n_{k+1}\dots n_d}\gamma_k\prod_{t=1}^q\left(\sigma_{r_k+1}(\mathbf{A}_{(k)})^2-\alpha_t^{(k)}\right)\sqrt{\sum_{i_k=r_k+1}^{\widehat{n}_k}\sigma_{i_k}(\mathbf{A}_{(k)})^2}
    \end{align*}
    with a probability at least the amount in (\ref{dash-rthosvd:app:eqn5}), where for each $t$, $\alpha_t^{(k)}$ denotes the shift value in the $t$-th shifted power iteration.
\end{theorem}

For each $k$, we now consider the upper bounds for $\|\mathbf{F}_k\|_2$ and $\|\prod_{t=1}^q(\mathbf{A}_{(k)}\mathbf{A}_{(k)}^\top-
\alpha_t^{(k)}\mathbf{I}_{n_k})\mathbf{A}_{(k)}\mathbf{\Omega}_k\mathbf{F}_k-\mathbf{A}_{(k)}\|^2$, which are summarized in the following theorem.
\begin{theorem}
\label{dash-rthosvd:app:thm2}
    For a given $k$, suppose that
    $\mathbf{\Omega}_k\in\mathbb{R}^{n_1\dots n_{k-1}n_{k+1}\dots n_d\times l_k}$ is any standard Gaussian matrix with $r_k<l_k=r_k+p_k<\widehat{n}_k$, and $\mathbf{A}_{(k)}\in\mathbb{R}^{n_k\times n_1\dots n_{k-1}n_{k+1}\dots n_d}$ is the mode-$k$ unfolding of $\mathcal{A}\in\mathbb{R}^{n_1\times n_2\times \cdots\times n_d}$. For $\beta_k>0$ and $\gamma_k>1$, let
    \begin{align}\label{dash-rthosvd:app:eqn6}
        1&-\frac{1}{\sqrt{2\pi(l_k-r_k+1)}}\left(\frac{e}{(l_k-r_k+1)\beta_k}\right)^{l_k-r_k+1}\nonumber\\
        &-\frac{1}{4(\gamma_k^2-
        1)\sqrt{\pi\widehat{n}_k\gamma_k^2}}\left(\frac{2\gamma_k^2}{e^{\gamma_k^2-1}}\right)^{\widehat{n}_k}
    \end{align}
    be nonnegative. For a given positive integer $q>0$, there exists an orthonormal matrix $\mathbf{U}_k\in\mathbb{R}^{n_k\times r_k}$ and $\mathbf{R}_k\in\mathbb{R}^{r_k\times (r_k+p_k)}$ such that
    \begin{align}\label{dash-rthosvd:app:eqn7}
        &\left\|\prod_{t=1}^q\left(\mathbf{A}_{(k)}\mathbf{A}_{(k)}^\top-\alpha_t^{(k)}\mathbf{I}_{n_k}\right)\mathbf{A}_{(k)}\mathbf{\Omega}_k\mathbf{F}_k-\mathbf{A}_{(k)}\right\|\nonumber\\
        &\leq\left(\sqrt{2\widehat{n}_k}\gamma_k
        \frac{\prod_{t=1}^q\left(\sigma_{r_k+1}(\mathbf{A}_{(k)})^2-\alpha_t^{(k)}\right)}{\prod_{t=1}^q\left(\sigma_{r_k}(\mathbf{A}_{(k)})^2-\alpha_t^{(k)}\right)}+1\right)\sqrt{\sum_{i_k=r_k+1}^{\widehat{n}_k}\sigma_{i_k}(\mathbf{A}_{(k)})^2}
    \end{align}
    and
    \begin{align}\label{dash-rthosvd:app:eqn8}
        \|\mathbf{F}_k\|_2\leq \frac{\sqrt{l_k}\beta_k}{\prod_{t=1}^q\left(\sigma_{r_k}(\mathbf{A}_{(k)})^2-\alpha_{t}^{(k)}\right)}
    \end{align}
    with a probability at least the amount in (\ref{dash-rthosvd:app:eqn6}), where for each $t$, $\alpha_t^{(k)}$ denotes the shift value in the $t$-th shifted power iteration.
\end{theorem}
\begin{remark}
    When $\alpha_t^{(k)}=0$, Theorem \ref{dash-rthosvd:app:thm2} is similar to Lemma 3.2 in \cite{rokhlin2010randomized}, where the former is based on the Frobenius norm and the latter is based on the spectral norm.
\end{remark}
\begin{proof}
    The detailed process for proving this theorem is followed by the proof of Lemma A.2 in \cite{rokhlin2010randomized} by constructing the matrix $\mathbf{F}_k$ to satisfy (\ref{dash-rthosvd:app:eqn7}) and (\ref{dash-rthosvd:app:eqn8}).

Suppose that $\mathbf{A}_{(k)}=\mathbf{Q}_k\mathbf{\Sigma}_k\mathbf{P}_k^\top$, where $\mathbf{Q}_k\in\mathbb{R}^{n_k\times \widehat{n}_k}$ and $\mathbf{P}_k\in\mathbb{R}^{n_1\dots n_{k-1}n_{k+1}\dots n_d\times \widehat{n}_k}$ are orthonormal, and the diagonal entries of the diagonal matrix $\mathbf{\Sigma}_k\in\mathbb{R}^{\widehat{n}_k\times \widehat{n}_k}$ are nonnegative and arranged in descending order. Then, we have
    \begin{align*}
        \prod_{t=1}^q\left(\mathbf{A}_{(k)}\mathbf{A}_{(k)}^\top-\alpha_t^{(k)}\mathbf{I}_{n_k}\right)\mathbf{A}_{(k)}
        =\mathbf{Q}_k\prod_{t=1}^q\left(\mathbf{\Sigma}_k^2-\alpha_t^{(k)}\mathbf{I}_{\widehat{n}_k}\right)\mathbf{\Sigma}_k\mathbf{P}_k^\top.
    \end{align*}
    According to the given $\mathbf{P}_k^\top$ and $\mathbf{\Omega}_k$, suppose that
    \begin{equation*}
        \mathbf{P}_k^\top\mathbf{\Omega}_k=
        \begin{pmatrix}
             \mathbf{H}_k\\
             \mathbf{\Phi}_k
        \end{pmatrix},
    \end{equation*}
    where $\mathbf{H}_k\in\mathbb{R}^{r_k\times l_k}$ and $\mathbf{\Phi}_k\in\mathbb{R}^{(n_1\dots n_{k-1}n_{k+1}\dots n_d-r_k)\times l_k}$. Since $\mathbf{\Omega}_k$ is a standard Gaussian matrix, and $\mathbf{P}_k$ is an orthonormal matrix, then $\mathbf{P}_k^\top\mathbf{\Omega}_k$ is also a standard Gaussian matrix, which implies that $\mathbf{H}_k$ and $\mathbf{\Phi}_k$ are standard Gaussian matrices. It follows from Lemma \ref{dash-rthosvd:lem2} that $\mathbf{H}_k\mathbf{H}_k^\top$ is invertible, and for any $\beta_k>0$, one has $\|\mathbf{H}_k^\dag\|_2\leq \sqrt{l_k}\beta_k$ with a probability at least
 \begin{align}
 \label{dash-rthosvd:app:eqn9}
    1-\frac{1}{\sqrt{2\pi(l_k-r_k+1)}}\left(\frac{e}{(l_k-r_k+1)\beta_k}\right)^{l_k-r_k+1},
 \end{align}
where the Moore-Penrose inverse of $\mathbf{H}_k$ is given by $\mathbf{H}_k^\dag=\mathbf{H}_k^{\top}(\mathbf{H}_k\mathbf{H}_k^\top)^{-1}$. 

For clarity, let $\mathbf{\Sigma}_{k,1}=\mathbf{\Sigma}_k(1:r_k,1:r_k)$ and $\mathbf{\Sigma}_{k,2}=\mathbf{\Sigma}_k(r_k+1:\widehat{n}_k,r_k+1:\widehat{n}_k)$. Suppose that $\mathbf{S}_k\in\mathbb{R}^{l_k\times r_k}$ is given by
\begin{equation*}
   \mathbf{S}_k=
      \begin{pmatrix}
         \mathbf{H}_k^\dag\prod_{t=1}^q(\mathbf{\Sigma}_{k,1}-\alpha_t^{(k)}\mathbf{I}_{r_k})^{-1}&
         \mathbf{0}_{1}
       \end{pmatrix},
\end{equation*}
where $\mathbf{0}_1$ is the zero matrix of size $l_k\times (\widehat{n}_k-r_k)$. Define $\mathbf{F}_k=\mathbf{S}_k\mathbf{P}_k^\top$. From Lemma \ref{dash-rthosvd:lem2}, for any $\beta_k>1$, the upper bound for  $\|\mathbf{F}_k\|_2$ is given by
    \begin{align*}
        \|\mathbf{F}_k\|_2&=
        \|\mathbf{S}_k\mathbf{P}_k^\top\|_2=
        \left\|\mathbf{H}_k^\dag\prod_{t=1}^q(\mathbf{\Sigma}_{k,1}-\alpha_t^{(k)}\mathbf{I}_{r_k})^{-1}\right\|_2\\
        &\leq\frac{1}{\prod_{t=1}^q(\sigma_{r_k}(\mathbf{A}_{(k)})^2-\alpha_{t}^{(k)})}\frac{1}{\sigma_{\min}(\mathbf{H}_k)}
        \leq\frac{\sqrt{l_k}\beta_k}{\prod_{t=1}^q(\sigma_{r_k}(\mathbf{A}_{(k)})^2-\alpha_{t}^{(k)})}
    \end{align*}
    with a probability at least the amount in (\ref{dash-rthosvd:app:eqn9}).

    Now, we will show that $\mathbf{F}_k$ satisfies (\ref{dash-rthosvd:app:eqn7}). For clarity, let
    \begin{align*}
        &\widehat{\mathbf{\Sigma}}_{k,1}=\prod_{t=1}^q\left(\mathbf{\Sigma}_{k,1}^2-\alpha_t^{(k)}\mathbf{I}_{r_k}\right),\
      \widehat{\mathbf{\Sigma}}_{k,2}=\prod_{t=1}^q\left(\mathbf{\Sigma}_{k,3}^2-\alpha_t^{(k)}\mathbf{I}_{\widehat{n}_k-r_k}\right).
    \end{align*}
    Note that we have
    \begin{align*}
        \prod_{t=1}^q&\left(\mathbf{A}_{(k)}\mathbf{A}_{(k)}^\top-\alpha_t^{(k)}\mathbf{I}_{n_k}\right)\mathbf{A}_{(k)}\mathbf{\Omega}_k\mathbf{F}_k-\mathbf{A}_{(k)}\\
        &=\mathbf{Q}_k\prod_{t=1}^q\left(\mathbf{\Sigma}_k^2-\alpha_t^{(k)}\mathbf{I}_{\widehat{n}_k}\right)\mathbf{\Sigma}_k
        \begin{pmatrix}
             \mathbf{H}_k\\
             \mathbf{\Phi}_k
        \end{pmatrix}
        \begin{pmatrix}
            \mathbf{H}_k^\dag\widehat{\mathbf{\Sigma}}_{k,1}^{-1}&
            \mathbf{0}_1
        \end{pmatrix}
        \mathbf{P}_k^\top-\mathbf{Q}_k\mathbf{\Sigma}_k\mathbf{P}_k^\top\\
        &=\mathbf{Q}_k\mathbf{\Sigma}_k\prod_{t=1}^q\left(\mathbf{\Sigma}_k^2-\alpha_t^{(k)}\mathbf{I}_{\widehat{n}_k}\right)
        \begin{pmatrix}
             \mathbf{H}_k\\
             \mathbf{\Phi}_k
        \end{pmatrix}
        \begin{pmatrix}
            \mathbf{H}_k^\dag\widehat{\mathbf{\Sigma}}_{k,1}^{-1}&
            \mathbf{0}_1
        \end{pmatrix}
        \mathbf{P}_k^\top-\mathbf{Q}_k\mathbf{\Sigma}_k\mathbf{P}_k^\top\\
        &=\mathbf{Q}_k\mathbf{\Sigma}_k\left(
        \begin{pmatrix}
             \mathbf{I}_{j_k} & \mathbf{0}_{r_k\times (\widehat{n}_k-r_k)} \\
             \widehat{\mathbf{\Sigma}}_{k,2}\mathbf{\Phi}_k\mathbf{H}_k^\dag\widehat{\mathbf{\Sigma}}_{k,1}^{-1}
              & \mathbf{0}_{(\widehat{n}_k-r_k)\times (\widehat{n}_k-r_k)}
        \end{pmatrix}
        -\mathbf{I}_{\widehat{n}_k}\right)\mathbf{P}_k^\top\\
        &=\mathbf{Q}_k
        \begin{pmatrix}
             \mathbf{0}_{j_k\times j_k} & \mathbf{0}_{r_k\times (\widehat{n}_k-r_k)} \\
             \mathbf{\Sigma}_{k,2}\widehat{\mathbf{\Sigma}}_{k,2}\mathbf{\Phi}_k\mathbf{H}_k^\dag\widehat{\mathbf{\Sigma}}_{k,1}^{-1} & -\mathbf{\Sigma}_{k,2}
        \end{pmatrix}
        \mathbf{P}_k^\top,
    \end{align*}
    which implies that
    \begin{align*}
       &\left\|\prod_{t=1}^q\left(\mathbf{A}_{(k)}\mathbf{A}_{(k)}^\top-\alpha_t^{(k)}\mathbf{I}_{n_k}\right)\mathbf{A}_{(k)}\mathbf{\Omega}_k\mathbf{F}_k-\mathbf{A}_{(k)}\right\| \\
       &\quad\quad \leq \|\mathbf{\Phi}_k\|_2\|\mathbf{H}_k^\dag\|_2\|\widehat{\mathbf{\Sigma}}_{k,2}\|_2\|\widehat{\mathbf{\Sigma}}_{k,1}^{-1}\|_2\|\mathbf{\Sigma}_{k,2}\| + \|\mathbf{\Sigma}_{k,2}\|.
    \end{align*}
Therefore,
    \begin{align*}
        \|\widehat{\mathbf{\Sigma}}_{k,1}^{-1}\|_2&
        =\prod_{t=1}^q\left(\sigma_{r_k}(\mathbf{A}_{(k)})^2-\alpha_t^{(k)}\right)^{-1},\
        \|\widehat{\mathbf{\Sigma}}_{k,2}\|_2
        =\prod_{t=1}^q\left(\sigma_{r_k+1}(\mathbf{A}_{(k)})^2-\alpha_t^{(k)}\right).
    \end{align*}
    According to Lemma \ref{dash-rthosvd:lem1}, for any $\gamma_k>1$, one has
    \begin{align*}
        \|\mathbf{\Phi}_k\|_2\leq\|\mathbf{\Omega}_k\|_2\leq \sqrt{2\widehat{n}_k}\gamma_k
    \end{align*}
    with a probability at least
    \begin{equation*}
	1-\frac{1}{4(\gamma_k^2-
        1)\sqrt{\pi\widehat{n}_k\gamma_k^2}}\left(\frac{2\gamma_k^2}{e^{\gamma_k^2-1}}\right)^{\widehat{n}_k}\geq 0.
    \end{equation*}
    Therefore,
    \begin{align*}
        &\left\|\prod_{t=1}^q\left(\mathbf{A}_{(k)}\mathbf{A}_{(k)}^\top-\alpha_t^{(k)}\mathbf{I}_{n_k}\right)\mathbf{A}_{(k)}\mathbf{\Omega}_k\mathbf{F}_k-\mathbf{A}_{(k)}\right\|_F\\
        &\leq\left(\sqrt{2\widehat{n}_k}\gamma_k
        \frac{\prod_{t=1}^q\left(\sigma_{r_k+1}(\mathbf{A}_{(k)})^2-\alpha_t^{(k)}\right)}{\prod_{t=1}^q\left(\sigma_{r_k}(\mathbf{A}_{(k)})^2-\alpha_t^{(k)}\right)}+1\right)\sqrt{\sum_{i_k=r_k+1}^{\widehat{n}_k}\sigma_{i_k}(\mathbf{A}_{(k)})^2}
    \end{align*}
    with a probability at least the amount in (\ref{dash-rthosvd:app:eqn6}).
\end{proof}
We now give a rigorous proof for Theorem \ref{dash-rthosvd:thm1}.
\begin{proof}
    According to (\ref{dash-rthosvd:thosvd-general-two}) and Lemma \ref{dash-rthosvd:lem5}, one has
    \begin{align*}
        &\|\mathcal{G}\times_1\mathbf{U}_1\dots\times_d\mathbf{U}_d-\mathcal{A}\|
        \leq2\cdot\sum_{k=1}^d\left\|\prod_{t=1}^q\left(\mathbf{A}_{(k)}\mathbf{A}_{(k)}^\top-\alpha_t^{(k)}\mathbf{I}_{n_k}\right)\mathbf{A}_{(k)}\mathbf{\Omega}_k\mathbf{F}_k-\mathbf{A}_{(k)}\right\|\\
        &\quad\quad+2\cdot\sum_{k=1}^d\|\mathbf{F}_k\|_2\left\|\mathbf{U}_k\mathbf{R}_k-\prod_{t=1}^q\left(\mathbf{A}_{(k)}\mathbf{A}_{(k)}^\top-\alpha_t^{(k)}\mathbf{I}_{n_k}\right)\mathbf{A}_{(k)}\mathbf{\Omega}_k\right\|,
    \end{align*}
    and noting that having two positive numbers $a,b\in\mathbb{R}$ implies $\sqrt{a^2+b^2}\leq a+b$,, it follows from Theorems \ref{dash-rthosvd:app:thm1} and \ref{dash-rthosvd:app:thm2} that
    \begin{align*}
        \|\mathcal{G}\times_1\mathbf{U}_1\dots\times_d\mathbf{U}_d-\mathcal{A}\|\leq \sum_{k=1}^d(f_k+g_k)\sqrt{\sum_{i_k=r_k+1}^{\widehat{n}_k}\sigma_{i_k}(\mathbf{A}_{(k)})^2}
    \end{align*}
    with a probability at least $1-\sum_{k=1}^d\Phi_k$, where for each $k$, $\Phi_k$ is given in (\ref{dash-rthosvd:main-probability}), and the expressions of $f_k$ and $g_k$ are given by
    \begin{equation*}
    f_k=2\cdot\left(\sqrt{2\widehat{n}_k}\gamma_k\frac{\prod_{t=1}^q\left(\sigma_{r_k+1}(\mathbf{A}_{(k)})^2-\alpha_t^{(k)}\right)}{\prod_{t=1}^q\left(\sigma_{r_k}(\mathbf{A}_{(k)})^2-\alpha_t^{(k)}\right)}+1\right)
    \end{equation*}
    and
    \begin{equation*}
    g_k=2\cdot\sqrt{2n_1\dots n_{k-1}n_{k+1}\dots n_dl_k}\beta_k\gamma_k\frac{\prod_{t=1}^q(\sigma_{r_k+1}(\mathbf{A}_{(k)})^2-\alpha_t^{(k)})}{\prod_{t=1}^q\left(\sigma_{r_k}(\mathbf{A}_{(k)})^2-\alpha_{t}^{(k)}\right)}.
    \end{equation*}
    
    Hence, the proof of Theorem \ref{dash-rthosvd:thm1} is complete.
\end{proof}

\subsection{Proof of Theorem \ref{dash-rthosvd:thm2}}
\label{dash-rthosvd:sec4:sub3}
We now consider the term $\|\mathcal{B}_k\times_k(\mathbf{U}_k\mathbf{U}_k^\top)-\mathcal{B}_k\|^2$ with $k=1,2,\cdots,d$ in the right-hand side of (\ref{dash-rthosvd:sthosvd-general-two}), where $\mathcal{B}_1=\mathcal{A}$ and $\mathcal{B}_k=\mathcal{A}\times_{1}\mathbf{U}_{1}^\top\dots\times_{k-1}\mathbf{U}_{k-1}^\top$ with $k=2,3,\dots,d$. For a given $k$, when the matrix $\mathbf{U}_k\in\mathbb{R}^{n_k\times r_k}$ is obtained by Algorithm \ref{dash-rthosvd:alg2:v2}, then the matrix $\mathbf{U}_k\mathbf{U}_k^\top\mathbf{A}_{(k)}$ is a good approximation to the matrix $\mathbf{B}_{k}$, where $\mathbf{B}_k$ is the mode-$k$ unfolding matrix of $\mathcal{B}_k$, provided that there exists matrices $\mathbf{\Omega}_k\in\mathbb{R}^{r_1\dots r_{k-1}n_{k+1}\dots n_d\times l_k}$ and $\mathbf{R}_k\in\mathbb{R}^{r_k\times l_k}$ such that:
\begin{enumerate}
    \item[(a)] $\mathbf{U}_k$ is orthonormal;
    \item[(b)] $\mathbf{U}_k\mathbf{R}_k$ is a good approximation of $\prod_{t=1}^q(\mathbf{B}_{k}\mathbf{B}_{k}^\top-\alpha_t^{(k)}\mathbf{I}_{n_k})\mathbf{B}_{k}\mathbf{\Omega}_k$;
    \item[(c)] there exists a matrix $\mathbf{F}_k\in\mathbb{R}^{l_k\times r_1\dots r_{k-1}n_{k+1}\dots n_d}$ such that $\prod_{t=1}^q(\mathbf{B}_{k}\mathbf{B}_{k}^\top-\alpha_t^{(k)}\mathbf{I}_{n_k})\mathbf{B}_{k}\mathbf{\Omega}_k\mathbf{F}_k$ is a good approximation of $\mathbf{A}_{(k)}$ and  $\|\mathbf{F}_k\|_2$ is not too large,
\end{enumerate}
where for each $t$, $\alpha_t^{(k)}$ denotes the shift value in the $t$-th shifted power iteration.
\begin{lemma}
\label{dash-rthosvd:lem5:v1}
    For a given $k$, suppose that $\mathbf{B}_{k}\in\mathbb{R}^{n_k\times r_1\dots r_{k-1}n_{k+1}\dots n_d}$,
    $\mathbf{U}_k\in\mathbb{R}^{n_k\times r_k}$ is orthonormal, $\mathbf{F}_k\in\mathbb{R}^{l_k\times r_1\dots r_{k-1}n_{k+1}\dots n_d}$, $\mathbf{\Omega}_k\in\mathbb{R}^{r_1\dots r_{k-1}n_{k+1}\dots n_d\times l_k}$ with $r_k<l_k=r_k+s_k\leq \widetilde{n}_k$ and  $\mathbf{R}_k\in\mathbb{R}^{r_k\times l_k}$. For a given positive integer $q>0$, we have
    \begin{align}\label{dash-rthosvd:app:eqn4:v1}
        \|\mathcal{B}_k\times_k(\mathbf{U}_k\mathbf{U}_k^\top)-\mathcal{B}_k\|^2
	      &\leq 2\left\|\prod_{t=1}^q\left(\mathbf{B}_{k}\mathbf{B}_{k}^\top-\alpha_t^{(k)}
        \mathbf{I}_{n_k}\right)\mathbf{B}_{k}\mathbf{\Omega}_k\mathbf{F}_k-\mathbf{B}_{k}\right\|^2\nonumber\\
          &+2\|\mathbf{F}_k\|_2^2\left\|\mathbf{U}_k\mathbf{R}_k-\prod_{t=1}^q\left(\mathbf{B}_{k}\mathbf{B}_{k}^\top
          -\alpha_t^{(k)}\mathbf{I}_{n_k}\right)\mathbf{B}_{k}\mathbf{\Omega}_k\right\|^2,
    \end{align}
    where for each $t$, $\alpha_t^{(k)}$ denotes the shift value in the $t$-th shifted power iteration.
\end{lemma}

By combining Lemmas \ref{dash-rthosvd:lem1}, \ref{dash-rthosvd:lem3} and \ref{dash-rthosvd:lem4}, the upper bound for the second part in the right-hand side of (\ref{dash-rthosvd:app:eqn4:v1}) is given in the following theorem.
\begin{theorem}
\label{dash-rthosvd:app:thm1:v1}
    For a given $k$, suppose that
    $\mathbf{\Omega}_k\in\mathbb{R}^{r_1\dots r_{k-1}n_{k+1}\dots n_d\times (r_k+p_k)}$ is any standard Gaussian matrix with $l_k=r_k+p_k<\widetilde{n}_k$ and $\mathbf{B}_{k}\in\mathbb{R}^{n_k\times r_1\dots r_{k-1}n_{k+1}\dots n_d}$ is the mode-$k$ unfolding of $\mathcal{B}_k\in\mathbb{R}^{r_1\times\cdots r_{k-1}\times n_{k} \cdots\times n_d}$. For $\gamma_k>1$, let
    \begin{align}\label{dash-rthosvd:app:eqn5:v1}
        1-\frac{1}{4(\gamma_k^2-1)\sqrt{\pi r_1\dots r_{k-1}n_{k+1}\dots n_d\gamma_k^2}}\left(\frac{2\gamma_k^2}{e^{\gamma_k^2}-1}\right)^{n_1\dots r_{k-1}r_{k+1}\dots n_d}>0.
    \end{align}
    For a given positive integer $q>0$, there exists an orthonormal matrix $\mathbf{U}_k\in\mathbb{R}^{n_k\times r_k}$ and $\mathbf{R}_k\in\mathbb{R}^{r_k\times l_k}$ such that
    \begin{align*}
        &\left\|\mathbf{U}_k\mathbf{R}_k-\prod_{t=1}^q\left(\mathbf{B}_{k}\mathbf{B}_{k}^\top-\alpha_t^{(k)}\mathbf{I}_{n_k}\right)\mathbf{B}_{k}\mathbf{\Omega}_k\right\|\\
        &\leq \sqrt{2r_1\dots r_{k-1}n_{k+1}\dots n_d}\gamma_k\prod_{t=1}^q\left(\sigma_{r_k+1}(\mathbf{B}_{k})^2-\alpha_t^{(k)}\right)\sqrt{\sum_{i_k=r_k+1}^{\widetilde{n}_k}\sigma_{i_k}(\mathbf{B}_{k})^2}
    \end{align*}
    with a probability at least the amount in (\ref{dash-rthosvd:app:eqn5:v1}), where for each $t$, $\alpha_t^{(k)}$ denotes the shift value in the $t$-th shifted power iteration.
\end{theorem}

For each $k$, we now consider the upper bounds for $\|\mathbf{F}_k\|_2$ and the first term in the right-hand side of (\ref{dash-rthosvd:app:eqn4:v1}), which are summarized in the following theorem.
\begin{theorem}
\label{dash-rthosvd:app:thm2:v1}
    For a given $k$, suppose that
    $\mathbf{\Omega}_k\in\mathbb{R}^{r_1\dots r_{k-1}n_{k+1}\dots n_d\times l_k}$ is any standard Gaussian matrix with $r_k<l_k=r_k+p_k<\widetilde{n}_k$, and $\mathbf{B}_{k}\in\mathbb{R}^{n_k\times r_1\dots r_{k-1}n_{k+1}\dots n_d}$ is the mode-$k$ unfolding of $\mathcal{B}_k\in\mathbb{R}^{r_1\times\cdots r_{k-1}\times n_{k} \cdots\times n_d}$. For $\beta_k>0$ and $\gamma_k>1$, let
    \begin{align}\label{dash-rthosvd:app:eqn6:v1}
        1&-\frac{1}{\sqrt{2\pi(l_k-r_k+1)}}\left(\frac{e}{(l_k-r_k+1)\beta_k}\right)^{l_k-r_k+1}\nonumber\\
        &-\frac{1}{4(\gamma_k^2-
        1)\sqrt{\pi\widetilde{n}_k\gamma_k^2}}\left(\frac{2\gamma_k^2}{e^{\gamma_k^2-1}}\right)^{\widetilde{n}_k}
    \end{align}
    be nonnegative. For a given positive integer $q>0$, there exists an orthonormal matrix $\mathbf{U}_k\in\mathbb{R}^{n_k\times r_k}$ and $\mathbf{R}_k\in\mathbb{R}^{r_k\times (r_k+p_k)}$ such that
    \begin{align}\label{dash-rthosvd:app:eqn7:v1}
        &\left\|\prod_{t=1}^q\left(\mathbf{B}_{k}\mathbf{B}_{k}^\top-\alpha_t^{(k)}\mathbf{I}_{n_k}\right)\mathbf{B}_{k}\mathbf{\Omega}_k\mathbf{F}_k-\mathbf{B}_{k}\right\|\nonumber\\
        &\leq\left(\sqrt{2\widetilde{n}_k}\gamma_k
        \frac{\prod_{t=1}^q\left(\sigma_{r_k+1}(\mathbf{B}_{k})^2-\alpha_t^{(k)}\right)}{\prod_{t=1}^q\left(\sigma_{r_k}(\mathbf{B}_{k})^2-\alpha_t^{(k)}\right)}+1\right)\sqrt{\sum_{i_k=r_k+1}^{\widetilde{n}_k}\sigma_{i_k}(\mathbf{B}_{k})^2}
    \end{align}
    and
    \begin{align}\label{dash-rthosvd:app:eqn8:v1}
        \|\mathbf{F}_k\|_2\leq \frac{\sqrt{l_k}\beta_k}{\prod_{t=1}^q\left(\sigma_{r_k}(\mathbf{B}_{k})^2-\alpha_{t}^{(k)}\right)}
    \end{align}
    with a probability at least the amount in (\ref{dash-rthosvd:app:eqn6:v1}), where for each $t$, $\alpha_t^{(k)}$ denotes the shift value in the $t$-th shifted power iteration.
\end{theorem}
\begin{remark} 
The proofs for Lemma \ref{dash-rthosvd:lem5:v1} and Theorem \ref{dash-rthosvd:app:thm2:v1} are similar to that for Lemma \ref{dash-rthosvd:lem5} and Theorem \ref{dash-rthosvd:app:thm2}, respectively. Hence, we omit them in this paper.
\end{remark}

Then, by combining Lemma \ref{dash-rthosvd:lem5:v1}, and Theorems \ref{dash-rthosvd:app:thm1:v1} and \ref{dash-rthosvd:app:thm2:v1}, the following inequality
\begin{align*}
    &\|\mathcal{G}\times_1\mathbf{U}_1\dots\times_d\mathbf{U}_d-\mathcal{A}\|\\
    &\leq2\cdot\sum_{k=1}^d\left(\sqrt{2\widetilde{n}_k}\gamma_k
    \frac{\prod_{t=1}^q\left(\sigma_{r_k+1}(\mathbf{B}_k)^2-\alpha_t^{(k)}\right)}{\prod_{t=1}^q\left(\sigma_{r_k}(\mathbf{B}_k)^2-\alpha_t^{(k)}\right)}+1\right)\sqrt{\sum_{i_k=r_k+1}^{\widetilde{n}_k}\sigma_{i_k}(\mathbf{B}_k)^2}\\
    &+2\cdot\sum_{k=1}^d\sqrt{2r_1\dots r_{k-1}n_{k+1}\dots n_dl_k}\beta_k\gamma_k\frac{\prod_{t=1}^q\left(\sigma_{r_k+1}(\mathbf{B}_k)^2-\alpha_t^{(k)}\right)}{\prod_{t=1}^q\left(\sigma_{r_k}(\mathbf{B}_k)^2-\alpha_{t}^{(k)}\right)}\sqrt{\sum_{i_k=r_k+1}^{\widetilde{n}_k}\sigma_{i_k}(\mathbf{B}_k)^2}
\end{align*}
holds with a probability at least $1-\sum_{k=1}^d\Psi_k$, where there exist two real numbers $\beta_k>0$ and $\gamma_k>1$ such that $0<\sum_{k=1}^d\Psi_k<1$, and for each $k$, $\Psi_k$ is given in 
\begin{align*}
    \Psi_k=&\frac{1}{\sqrt{2\pi(l_k-r_k+1)}}\left(\frac{e}{(l_k-r_k+1)\beta_k}\right)^{l_k-r_k+1}+\frac{1}{4(\gamma_k^2-
    1)\sqrt{\pi\widetilde{n}_k\gamma_k^2}}\left(\frac{2\gamma_k^2}{e^{\gamma_k^2-1}}\right)^{\widetilde{n}_k}\nonumber\\
    &+\frac{1}{4(\gamma_k^2-1)\sqrt{\pi r_1\dots r_{k-1}n_{k+1}\dots n_d\gamma_k^2}}\left(\frac{2\gamma_k^2}{e^{\gamma_k^2}-1}\right)^{r_1\dots r_{k-1}n_{k+1}\dots n_d}.
\end{align*}

\begin{remark}
    Note that for each $k$, one has $\mathbf{B}_k=\mathbf{A}_{(k)}(\mathbf{U}_1\otimes \dots \otimes \mathbf{U}_{k-1}\otimes \mathbf{I}_{n_{k+1}}\otimes\dots \otimes \mathbf{I}_{n_{d}})$. It follows from \cite[Lemma 9]{che2025efficient-acom} that $\sigma_{i_k}(\mathbf{B}_k)\leq \sigma_{i_k}(\mathbf{A}_{(k)})$ with $i_k=1,2,\dots,\min\{\widehat{n}_k,\widetilde{n}_k\}$, and \begin{align*}
       \sum_{i_k=r_k+1}^{\widetilde{n}_k}\sigma_{i_k}(\mathbf{B}_k)^2\leq \sum_{i_k=r_k+1}^{\widehat{n}_k}\sigma_{i_k}({\color{red}\mathbf{A}_{(k)}})^2.
    \end{align*}
    Combining the choice of $\alpha_t^{(k)}$ in Algorithm \ref{dash-rthosvd:alg2:v2}, one has $\sigma_{i_k}(\mathbf{B}_k)^2-\alpha_t^{(k)})\leq \sigma_{i_k}(\mathbf{A}_{(k)})^2-\alpha_t^{(k)})$. For each $k$, we cannot ensure that
    \begin{align*}
       \frac{\prod_{t=1}^q\left(\sigma_{r_{k+1}}(\mathbf{B}_k)^2-\alpha_t^{(k)}\right)}{\prod_{t=1}^q\left(\sigma_{r_k}(\mathbf{B}_k)^2-\alpha_{t}^{(k)}\right)}\leq \frac{\prod_{t=1}^q\left(\sigma_{r_{k+1}}(\mathbf{A}_{(k)})^2-\alpha_t^{(k)}\right)}{\prod_{t=1}^q\left(\sigma_{r_k}(\mathbf{A}_{(k)})^2-\alpha_{t}^{(k)}\right)}.
    \end{align*}

However, it follows from $\sigma_{r_{k+1}}(\mathbf{B}_k)^2-\alpha_t^{(k)})\leq \sigma_{r_k}(\mathbf{B}_{(k)})^2-\alpha_t^{(k)})$ that
{\color{red}
\begin{align*}
   \frac{\prod_{t=1}^q\left(\sigma_{r_{k+1}}(\mathbf{B}_k)^2-\alpha_t^{(k)}\right)}{\prod_{t=1}^q\left(\sigma_{r_k}(\mathbf{B}_k)^2-\alpha_{t}^{(k)}\right)}\leq 1.
\end{align*}
}
\end{remark}
\begin{remark}
Let $g(x)=cx^{-1/2}a^x$ for $x> 1$, with
\begin{equation*}
 c=\frac{1}{4(\gamma_k^2-1)\sqrt{\pi\gamma_k^2}},\quad a=\frac{2\gamma_k^2}{e^{\gamma_k^2-1}}.
 \end{equation*}
 By calculation, the differentiation of $g(x)$ is
 \begin{equation*}
 g'(x)=-\frac{c}{2}x^{-3/2}a^x+cx^{-1/2}a^x=cx^{-1/2}a^x\frac{x-1}{2x}>0.
 \end{equation*}
 Hence, $g(x)$ is strictly monotonically increasing with respect to $x$, which implies that $\Psi_k\leq \Phi_k$ with $\widetilde{n}_k\leq \widehat{n}_k$.
\end{remark}

Hence, the proof of Theorem \ref{dash-rthosvd:thm2} is complete according the above descriptions.
\section{Numerical examples}
\label{dash-rthosvd:sec5}

We use MATLAB R2023b and MATLAB Tensor Toolbox (see \cite{bader2006algorithm,bader2023matlab}) to implement {\color{red}Algorithms \ref{dash-rthosvd:alg2} and \ref{dash-rthosvd:alg2:v2}}. We compare the proposed algorithms with state-of-the-art methods for fixed-rank Tucker approximation, including deterministic and randomized variants. The numerical experiments were conducted on a personal computer (PC) with an Intel Core i7-10700 CPU (2.90GHz) and 64GB RAM. For run time comparison, all methods are run under the same running settings to maintain fairness. The numerical results for each method are the average of ten repeated experiments to ensure statistical reliability.

{\color{red}For any $d$-tuple $\{r_1,\dots,r_d\}$, let $\{\mathcal{G};\mathbf{U}_1,\dots,\mathbf{U}_d\}$ be any solution for Problem \ref{dash-rthosvd:prob3}}, then the relative error (RE) is defined as
\begin{equation}
   {\rm RE} = \|\mathcal{A}-\widehat{\mathcal{A}}\|/\|\mathcal{A}\|
   \label{dash-rthosvd:eqn:re}
\end{equation}
with $\widehat{\mathcal{A}}=\mathcal{G}\times_{1}{\bf U}_{1}\times_{2}{\bf U}_{2}\dots\times_{d}{\bf U}_{d}$.

The state-of-the-art algorithms used to solve Problem \ref{dash-rthosvd:prob1} for comparison are listed as follows:
\begin{enumerate}
    \item[(1)] HOOI (denoted by Tucker-ALS and also called the alternating least squares): implemented by the function ``tucker\_als'' (see \cite{bader2023matlab});
    \item[(2)] several randomized variants of HOOI: Tucker-TS and Tucker-TTMTS\footnote{The MATLAB codes for Tucker-TS and Tucker-TTMTS can be found in https://github.com/OsmanMalik/tucker-tensorsketch.} (see \cite{malik2018low}), Sketch-Tucker-ALS (see \cite{ma2021fast}) and Randomized-Tucker-ALS (see \cite{ahmadi2021randomized});
    \item[(3)] T-HOSVD and ST-HOSVD: implemented by the function ``mlsvd''  (see \cite{tensorlab}) with different parameters (in particular, for the case of ST-HOSVD, a faster but possibly less accurate eigenvalue decomposition is used to compute the factor matrices);
    \item[(4)] the randomized variants of T-HOSVD and ST-HOSVD in \cite{minster2020randomized} (denoted by R-T-HOSVD and R-ST-HOSVD): implemented by the function ``mlsvd\_rsi'' (see \cite{tensorlab});
    \item[(5)] the randomized variants of T-HOSVD and ST-HOSVD in \cite{minster2024parallel} (denoted by rSTHOSVDkron and rHOSVDkronreuse);
    \item[(6)] the randomized variants of T-HOSVD and ST-HOSVD in \cite{sun2020low} (denoted by Tucker-Sketch-SP and Tucker-Sketch-TP);
    \item[(7)] the randomized variants of T-HOSVD and ST-HOSVD in \cite{ahmadi2021randomized} (denoted by RP-HOSVD).
\end{enumerate}

{\color{red}Without loss of generality, when finding the mode-$k$ factor matrices $\{\mathbf{U}_1,\cdots,\mathbf{U}_d\}$, the processing order in all randomized variants of ST-HOSVD  is set to $\{1,2,\dots,d\}$.}

\subsection{Several test tensors}
\label{dash-rthosvd:sec5:subsec1}
In this section, we introduce several testing tensors from synthetic and real-world data. The first is a sparse tensor $\mathcal{B}\in\mathbb{R}^{600\times 600\times 600}$ {\rm(cf. \cite{ahmadi2021randomized,minster2020randomized,saibaba2016hoid})}, which is given by
\begin{equation*}
    \mathcal{B}=\sum_{i=1}^{50}\frac{\gamma}{i}\mathbf{x}_i\circ\mathbf{y}_i\circ\mathbf{z}_i
    +\sum_{i=51}^{600}\frac{1}{i}\mathbf{x}_i\circ\mathbf{y}_i\circ\mathbf{z}_i,
\end{equation*}
where $\mathbf{x}_i\in \mathbb{R}^{600}$, $\mathbf{y}_i\in\mathbb{R}^{600}$ and $\mathbf{z}_i\in\mathbb{R}^{600}$ are threesparse vectors each with a sparsity ratio of 0.05, and the symbol ``$\circ$'' represents the vector outer product. Here we set $\gamma=1000$.

The second tensor $\mathcal{C}\in \mathbb{R}^{600\times 600\times 600}$ is constructed as
\begin{align*}
    \mathcal{C}={\rm tendiag}(\mathbf{v},[600,600,600])\times_1\mathbf{A}_1\times_2\mathbf{A}_2\times_3\mathbf{A}_3,
\end{align*}
where ``{\rm tendiag}'' is the function in the MATLAB Tensor Toolbox  \cite{bader2023matlab}, which converts $\mathbf{v}$ to a diagonal tensor of the size $600\times 600\times 600$, and $\mathbf{A}_n$ is an orthogonal basis for the column space of a standard Gaussian matrix $\mathbf{B}_n\in\mathbb{R}^{600\times 600}$. Three kinds of tensors were tested to representing different distribution patterns of the vector $\mathbf{v}\in\mathbb{R}^{600}$, where three kinds of $\mathbf{v}$ are given by
\begin{enumerate}
    \item[(a)] {\bf (Slow decay)} $v_i=1/i^2$ with $i=1,2,\dots,600$;
    \item[(b)] {\bf (Fast decay)} $v_i=\exp(-i/7)$ with $i=1,2,\dots,600$;
    \item[(c)] {\bf (S-shaped decay)} $v_i=0.001+(1+\exp(i-29))^{-1}$ with $i=1,2,\dots,600$.
\end{enumerate}

The extended Yale B database\footnote{The extended Yale Face Database B is available at \url{http://vision.ucsd.edu/~iskwak/ExtYaleDatabase/ExtYaleB.html}.} {\rm(see \cite{georghiades2001from})} contains $560$ images with each image containing $480\times 640$ pixels in the gray-scale range. This data is extracted into a tensor $\mathcal{A}_{\rm YaleB}\in\mathbb{R}^{480\times 640\times 560}$. The Columbia object image library COIL-100\footnote{COIL-100 can be downloaded at \url{https://www.cs.columbia.edu/CAVE/software/softlib/coil-100.php}.} consists of $7200$ color images that contain 100 objects under 72 different rotations. Each image has $128\times 128$ pixels in the RGB color range. We reshape these data to a third-order tensor $\mathcal{A}_{\rm COIL}$ of size $1152\times 1024\times 300$. In the Washington DC Mall database\footnote{The Washington DC Mall database is available at \url{https://engineering.purdue.edu/biehl/MultiSpec/hyperspectral.html}.}, the sensor system used in this case measures pixel response in 210 bands in the 0.4 to 2.4 $\mu$m region of the visible and infrared spectrum. Bands in the 0.9 and 1.4 $\mu$m region where the atmosphere is opaque have been omitted from the data set, leaving 191 bands remaining. The data set contains 1280 scan lines with 307 pixels in each scan line. This database is constructed as a tensor $\mathcal{A}_{\rm DCmall}$ of size $1280\times 307\times 191$.

\begin{figure}[htb]
    \setlength{\tabcolsep}{4pt}
    \renewcommand\arraystretch{1}
    \centering
    \subfigure[Algorithm \ref{dash-rthosvd:alg2}]{\includegraphics[width=4.5in,height=1.6in]{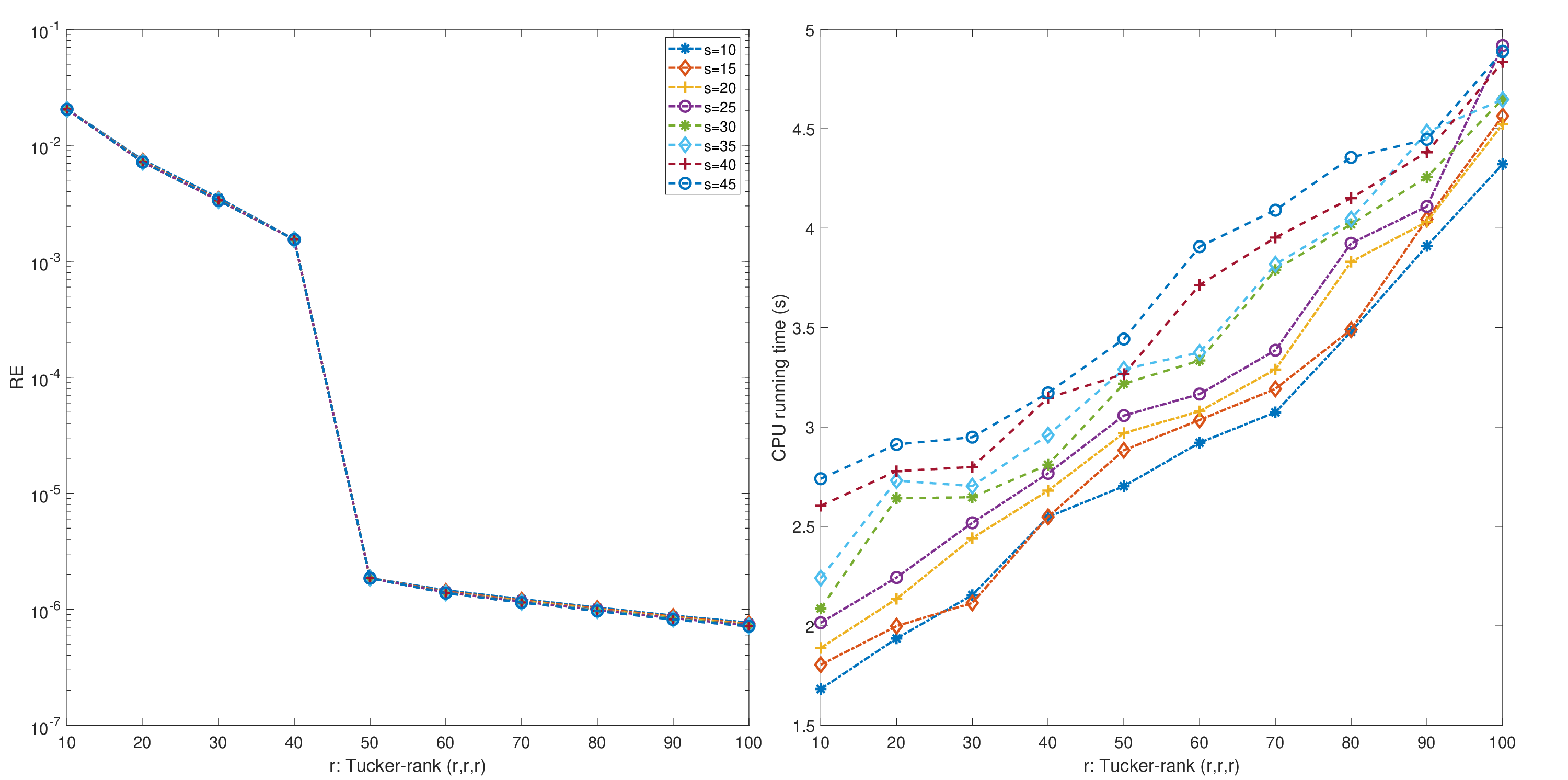}}\\
	\subfigure[Algorithm
\ref{dash-rthosvd:alg2:v2}]{\includegraphics[width=4.5in,height=1.6in]{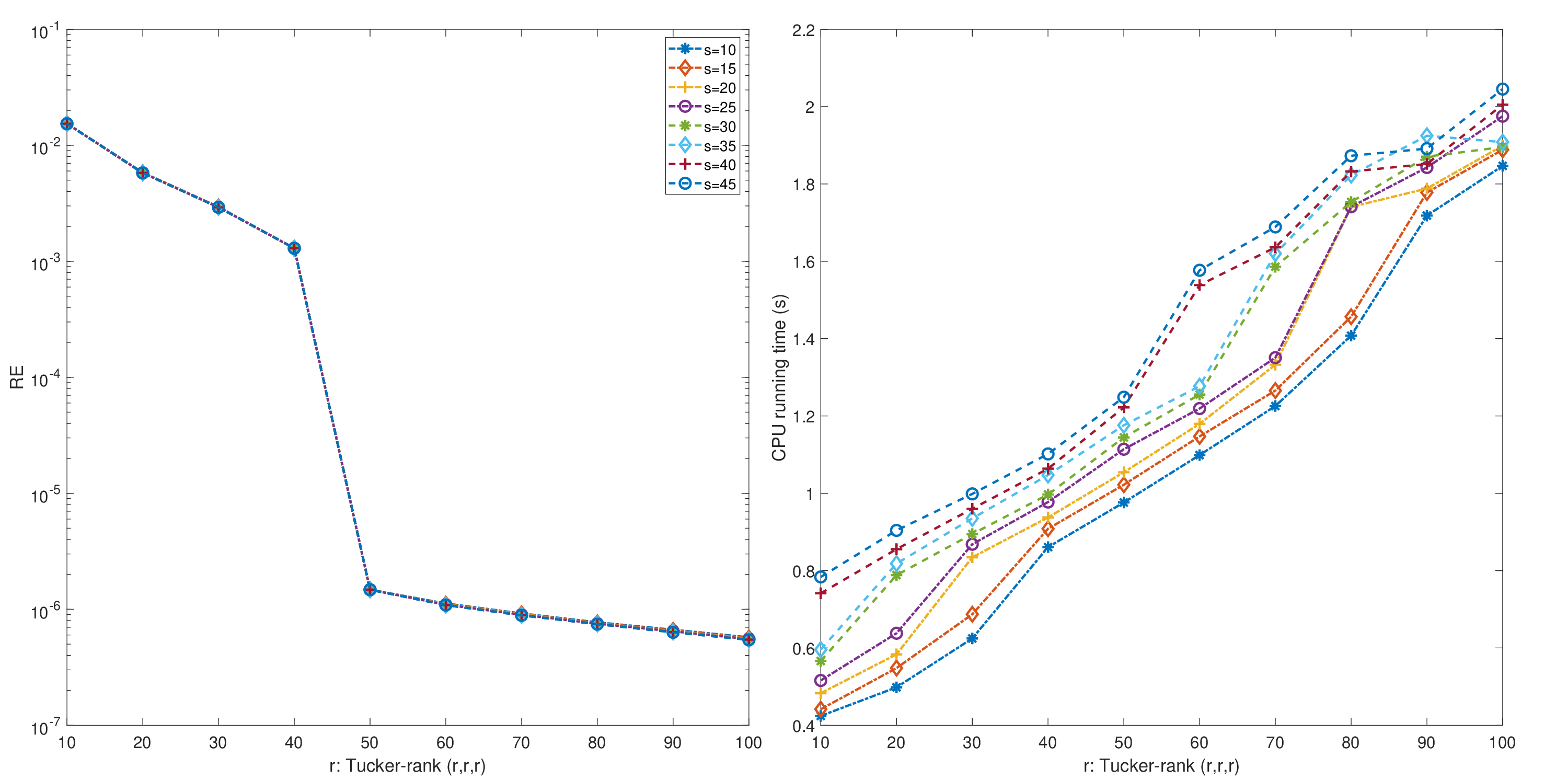}}\\
    \caption{For fixed $q=1$, numerical simulation results of Algorithms \ref{dash-rthosvd:alg2} and \ref{dash-rthosvd:alg2:v2}
    with different Tucker-ranks $(r,r,r)$ and $s$ on the test tensor $\mathcal{B}$.}\label{dash-rthosvd:figure1}
\end{figure}

\begin{figure}[htb]
    \setlength{\tabcolsep}{4pt}
    \renewcommand\arraystretch{1}
    \centering
    \subfigure[Algorithm \ref{dash-rthosvd:alg2}]{\includegraphics[width=4.5in,height=1.6in]{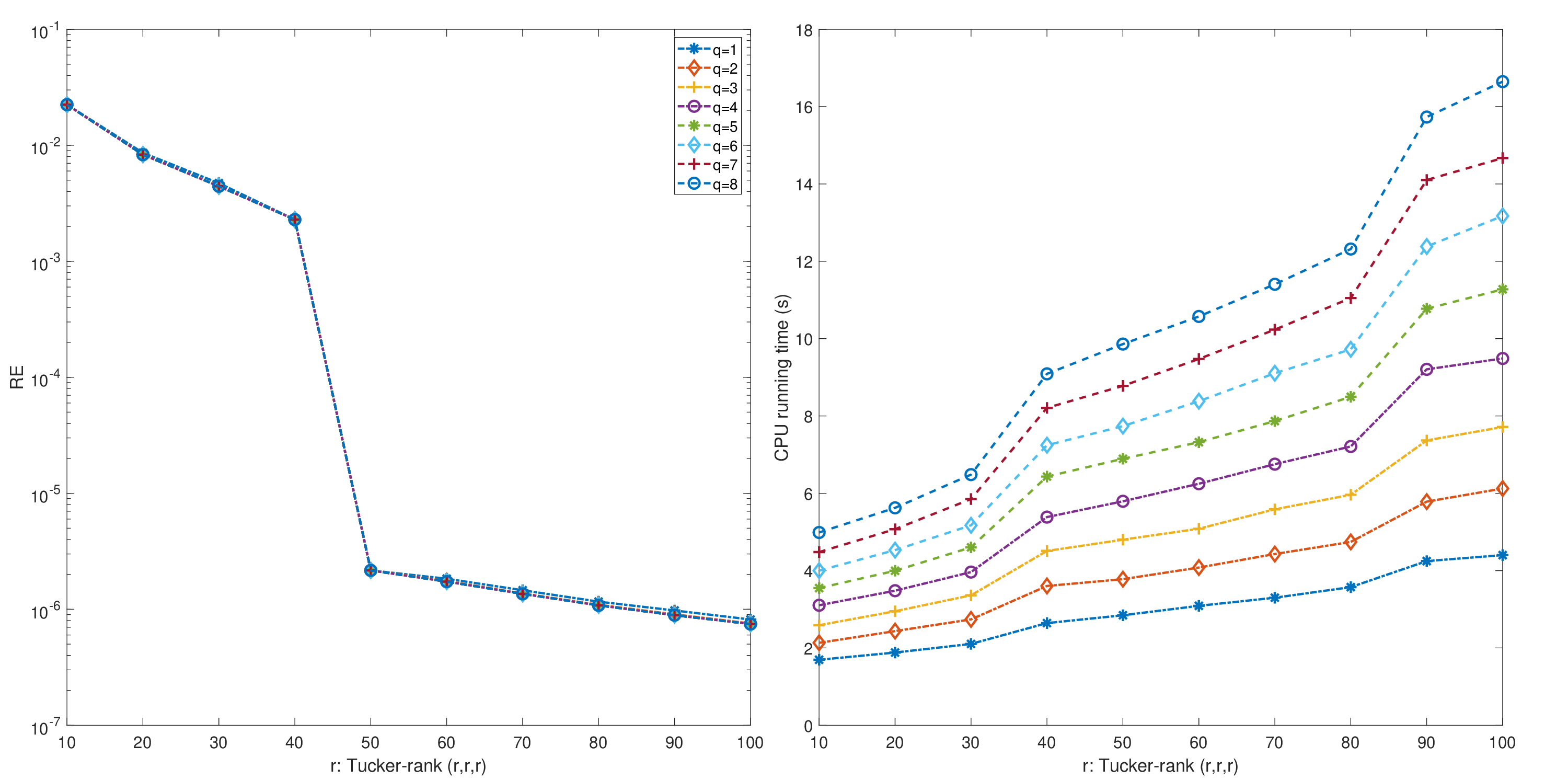}}\\
	\subfigure[Algorithm
\ref{dash-rthosvd:alg2:v2}]{\includegraphics[width=4.5in,height=1.6in]{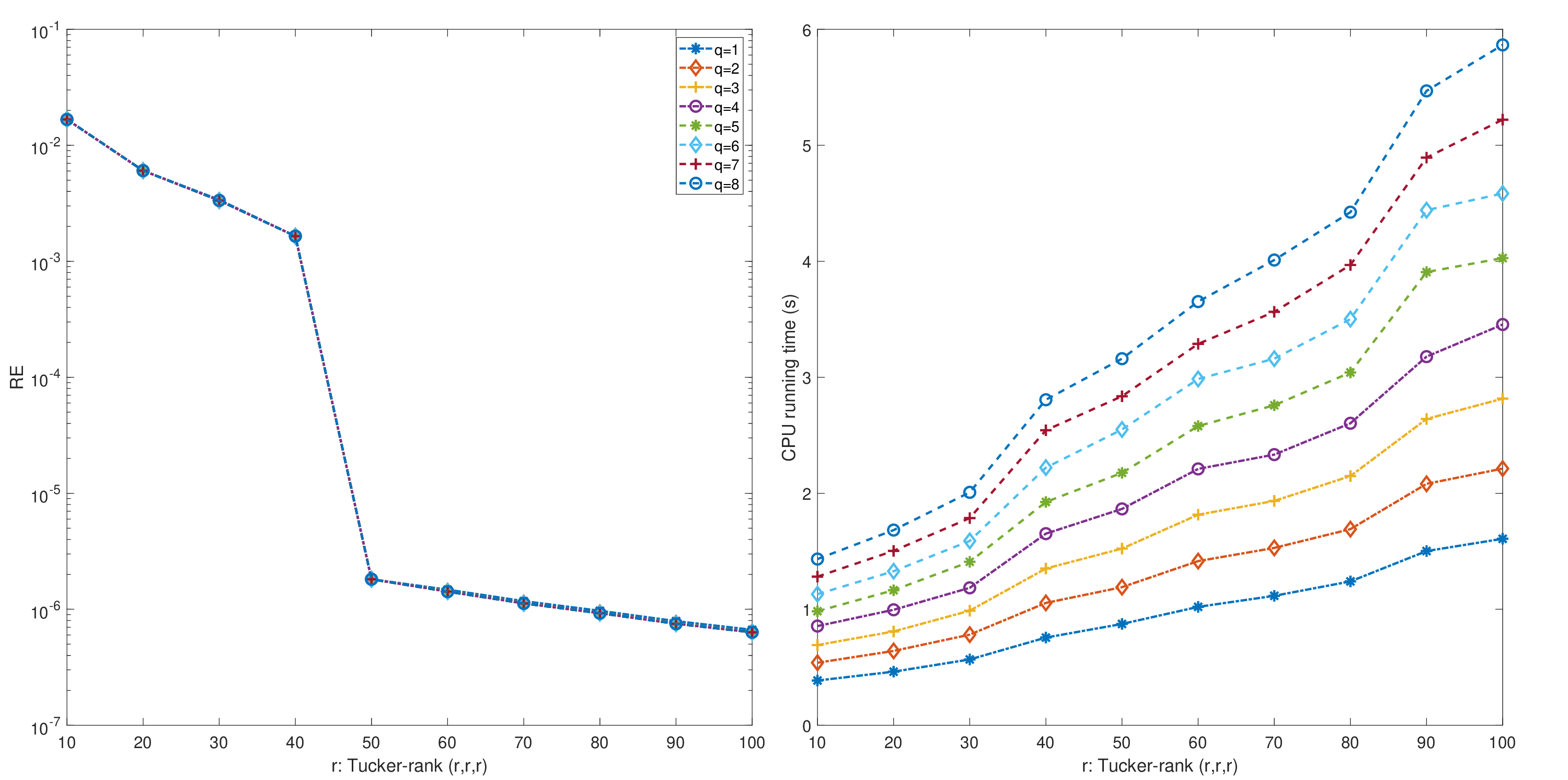}}\\
    \caption{For fixed $s=10$, numerical simulation results of Algorithms \ref{dash-rthosvd:alg2} and \ref{dash-rthosvd:alg2:v2}
    with different Tucker-ranks $(r,r,r)$ and $q$ on the test tensor $\mathcal{B}$.}\label{dash-rthosvd:figure2}
\end{figure}
\subsection{Parameter selection}
The desired Tucker-rank $\mathbf{r}=(r_1,r_2,\dots,r_d)$, the $d$-tuple of oversampling parameters $\mathbf{s}=(s_1,s_2,\dots,s_d)$ and the power parameter $q\geq 1$ are important inputs for Algorithms \ref{dash-rthosvd:alg2} and \ref{dash-rthosvd:alg2:v2}, which $\mathbf{s}$ and $q$ will affect the efficiency of these algorithms. For clarity, we assume that $s_k=s$ and $r_k=r$ with $k=1,2,\dots,d$.

With the same $r$, we compare RE and CPU run time of Algorithms \ref{dash-rthosvd:alg2} and \ref{dash-rthosvd:alg2:v2} with different $s$ and $q$. For fixed $q=1$, Figure \ref{dash-rthosvd:figure1} shows the comparison of RE and CPU run time obtained by applying Algorithms \ref{dash-rthosvd:alg2} and \ref{dash-rthosvd:alg2:v2} with different $r$ and $s$ values to the test tensor $\mathcal{B}$. For fixed $r=10$, the RE and CPU run time of Algorithms \ref{dash-rthosvd:alg2} and \ref{dash-rthosvd:alg2:v2} with different $r$ and $q$ values are illustrated in Figure \ref{dash-rthosvd:figure2}.

From these two figures, we conclude that: (a) with the same $r$, different choices of $s$ and $q$ are comparable in terms of RE; and (b) when fixing $q=1$ (or $s=10$), with the same $r$, CPU run time increases as $s$ (or $q$) increases. Hence, if there are no special requirements, we assume that $s=10$ and $q=1$.

\subsection{Comparison with different processing orders in Algorithm \ref{dash-rthosvd:alg2:v2}}

\begin{figure}
    \setlength{\tabcolsep}{4pt}
    \renewcommand\arraystretch{1}
    \centering
    \includegraphics[width=4.5in,height=1.6in]{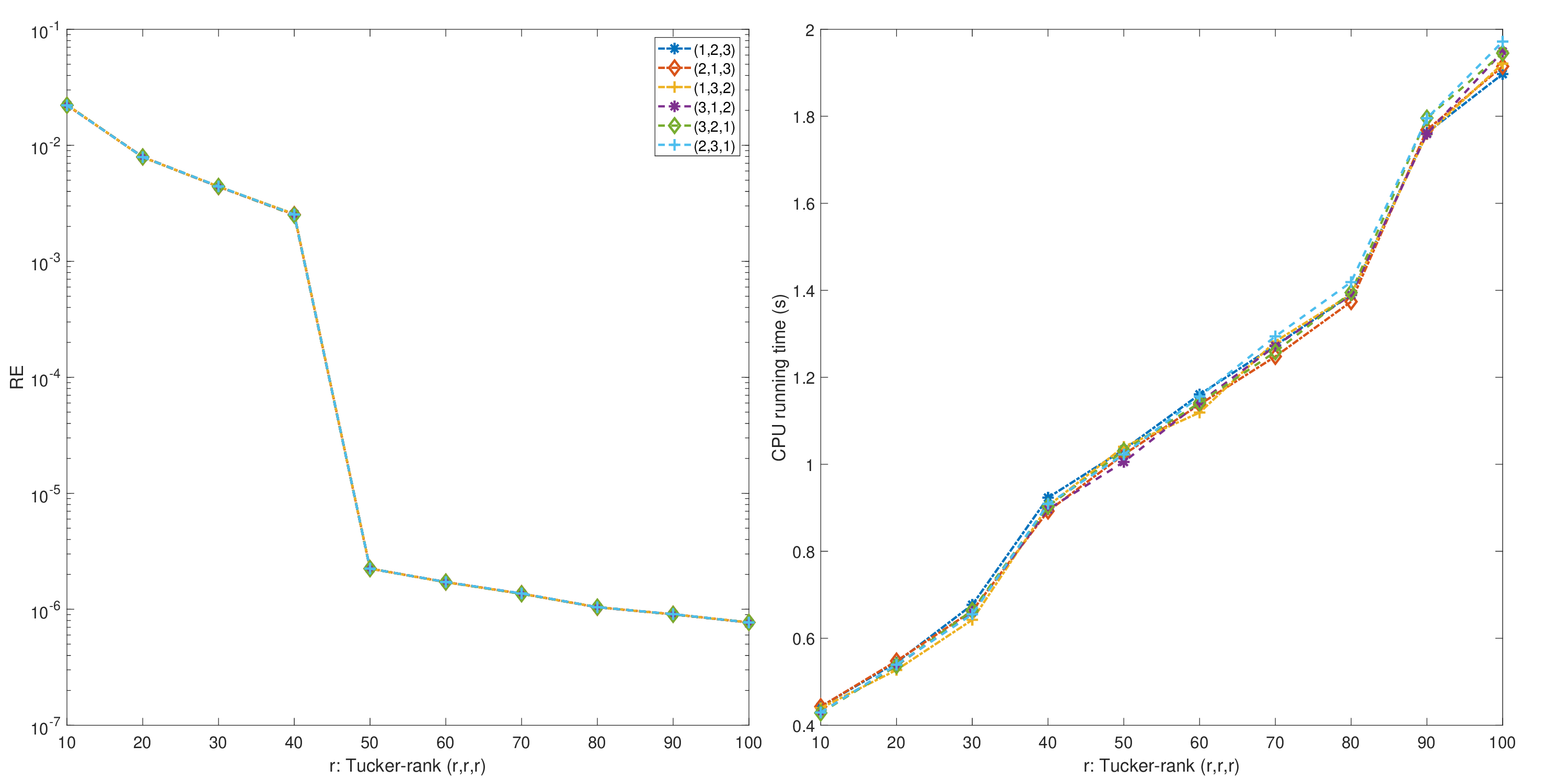}\\
    \caption{For fixed $s=10$ and $q=1$, numerical simulation results of Algorithm \ref{dash-rthosvd:alg2:v2}
    with different Tucker-ranks $(r,r,r)$ and processing orders on the test tensor $\mathcal{B}$.}\label{dash-rthosvd:figure3}
\end{figure}

The processing order in Algorithm \ref{dash-rthosvd:alg2:v2} was set to $(1,2,\dots,d)$. In general, any processing order in $\mathbb{S}_d$ is suitable for Algorithm \ref{dash-rthosvd:alg2:v2}. With the same $r$, we consider the efficiency of Algorithm \ref{dash-rthosvd:alg2:v2} with different processing orders via the test tensor $\mathcal{B}$. Figure \ref{dash-rthosvd:figure3} shows that the comparison of RE and CPU run time for Algorithm \ref{dash-rthosvd:alg2:v2}. The difference of the RE and CPU run time of Algorithm \ref{dash-rthosvd:alg2:v2} with different processing orders is not significant.
\subsection{Comparison with different algorithms}
\label{dash-rthosvd:sec5:subsec4}
In this section, we compare the efficiency of Algorithms \ref{dash-rthosvd:alg2} and \ref{dash-rthosvd:alg2:v2} with three commonly used algorithms (i.e., Tucker-ALS, T-HOSVD and ST-HOSVD) and their randomized variants via several test tensors in Section \ref{dash-rthosvd:sec5:subsec1}.

\begin{figure}
    \setlength{\tabcolsep}{4pt}
    \renewcommand\arraystretch{1}
    \centering
    \subfigure[The tensor $\mathcal{B}$]{\includegraphics[width=4.5in,height=1.6in]{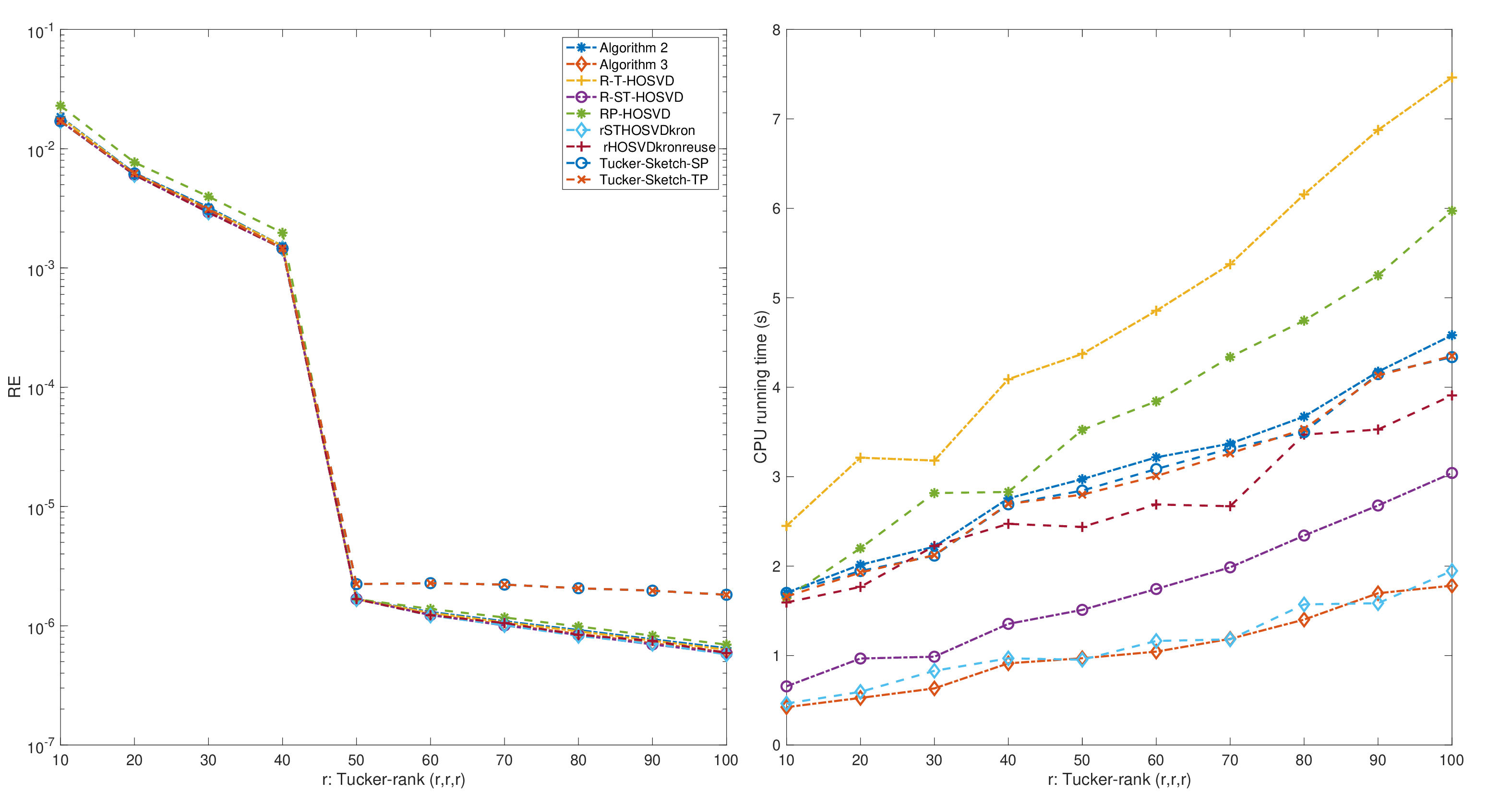}}\\
    \subfigure[The tensor $\mathcal{C}$ with Slow decay]{\includegraphics[width=4.5in,height=1.6in]{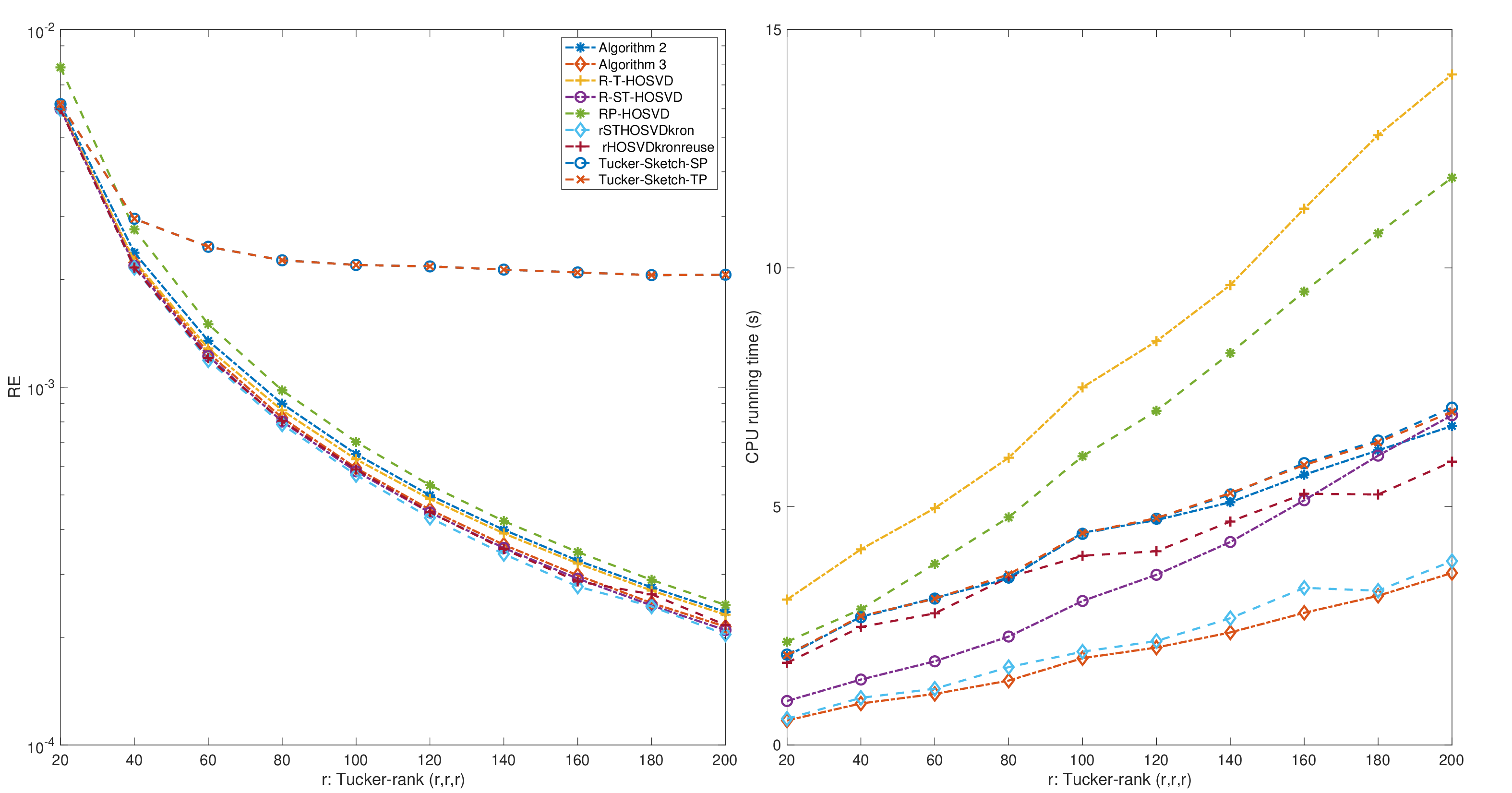}}\\
	\subfigure[The tensor $\mathcal{C}$ with Fast decay]{\includegraphics[width=4.5in,height=1.6in]{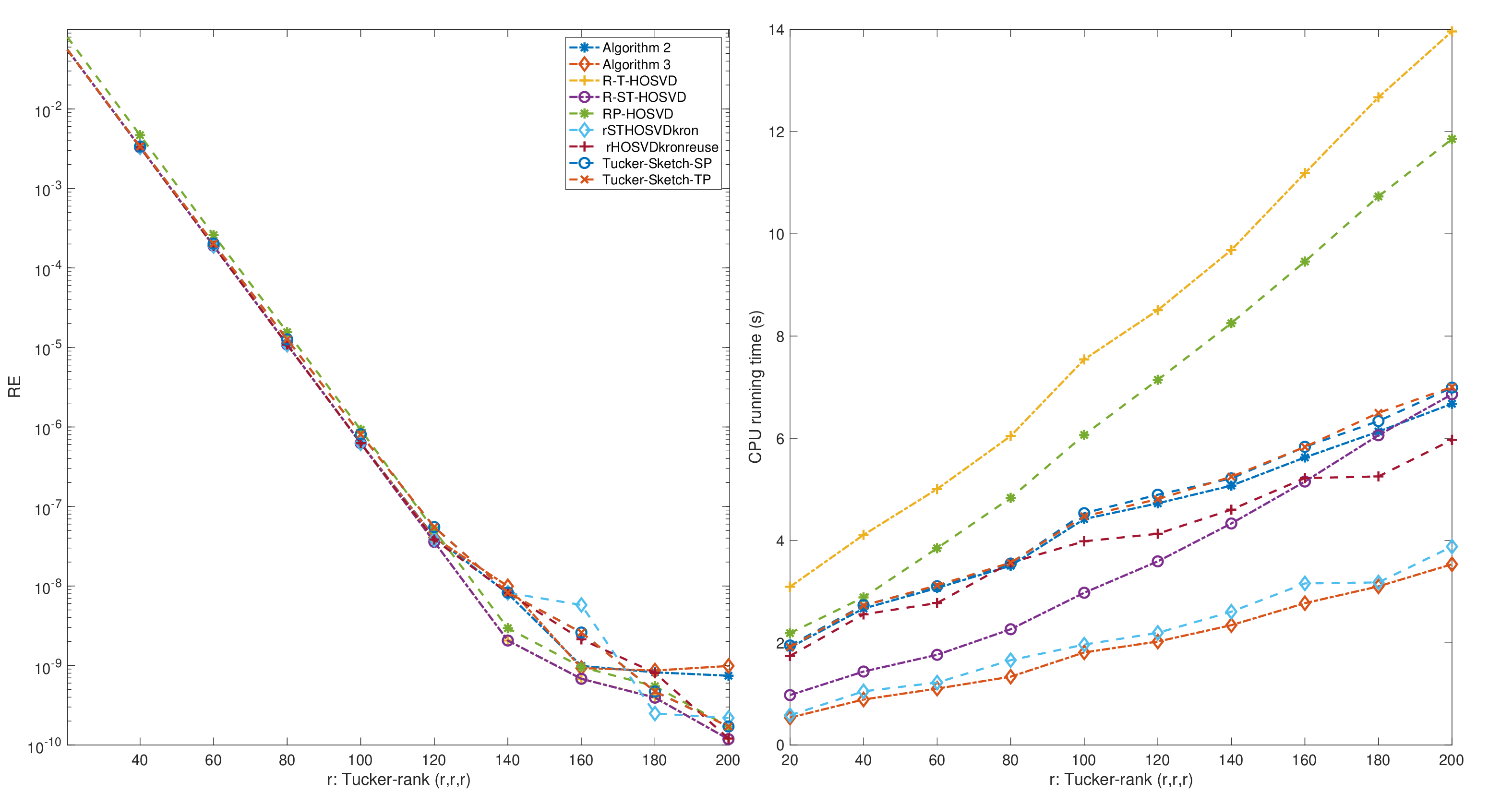}}\\
	\subfigure[The tensor $\mathcal{C}$ with S-shape decay]{\includegraphics[width=4.5in,height=1.6in]{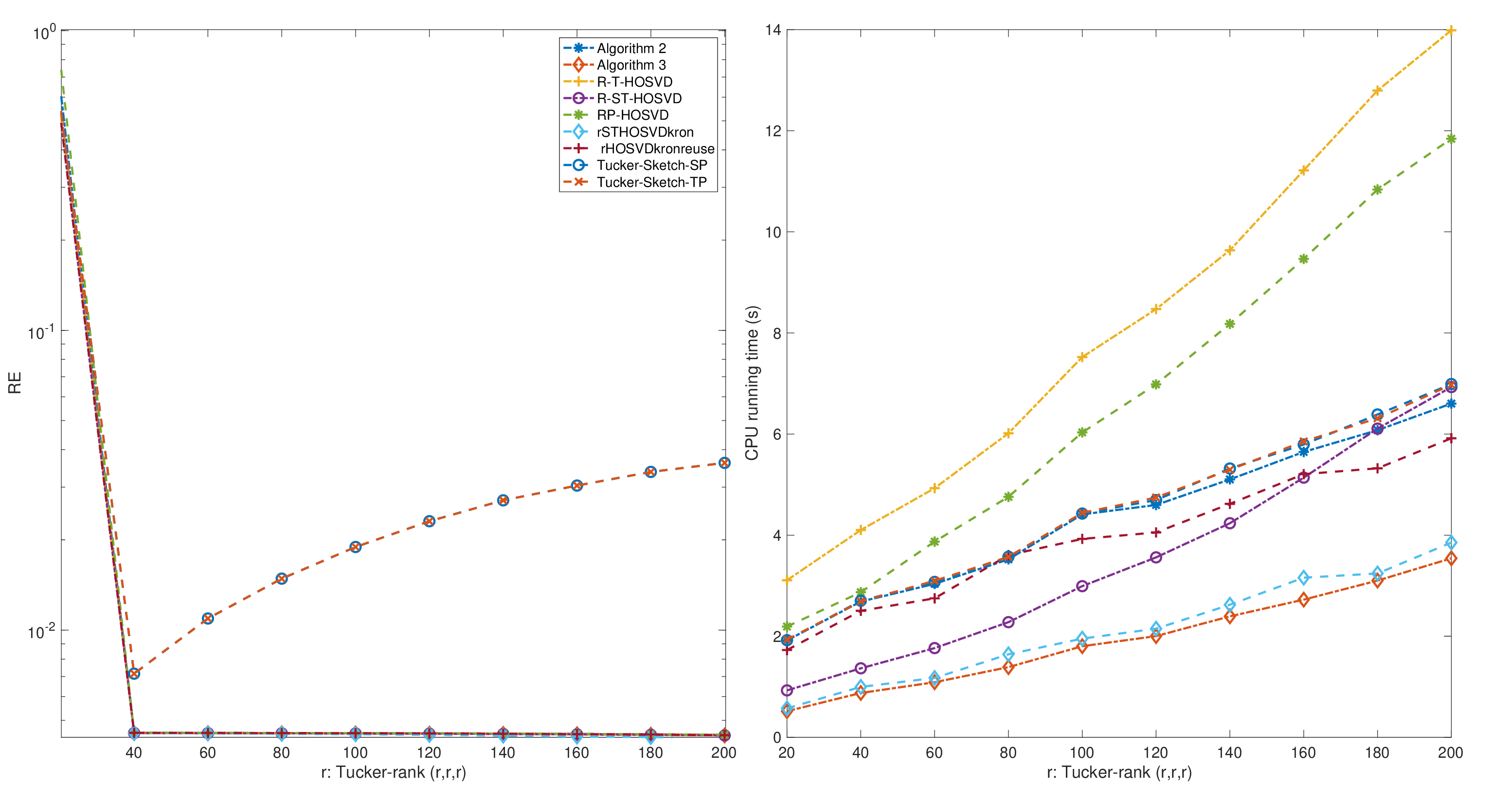}}\
    \caption{For fixed $s=10$ and $q=1$, numerical simulation results of Algorithms \ref{dash-rthosvd:alg2} and \ref{dash-rthosvd:alg2:v2}, and the existing randomized variants of T-HOSVD and ST-HOSVD on the test tensors $\mathcal{B}$ and $\mathcal{C}$.}\label{dash-rthosvd:figure5}
\end{figure}

\begin{figure}
    \setlength{\tabcolsep}{4pt}
    \renewcommand\arraystretch{1}
    \centering
    \subfigure[The extended Yale B database]{\includegraphics[width=4.5in,height=1.6in]{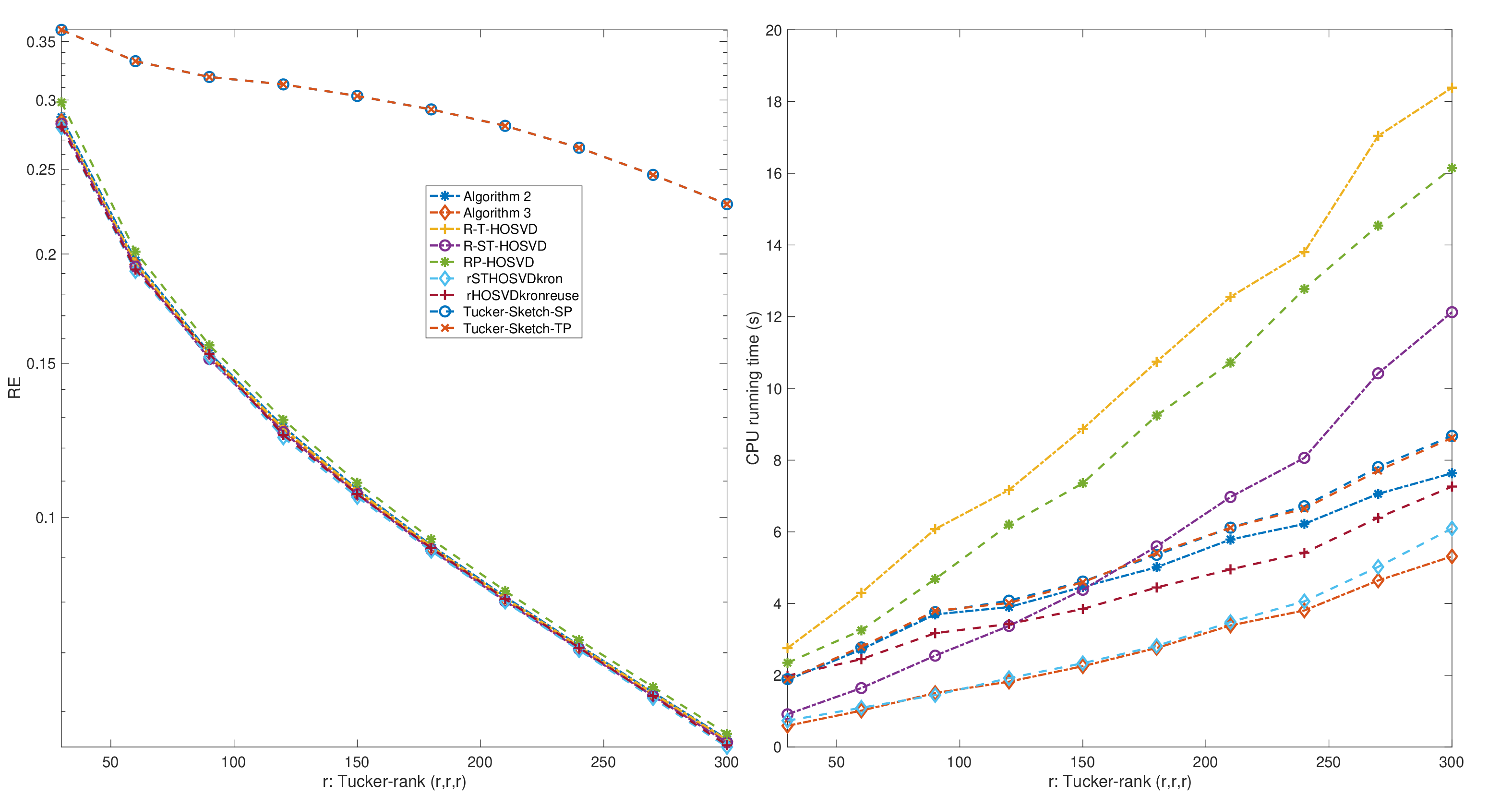}}\\
	\subfigure[COIL-100]{\includegraphics[width=4.5in,height=1.6in]{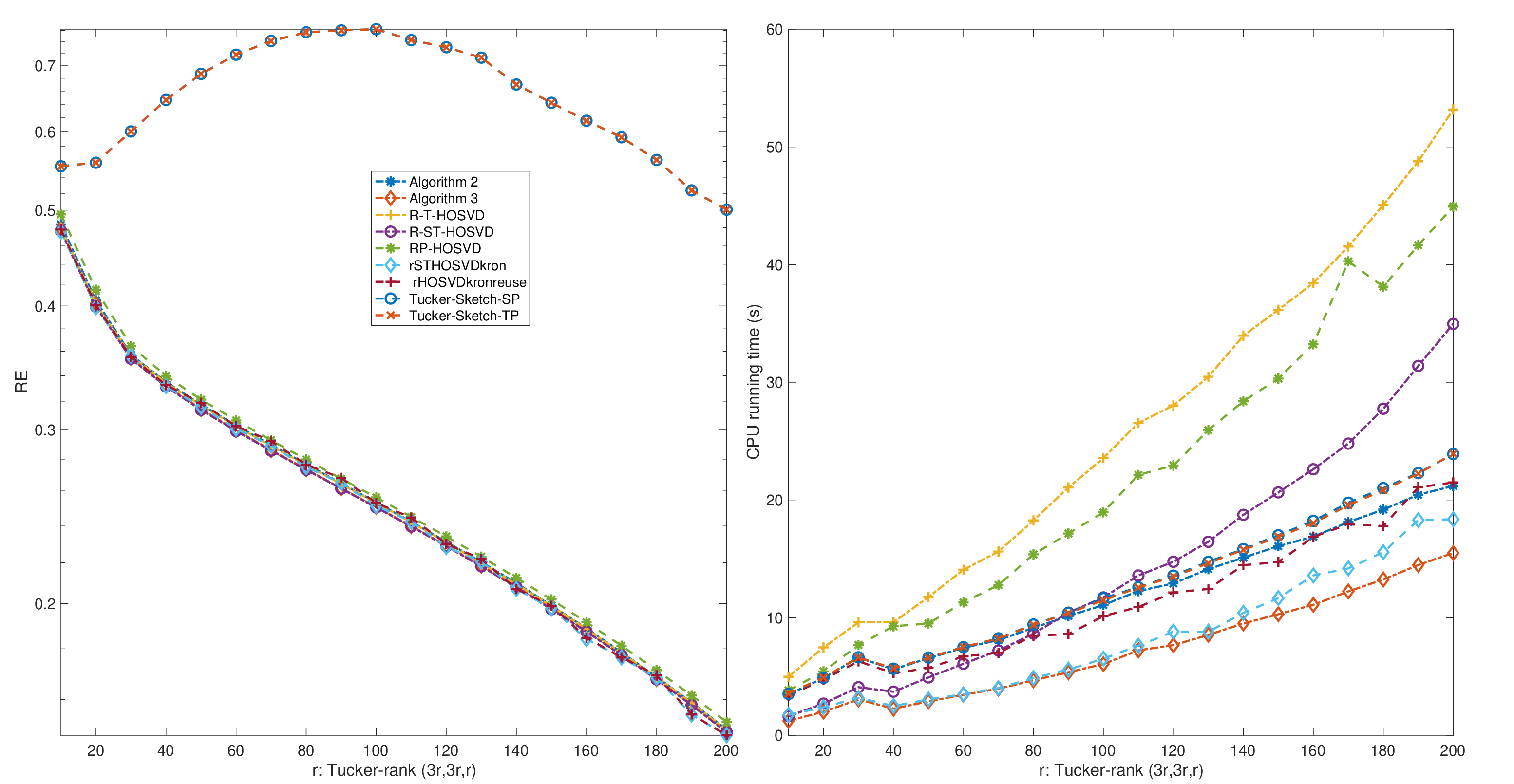}}\\
	\subfigure[The Washington DC Mall database]{\includegraphics[width=4.5in,height=1.6in]{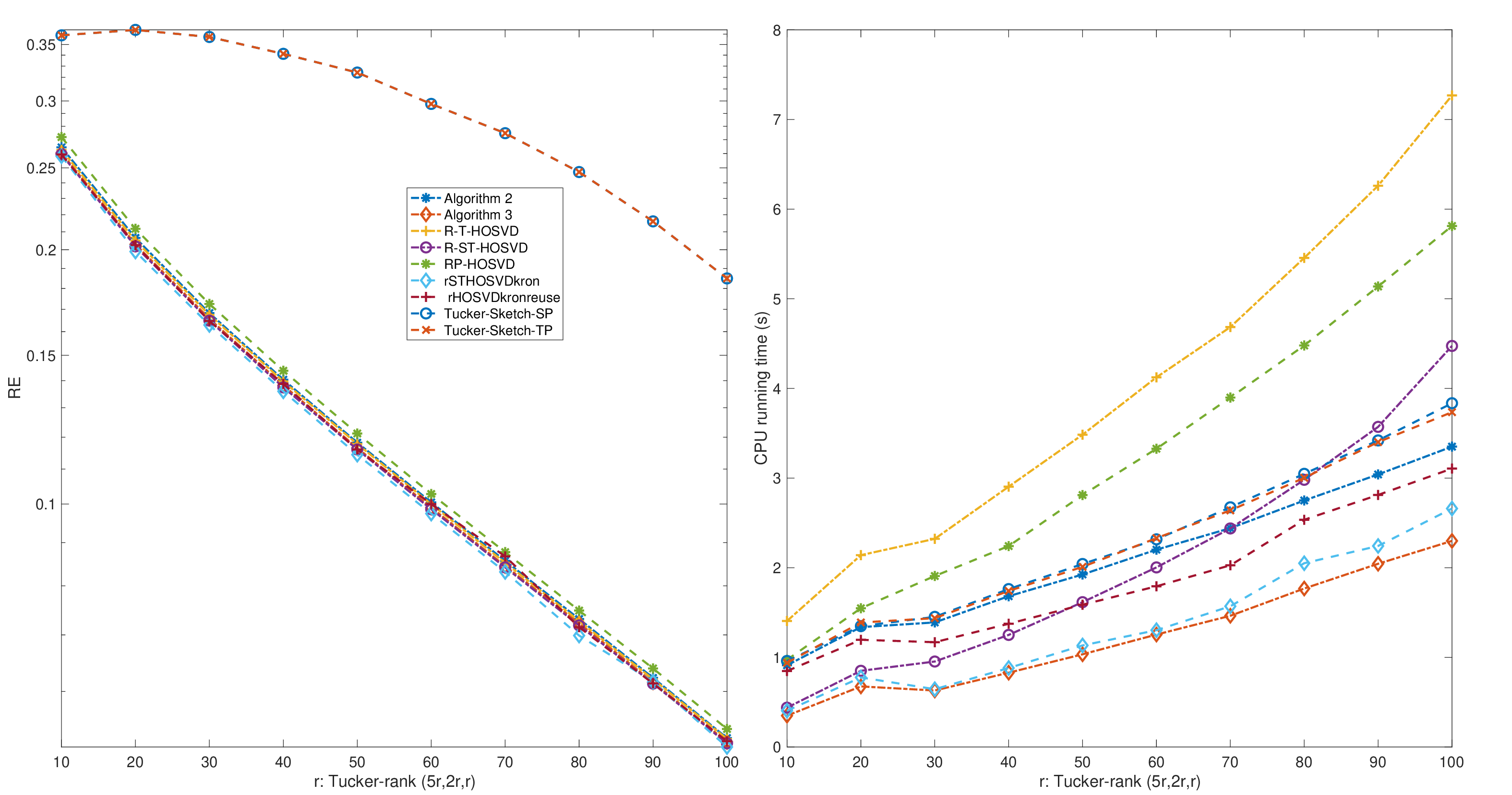}}\\
    \caption{For fixed $s=10$ and $q=1$, numerical simulation results of Algorithms \ref{dash-rthosvd:alg2} and \ref{dash-rthosvd:alg2:v2}, and the existing randomized variants of T-HOSVD and ST-HOSVD on the test tensors from three real databases.}\label{dash-rthosvd:figure6}
\end{figure}

We first compared Algorithms \ref{dash-rthosvd:alg2} and \ref{dash-rthosvd:alg2:v2} with existing randomized variants of T-HOSVD and ST-HOSVD via the test tensors $\mathcal{B}$, $\mathcal{C}$ and three tensors from three real databases. Following the descriptions above, we set $q=1$ in Algorithms \ref{dash-rthosvd:alg2} and \ref{dash-rthosvd:alg2:v2}. Hence, the power parameter in the randomized variants of T-HOSVD and ST-HOSVD was also set to 1. We set $r=10,20,\dots,100$ for $\mathcal{B}$, and $r=20,40,\dots,200$ for $\mathcal{C}$. The desired Tucker-rank for $\mathcal{A}_{\rm YaleB}$ is denoted as $(r,r,r)$ with $r=30,60,\dots,300$, the desired Tucker-rank for $\mathcal{A}_{\rm COIL}$ is denoted as $(3r,3r,r)$ with $r=10,20,\dots,200$, and the desired Tucker-rank for $\mathcal{A}_{\rm DCmall}$ is denoted as $(5r,2r,r)$ with $r=10,20,\dots,100$. The corresponding results are shown in Figures \ref{dash-rthosvd:figure5} and \ref{dash-rthosvd:figure6}.

{\color{red}According to these two figures}, we conclude that (a) in terms of RE, Algorithms \ref{dash-rthosvd:alg2} and \ref{dash-rthosvd:alg2:v2} are comparable to R-T-HOSVD, R-ST-HOSVD, RP-HOSVD, rSTHOSVDkron and rHOSVDkronreuse, and better than Tucker-Sketch-SP and Tucker-Sketch-TP with Slow decay and S-shape decay; and (b) in terms of CPU run time, Algorithm \ref{dash-rthosvd:alg2:v2} is the fastest one, and Algorithm \ref{dash-rthosvd:alg2} is faster than R-T-HOSVD and RP-HOSVD, similar to Tucker-Sketch-SP and Tucker-Sketch-TP, and slower than R-ST-HOSVD, rSTHOSVDkron and rHOSVDkronreuse.

For comparison of real databases, when we compared the proposed algorithms with Tucker-ALS and its randomized variants, we set the desired Tucker-ranks $(r_1,r_2,r_3)$ for Yale, COIL-100 and Washington DC as $(50,50,50)$, $(30,30,15)$ and $(40,30,20)$, respectively, and when we compared the proposed algorithms with T-HOSVD, ST-HOSVD and their randomized variants, we set the desired Tucker-ranks $(r_1,r_2,r_3)$ for Yale, COIL-100 and Washington DC as $(200,200,200)$, $(450,450,120)$ and $(500,100,50)$, respectively. The corresponding results for RE and CPU run times are shown in Tables \ref{dash-rthosvd:tab1} and \ref{dash-rthosvd:tab2}, and demonstrate that: (a) Algorithm \ref{dash-rthosvd:alg2:v2} is comparable to and faster than Algorithm \ref{dash-rthosvd:alg2}, Tucker-ALS and its randomized variants; (b) in terms of RE, Algorithms \ref{dash-rthosvd:alg2} and \ref{dash-rthosvd:alg2:v2} are comparable to T-HOSVD, ST-HOSVD and the randomized algorithms (i.e, R-T-HOSVD, R-ST-HOSVD, R-PET, RP-HOSVD, rSTHOSVD-kron and rHOSVDkronreuse) with $q=1$, and better than these randomized algorithms with $q=0$; and (c) {\color{red}in terms of CPU run time, Algorithms \ref{dash-rthosvd:alg2:v2} is slightly slower than ST-HOSVD for Yale and slightly faster than ST-HOSVD for COIL-100 and Washington-DC}, Algorithms \ref{dash-rthosvd:alg2} and \ref{dash-rthosvd:alg2:v2} are slower than RP-HOSVD, rSTHOSVDkron and rSTHOSVDkronreuse with $q=0$ and faster than T-HOSVD and the randomized algorithms (i.e, R-T-HOSVD, R-ST-HOSVD, R-PET, RP-HOSVD, rSTHOSVD-kron and rHOSVDkronreuse) with $q=1$.

\begin{sidewaystable}
    \scriptsize
    \centering
    \begin{tabular}{|c|cc|cc|cc|}
        \hline
        \multirow{2}{*}{Algorithms}  & \multicolumn{2}{c|}{Yale} & \multicolumn{2}{c|}{COIL-100} & \multicolumn{2}{c|}{Washington DC}  \\
        \cline{2-7}
        &  RE      &  CPU run time (s)   &   RE      &  CPU run time (s)  &   RE      &  CPU run time (s)     \\
        \hline
        Algorithm \ref{dash-rthosvd:alg2}
        & 2.17e-1 & 2.31
        & 4.53e-1 & 3.53
        & 2.49e-1 & 1.02    \\
        \hline
        Algorithm \ref{dash-rthosvd:alg2:v2}
        & 2.15e-1 & 0.86
        & 4.46e-1 & 1.25
        & 2.40e-1 & 0.33    \\
        \hline
        Tucker-ALS
        & 3.45e-1 & 1.06
        & 4.34e-1 & 6.43
        & 2.30e-1 & 1.54      \\
        \hline
        Randomized-Tucker-ALS
        & 3.61e-1 & 6.12
        & 4.48e-1 & 26.12
        & 2.36e-1 & 3.51      \\
        \hline
        Tucker-TS
        & 4.14e-1 & 501.03
        & 5.10e-1 & 2.04e+3
        & 2.72e-1 & 3.56e+4    \\
        \hline
        Tucker-TTMTS
        & 6.65e-1 & 9.34
        & 7.23e-1 & 32.54
        & 5.99e-1 & 67.72 \\
        \hline
        Sketch-Tucker-ALS
        & 3.79e-1 & 60.15
        & 4.58e-1 & 137.01
        & 2.72e-1 & 65.21    \\
    \hline
    \end{tabular}
    \caption{Numerical simulation results of Algorithms \ref{dash-rthosvd:alg2} and \ref{dash-rthosvd:alg2:v2} with Tucker-ALS, randomized Tucker-ALS, Tucker-TS, Tucker-TTMTS, and Sketch-Tucker-ALS on the tensors from three real databases.}
    \label{dash-rthosvd:tab1}
\end{sidewaystable}

\begin{sidewaystable}
    \scriptsize
    \centering
    \begin{tabular}{|c|cc|cc|cc|}
        \hline
        \multirow{2}{*}{Algorithms}  & \multicolumn{2}{c|}{Yale} & \multicolumn{2}{c|}{COIL-100} & \multicolumn{2}{c|}{ Washington DC}  \\
        \cline{2-7}
        &  RE      &  CPU run time (s)   &   RE      &  CPU run time (s)  &   RE      &  CPU run time (s)     \\
        \hline
        Algorithm \ref{dash-rthosvd:alg2}
        & 8.48e-2 & 5.80
        & 2.22e-1 & 15.37
        & 1.12e-1 & 2.63    \\
        \hline
        Algorithm \ref{dash-rthosvd:alg2:v2}
        & 8.40e-2 & 3.34
        & 2.21e-1 & 10.19
        & 1.12e-1 & 1.86    \\
        \hline
        T-HOSVD
        & 7.95e-2 & 25.92
        & 2.08e-1 & 78.58
        & 1.06e-1 & 9.64     \\
        \hline
        ST-HOSVD
        & 7.94e-2 & 2.99
        & 2.08e-1 & 10.24
        & 1.06e-1 & 1.85 \\
        \hline
        R-T-HOSVD ($q=0$)
        & 1.34e-1 & 10.62
        & 2.87e-1 & 52.19
        & 1.61e-1 & 4.46       \\
        \hline
        R-T-HOSVD ($q=1$)
        & 8.46e-2 & 12.68
        & 2.22e-1 & 38.85
        & 1.12e-1 & 5.39     \\
        \hline
        R-ST-HOSVD ($q=0$)
        & 1.29e-1 & 5.20
        & 2.81e-1 & 16.48
        & 1.58e-1 & 2.79       \\
        \hline
        R-ST-HOSVD ($q=1$)
        & 8.38e-2 & 9.01
        & 2.21e-1 & 28.04
        & 1.11e-1 & 4.72 \\
        \hline
        RP-HOSVD ($q=0$)
        & 1.40e-1 & 3.46
        & 2.94e-1 & 10.65
        & 1.69e-1 & 1.69       \\
        \hline
        RP-HOSVD ($q=1$)
        & 8.61e-2 & 10.63
        & 2.25e-1 & 42.27
        & 1.15e-1 & 4.02 \\
        \hline
        R-PET ($q=0$)
        & 1.89e-1 & 3.97
        & 4.02e-1 & 20.22
        & 2.25e-1 & 2.00       \\
        \hline
        R-PET ($q=1$)
        & 1.17e-1 & 11.42
        & 3.08e-1 & 41.19
        & 1.53e-1 & 4.87   \\
        \hline
        rSTHOSVDkron ($q=0$)
        & 1.32e-1 & 2.12
        & 3.10e-1 & 7.07
        & 1.82e-1 & 1.22       \\
        \hline
        rSTHOSVDkron ($q=1$)
        & 8.27e-2 & 6.62
        & 2.20e-1 & 21.09
        & 1.11e-1 & 3.03 \\
        \hline
        rSTHOSVDkronreuse ($q=0$)
        & 1.36e-1 & 2.98
        & 3.07e-1 & 10.12
        & 2.05e-1 & 1.50      \\
        \hline
        rSTHOSVDkronreuse ($q=1$)
        & 8.33e-2 & 11.38
        & 2.21e-1 & 49.26
        & 1.12e-1 & 4.32 \\
        \hline
    \end{tabular}
    \caption{Numerical simulation results of Algorithms \ref{dash-rthosvd:alg2} and \ref{dash-rthosvd:alg2:v2} with T-HOSVD, ST-HOSVD, and several randomized variants to the tensors from three real databases.}
    \label{dash-rthosvd:tab2}
\end{sidewaystable}

\subsection{Numerical results for Algorithm \ref{dash-rthosvd:alg3:app}}
\label{dash-rthosvd:sec5:subsec6}

For clarity, we set $s_1=s_2=s_3=10$, {\color{red}$q_{\max}=500$} and $r_1=r_2=r_3:=r$. For each $(r,r,r)$, when we apply Algorithm \ref{dash-rthosvd:alg3:app} with different $\textit{tol}$ to the tensors $\mathcal{B}$ and $\mathcal{C}$ in Section \ref{dash-rthosvd:sec5:subsec1}, the values of $(q_1,q_2,q_3)$ and RE are shown in Table \ref{dash-rthosvd:tab1:appsec3}, which implies that the case of $\textit{tol}=0.5$ is suitable when applying Algorithm \ref{dash-rthosvd:alg3:app} to find a quasi-optimal solution for Problem \ref{dash-rthosvd:prob1}.
\begin{sidewaystable}
	\scriptsize
	\centering
	\begin{tabular}{|c|c|cc|cc|cc|cc|}
            \hline
            \multirow{2}{*}{Tensors}  &  \multirow{2}{*}{r} & \multicolumn{2}{c|}{0.5} &  \multicolumn{2}{c|}{0.1} & \multicolumn{2}{c|}{0.05} & \multicolumn{2}{c|}{0.01} \\
            \cline{3-10}
            & & $(q_1,q_2,q_3)$ & RE & $(q_1,q_2,q_3)$ & RE & $(q_1,q_2,q_3)$ & RE & ${\color{red}(q_1,q_2,q_3)}$ & RE \\
            \hline
            \multirow{10}{*}{$\mathcal{B}$} & 10
            & (2,2,2) & 1.81e-2 & (3,2,2) & 1.81e-2 & (3,3,2) & 1.81e-2
            & (3,3,3) & 1.81e-2 \\
            \cline{2-10}
            \multirow{10}{*}{} & 20
            & (2,2,2) & 6.26e-3 & (3,3,3) & 6.26e-3 & (3,3,3) & 6.26e-3
            & (3,3,3) & 6.26e-3 \\
            \cline{2-10}
            \multirow{10}{*}{} & 30
            & (2,2,2) & 3.11e-3 & (3,3,3) & 3.11e-3 & (3,3,3) & 3.11e-3
            & (4,3,3) & 3.11e-3 \\
            \cline{2-10}
            \multirow{10}{*}{} & 40
            & (2,2,2) & 1.62e-3 & (2,2,2) & 1.62e-3 & (2,2,2) & 1.62e-3
            & (2,2,2) & 1.62e-3 \\
            \cline{2-10}
            \multirow{10}{*}{} & 50
            & (3,3,3) & 1.66e-6 & (3,3,3) & 1.66e-6 & (3,3,4) & 1.66e-6
            & (3,4,8) & 1.66e-6 \\
            \cline{2-10}
            \multirow{10}{*}{} & 60
            & (3,2,2) & 1.27e-6 & (3,3,3) & 1.27e-6 & (4,3,3) & 1.27e-6
            & (5,3,4) & 1.27e-6 \\
            \cline{2-10}
            \multirow{10}{*}{} & 70
            & (3,2,2) & 1.01e-6 & (3,3,3) & 1.01e-6 & (4,3,4) & 1.01e-6
            & (5,5,7) & 1.01e-6 \\
            \cline{2-10}
            \multirow{10}{*}{} & 80
            & (3,3,3) & 8.18e-7 & (3,3,6) & 8.18e-7 & (4,4,16) & 8.17e-7
            & (54,109,87) & 6.00e-4 \\
            \cline{2-10}
            \multirow{10}{*}{} & 90
            & (3,3,3) & 6.96e-7 & (3,4,10) & 6.95e-7 & (5,7,20) & 6.98e-7
            & (76,134,87) & 6.65e-4 \\
            \cline{2-10}
            \multirow{10}{*}{} & 100
            & (3,3,13) & 5.68e-7 & (5,16,28) & 6.21e-7 & (28,28,47) & 8.79e-7
            & (92,148,110) & 6.75e-4 \\
            \hline

            \multirow{10}{*}{$\mathcal{C}$ with Slow decay} & 20
            & (2,2,5) & 5.98e-3 & (3,2,7) & 5.98e-3 & (3,2,6) & 5.98e-3
            & (3,3,11) & 5.98e-3 \\
            \cline{2-10}
            \multirow{10}{*}{} & 40
            & (3,2,3) & 2.15e-3 & (3,2,6) & 2.15e-3 & (3,2,8) & 2.15e-3
            & (4,3,11) & 2.15e-3 \\
            \cline{2-10}
            \multirow{10}{*}{} & 60
            & (3,2,5) & 1.18e-3 & (3,2,5) & 1.18e-3 & (4,3,10) & 1.18e-3
            & (5,3,12) & 1.18e-3 \\
            \cline{2-10}
            \multirow{10}{*}{} & 80
            & (3,2,4) & 7.70e-4 & (3,2,6) & 7.69e-4 & (4,3,7) & 7.68e-4
            & (5,3,11) & 7.67e-4 \\
            \cline{2-10}
            \multirow{10}{*}{} & 100
            & (3,2,4) & 5.52e-4 & (3,2,6) & 5.11e-4 & (4,3,7) & 5.50e-4
            & (5,3,13) & 5.50e-4 \\
            \cline{2-10}
            \multirow{10}{*}{} & 120
            & (3,2,3) & 4.20e-4 & (3,3,7) & 4.20e-4 & (4,3,8) & 4.18e-4
            & (6,3,12) & 4.18e-4 \\
            \cline{2-10}
            \multirow{10}{*}{} & 140
            & (3,2,3) & 3.34e-4 & (4,2,5) & 3.32e-4 & (4,3,8) & 3.32e-4
            & (6,3,14) & 3.31e-4 \\
            \cline{2-10}
            \multirow{10}{*}{} & 160
            & (3,2,4) & 2.73e-4 & (4,3,6) & 2.71e-4 & (4,3,8) & 2.71e-4
            & (46,12,11) & 8.47e-4 \\
            \cline{2-10}
            \multirow{10}{*}{} & 180
            & (3,2,3) & 2.28e-4 & (4,2,7) & 2.27e-4 & (4,3,8) & 2.26e-4
            & (27,8,13) & 4.89e-4 \\
            \cline{2-10}
            \multirow{10}{*}{} & 200
            & (3,2,4) & 1.94e-4 & (4,3,8) & 1.93e-4 & (4,3,11) & 1.92e-4
            & (76,21,20) & 1.05e-3 \\
            \hline

            \multirow{10}{*}{$\mathcal{C}$ with Fast decay} & 20
            & (2,2,11) & 5.74e-2 & (3,8,15) & 5.74e-2 & (3,22,15) & 5.74e-2
            & (3,24,20) & 5.74e-2 \\
            \cline{2-10}
            \multirow{10}{*}{} & 40
            & (2,2,10) & 3.30e-3 & (3,10,13) & 3.30e-3 & (3,18,16) & 3.30e-3
            & (3,21,19) & 3.30e-3 \\
            \cline{2-10}
            \multirow{10}{*}{} & 60
            & (2,8,9) & 1.89e-4 & (3,16,13) & 1.89e-4 & (3,15,15) & 1.89e-4
            & (3,18,19) & 1.89e-4 \\
            \cline{2-10}
            \multirow{10}{*}{} & 80
            & (3,9,9) & 1.09e-5 & (3,14,13) & 1.09e-5 & (3,15,15) & 1.09e-5
            & (3,18,19) & 1.09e-5 \\
            \cline{2-10}
            \multirow{10}{*}{} & 100
            & (2,10,9) & 6.25e-7 & (3,13,13) & 6.25e-7 & (3,15,15) & 6.25e-7
            & (3,19,19) & 6.25e-7 \\
            \cline{2-10}
            \multirow{10}{*}{} & 120
            & (8,14,10) & 4.66e-8 & (55,18,17) & 1.57e-7 & (75,19,20) & 2.78e-7
            & (104,24,26) & 5.40e-7 \\
            \cline{2-10}
            \multirow{10}{*}{} & 140
            & (33,17,12) & 4.66e-8 & (63,18,16) & 1.37e-7 & (80,18,20) & 2.42e-7
            & (115,25,26) & 5.67e-7 \\
            \cline{2-10}
            \multirow{10}{*}{} & 160
            & (33,16,10) & 4.01e-8 & (65,20,18) & 1.42e-7 & (80,20,20) & 2.51e-7
            & (117,25,27) & 5.88e-7 \\
            \cline{2-10}
            \multirow{10}{*}{} & 180
            & (33,18,11) & 4.53e-8 & (67,18,17) & 1.47e-7 & (77,20,20) & 2.33e-7
            & (118,25,27) & 6.00e-7 \\
            \cline{2-10}
            \multirow{10}{*}{} & 200
            & (33,15,11) & 4.72e-8 & (65,19,18) & 1.48e-7 & (79,19,20) & 2.42e-7
            & (115,25,27) & 6.00e-7 \\
            \hline

            \multirow{10}{*}{$\mathcal{C}$ with S-shape decay} & 20
            & (2,2,98) & 5.30e-1 & (2,2,119) & 5.28e-1 & (2,2,120) & 5.25e-1
            & (2,2,126) & 5.28e-1 \\
            \cline{2-10}
            \multirow{10}{*}{} & 40
            & (4,3,3) & 4.52e-3 & (6,3,24) & 4.51e-3 & (6,3,20) & 4.51e-3
            & (9,3,24) & 4.51e-3 \\
            \cline{2-10}
            \multirow{10}{*}{} & 60
            & (2,2,2) & 4.51e-3 & (5,5,3) & 4.51e-3 & (6,6,3) & 4.50e-3
            & (8,26,4) & 4.50e-3 \\
            \cline{2-10}
            \multirow{10}{*}{} & 80
            & (3,3,3) & 4.51e-3 & (5,4,3) & 4.50e-3 & (5,5,4) & 4.50e-3
            & (8,25,3) & 4.50e-3 \\
            \cline{2-10}
            \multirow{10}{*}{} & 100
            & (3,3,3) & 4.50e-3 & (4,4,4) & 4.49e-3 & (5,5,4) & 4.48e-3
            & (7,27,4) & 4.49e-3 \\
            \cline{2-10}
            \multirow{10}{*}{} & 120
            & (3,3,3) & 4.48e-3 & (4,4,4) & 4.47e-3 & (5,5,4) & 4.46e-3
            & (7,29,3) & 4.47e-3 \\
            \cline{2-10}
            \multirow{10}{*}{} & 140
            & (3,3,3) & 4.47e-3 & (4,4,4) & 4.45e-3 & (5,5,4) & 4.44e-3
            & (7,24,4) & 4.44e-3 \\
            \cline{2-10}
            \multirow{10}{*}{} & 160
            & (3,3,3) & 4.45e-3 & (4,4,4) & 4.43e-3 & (5,5,4) & 4.41e-3
            & (7,18,4) & 4.40e-3 \\
            \cline{2-10}
            \multirow{10}{*}{} & 180
            & (3,3,3) & 4.42e-3 & (4,4,4) & 4.40e-3 & (5,4,5) & 4.40e-3
            & (7,10,4) & 4.35e-3 \\
            \cline{2-10}
            \multirow{10}{*}{} & 200
            & (3,3,3) & 4.39e-3 & (4,4,4) & 4.36e-3 & (5,5,5) & 4.35e-3
            & (6,7,5) & 4.31e-3 \\
        \hline
	\end{tabular}
 \caption{
       Fixed $s_1=s_2=s_3=10$, for different Tucker-ranks $(r,r,r)$, the values of the 3-tuple $(q_1,q_2,q_3)$ and RE obtained by applying Algorithm \ref{dash-rthosvd:alg3:app} with $\textit{tol}$ on two test tensors $\mathcal{B}$ and $\mathcal{C}$ in Section \ref{dash-rthosvd:sec5:subsec1}.}	
 \label{dash-rthosvd:tab1:appsec3}
\end{sidewaystable}

\begin{figure}
	\setlength{\tabcolsep}{4pt}
	\renewcommand\arraystretch{1}
	\centering
    \subfigure[The tensor $\mathcal{A}$]{\includegraphics[width=4.5in,height=1.6in]{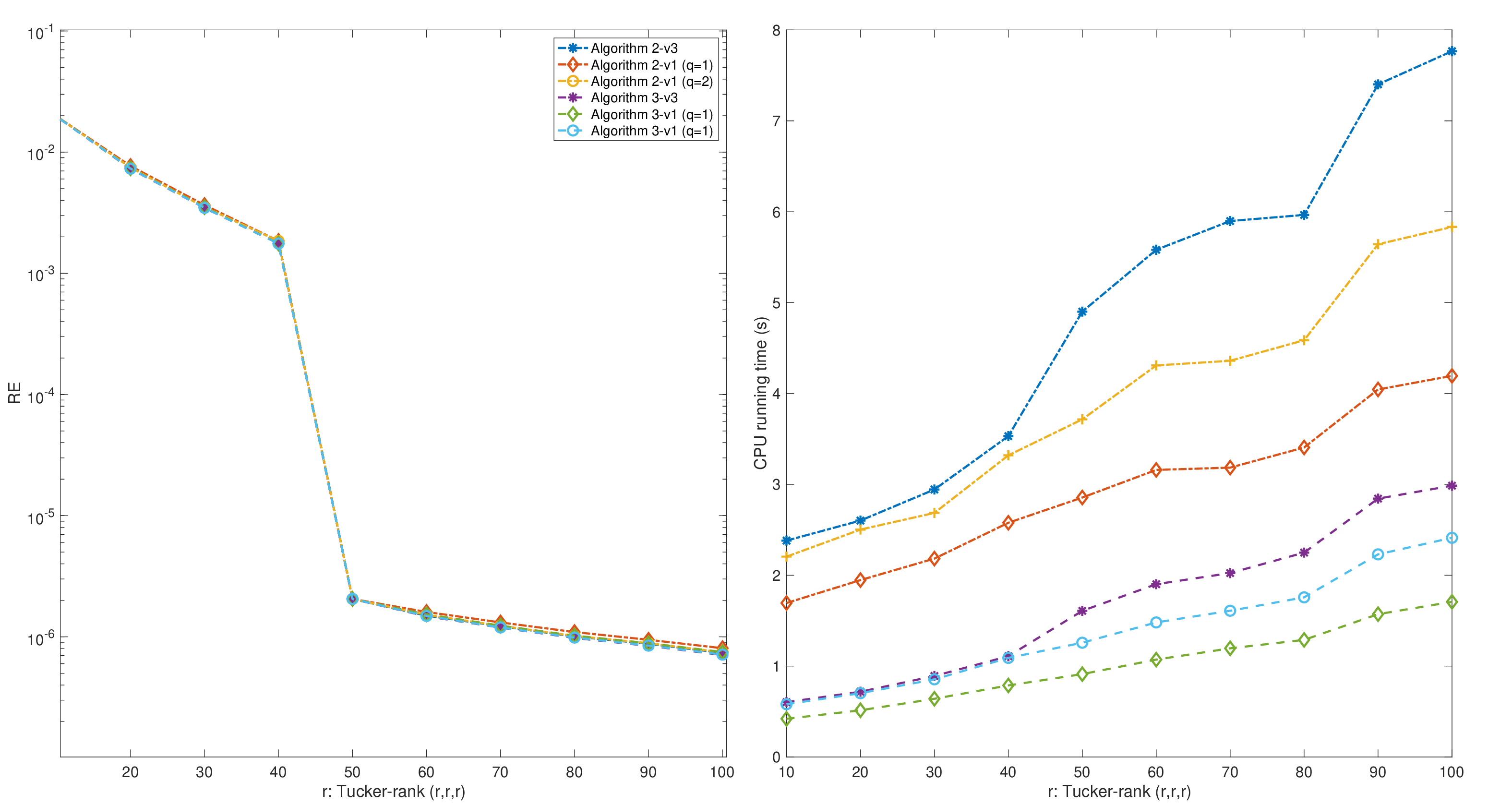}}\\
	\subfigure[The tensor $\mathcal{B}$ with Slow decay]{\includegraphics[width=4.5in,height=1.6in]{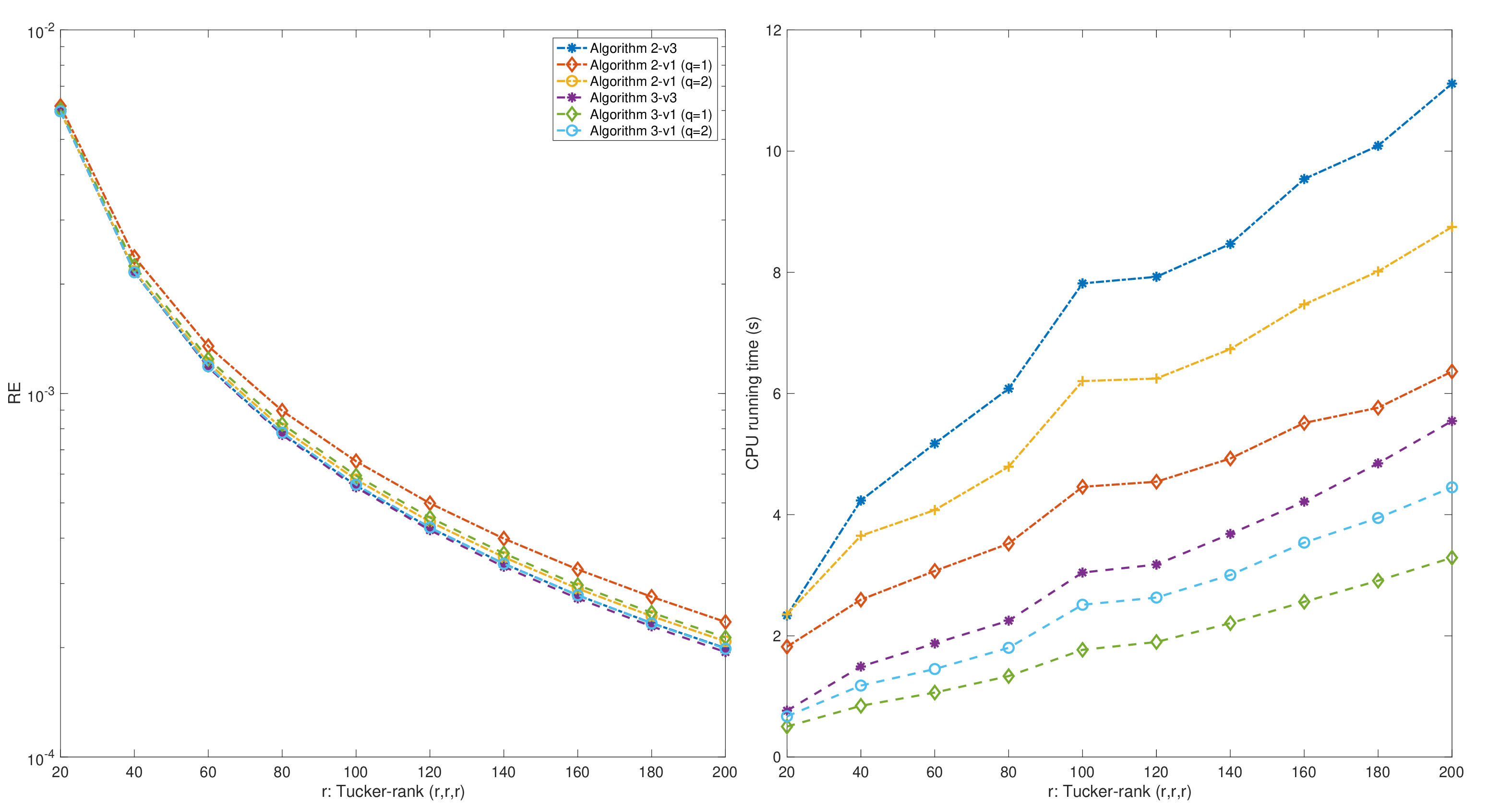}}\\
	\subfigure[The tensor $\mathcal{B}$ with Fast decay]{\includegraphics[width=4.5in,height=1.6in]{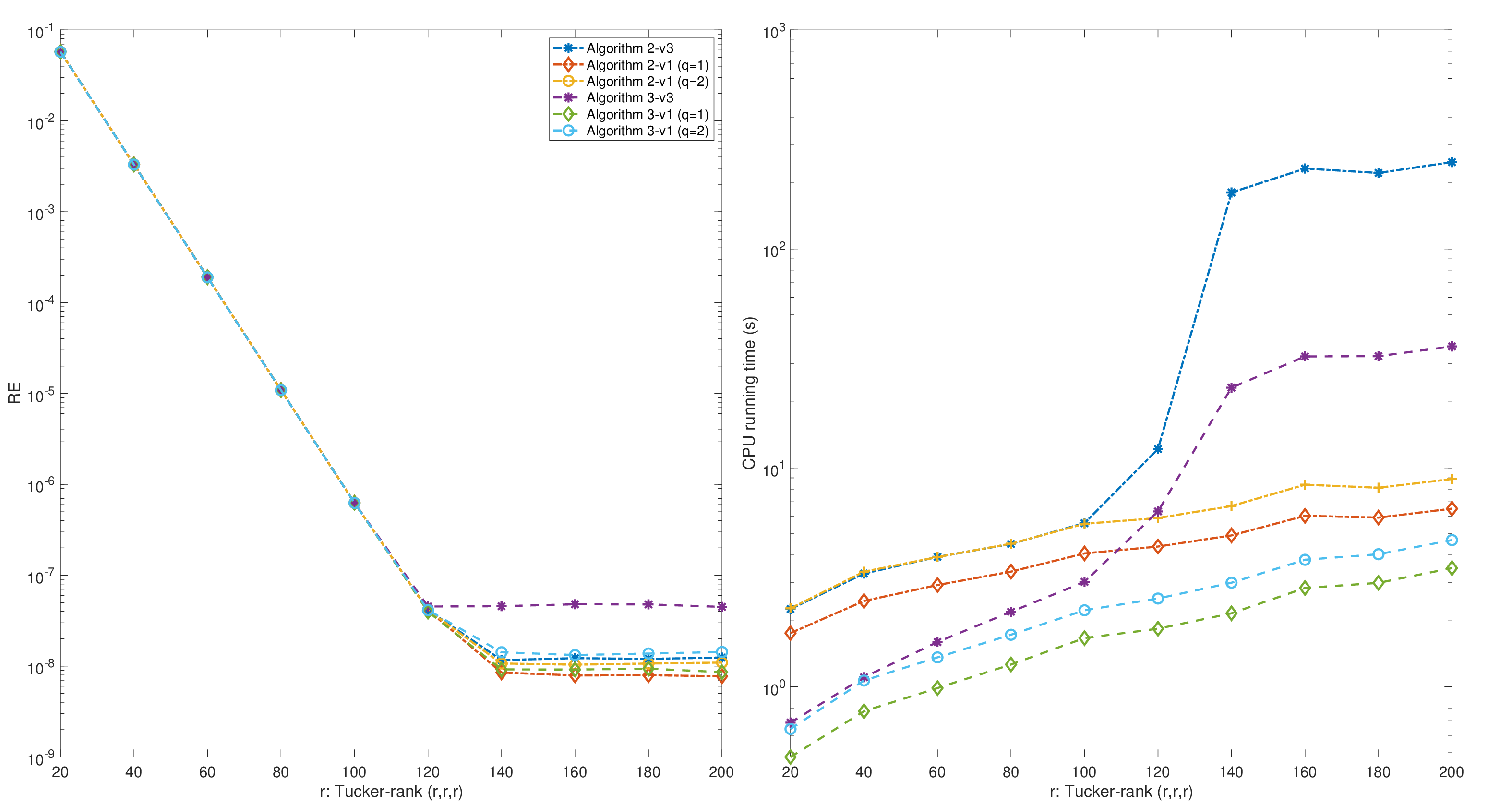}}\\
    \subfigure[The tensor $\mathcal{B}$ with S-shape decay]{\includegraphics[width=4.5in,height=1.6in]{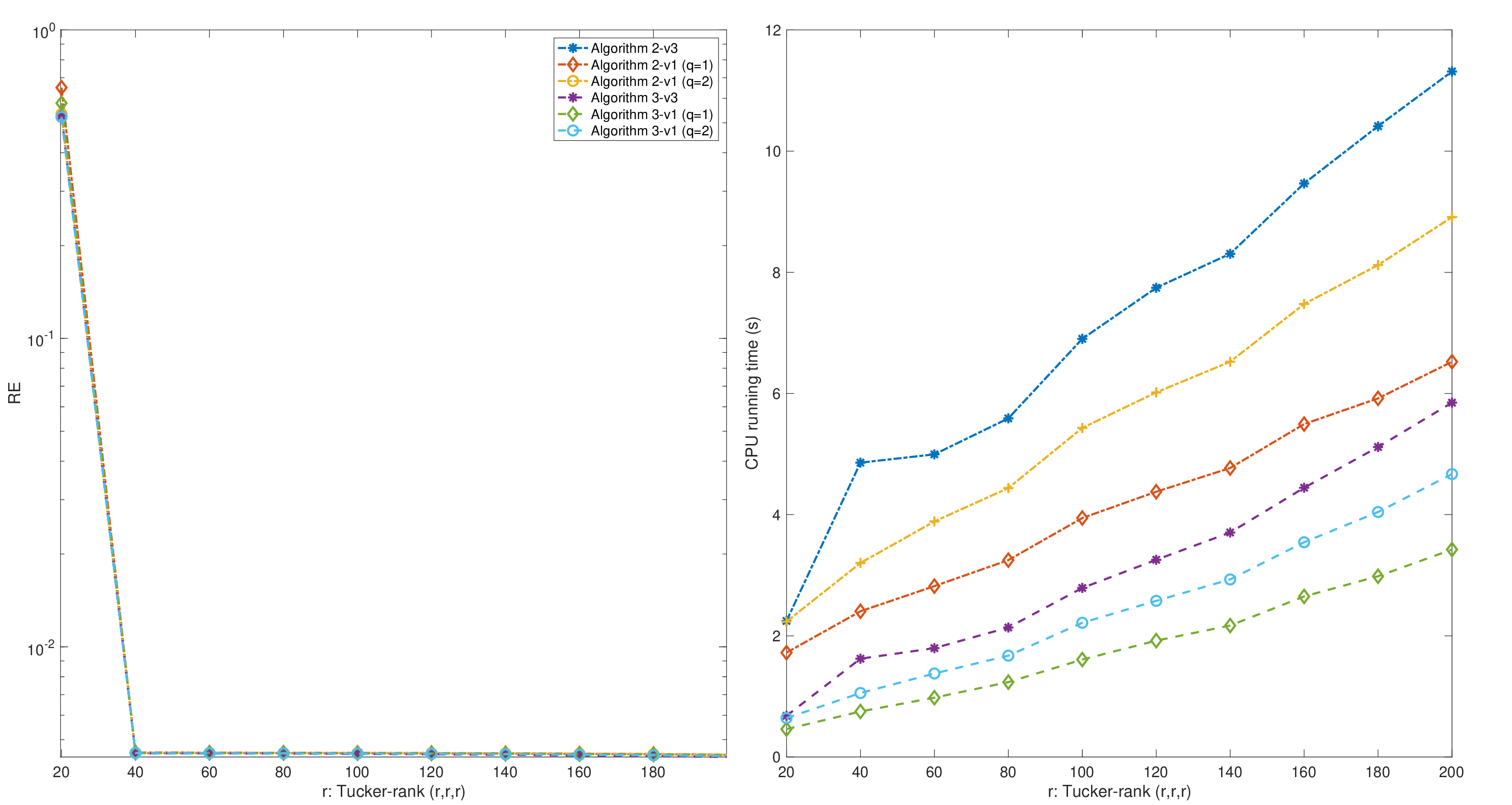}}\\
	\caption{For fixed $s_1=s_2=s_3=10$, the changes of values of RE and CPU run time by applying
    Algorithms \ref{dash-rthosvd:alg2}-v1/v3 and Algorithms \ref{dash-rthosvd:alg2:v2}-v1/v3
    using different Tucker-ranks $(r,r,r)$ on two test tensors in Section \ref{dash-rthosvd:sec5:subsec1}.}\label{dash-rthosvd:figure1:appsec3}
\end{figure}

Finally, with $q=1,2$ and $\textit{tol}=0.5$, we now compare the efficiency of Algorithm \ref{dash-rthosvd:alg2}-v1/v3 and Algorithms \ref{dash-rthosvd:alg2:v2}-v1/v3 via two test tensors $\mathcal{B}$ and $\mathcal{C}$ (see Section \ref{dash-rthosvd:sec5:subsec1}). For each 3-tuple $(r,r,r)$, the related results are shown in Figure \ref{dash-rthosvd:figure1:appsec3}. According to this figure, it can be seen that Algorithms \ref{dash-rthosvd:alg2} and \ref{dash-rthosvd:alg2:v2} with $q=1$ are suitable for the fixed-rank approximation in the Tucker decomposition.

\subsection{Comparison of different types of random matrices}
\label{dash-rthosvd:sec5:subsec5}
The proposed algorithms are based on standard Gaussian matrices. We now compare the efficiency of Algorithms \ref{dash-rthosvd:alg2} and \ref{dash-rthosvd:alg2:v2} under other types of random matrices, including any uniform random matrix (i.e., a matrix $\mathbf{\Omega}\in\mathbb{R}^{m\times k}$ is a uniform random matrix if its entries are drawn independently from the continuous uniform distribution in the interval $[-1,1]$), the Khatri-Rao product of standard Gaussian matrices, and the Khatri-Rao product of uniform random matrices.

\begin{figure}
	\setlength{\tabcolsep}{4pt}
	\renewcommand\arraystretch{1}
	\centering
	\subfigure[The tensor $\mathcal{B}$]{\includegraphics[width=4.5in,height=1.6in]{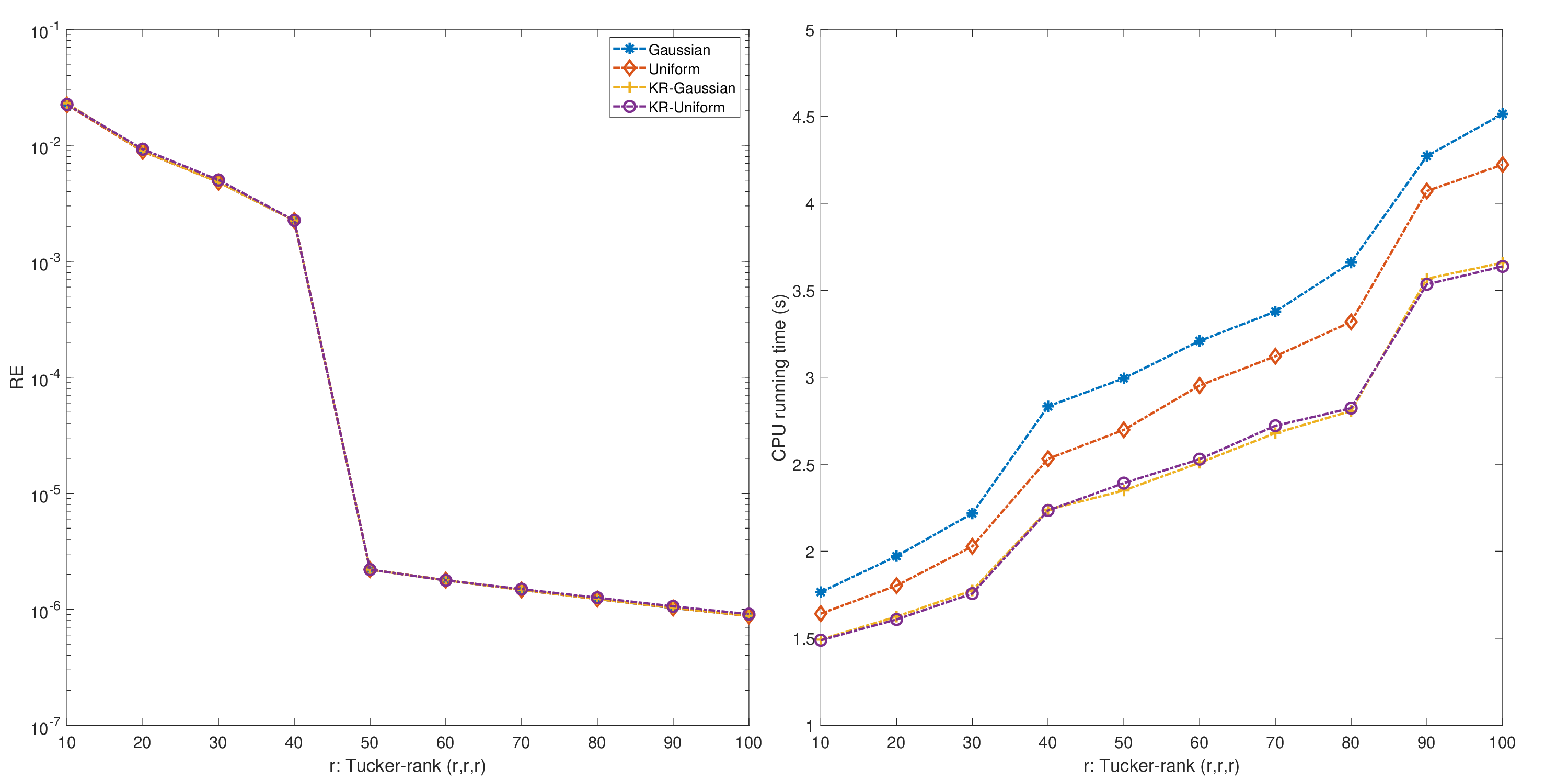}}\\
	\subfigure[The tensor $\mathcal{C}$ with Slow decay]{\includegraphics[width=4.5in,height=1.6in]{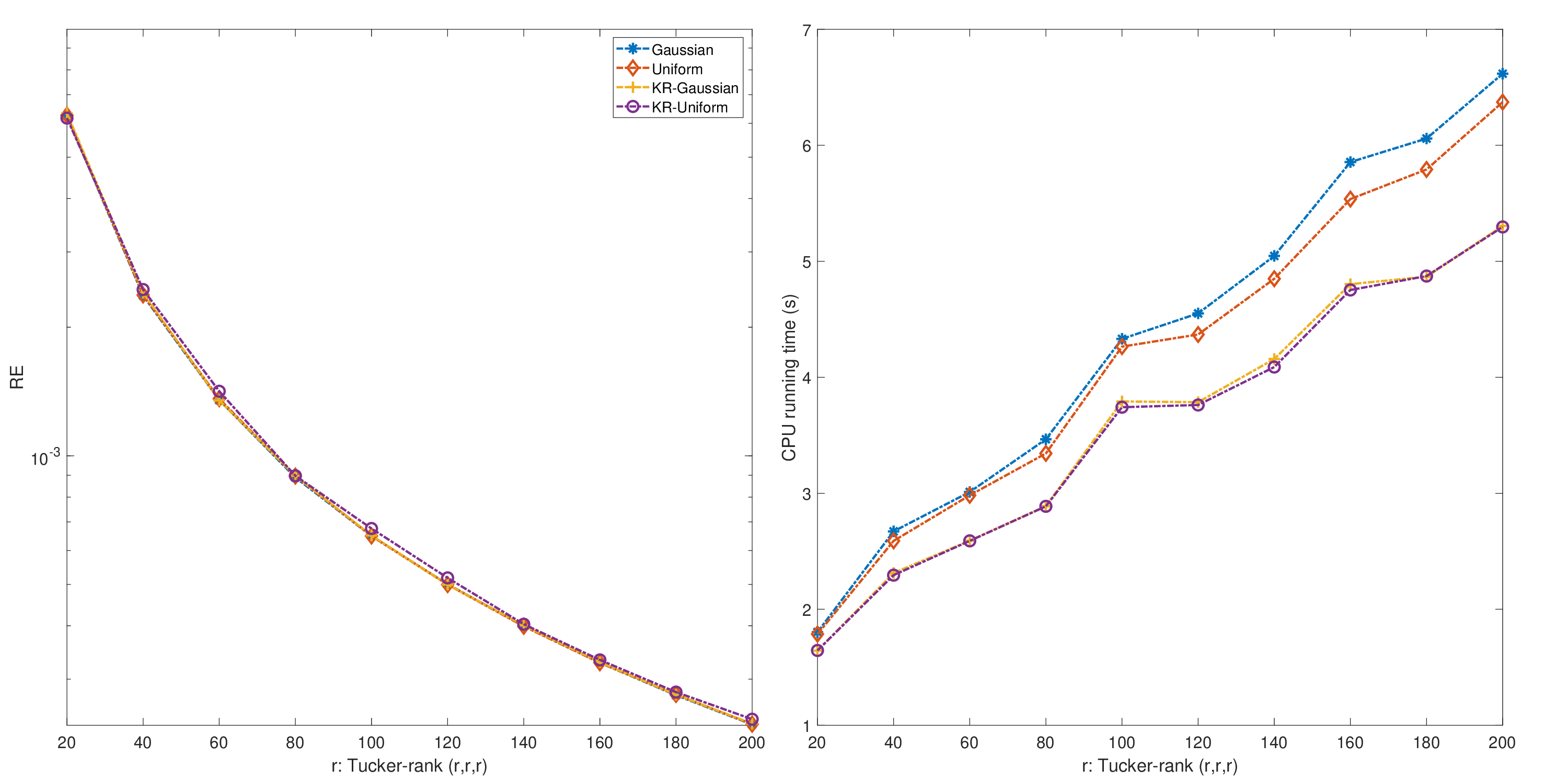}}\\
	\subfigure[The tensor $\mathcal{C}$ with Fast decay]{\includegraphics[width=4.5in,height=1.6in]{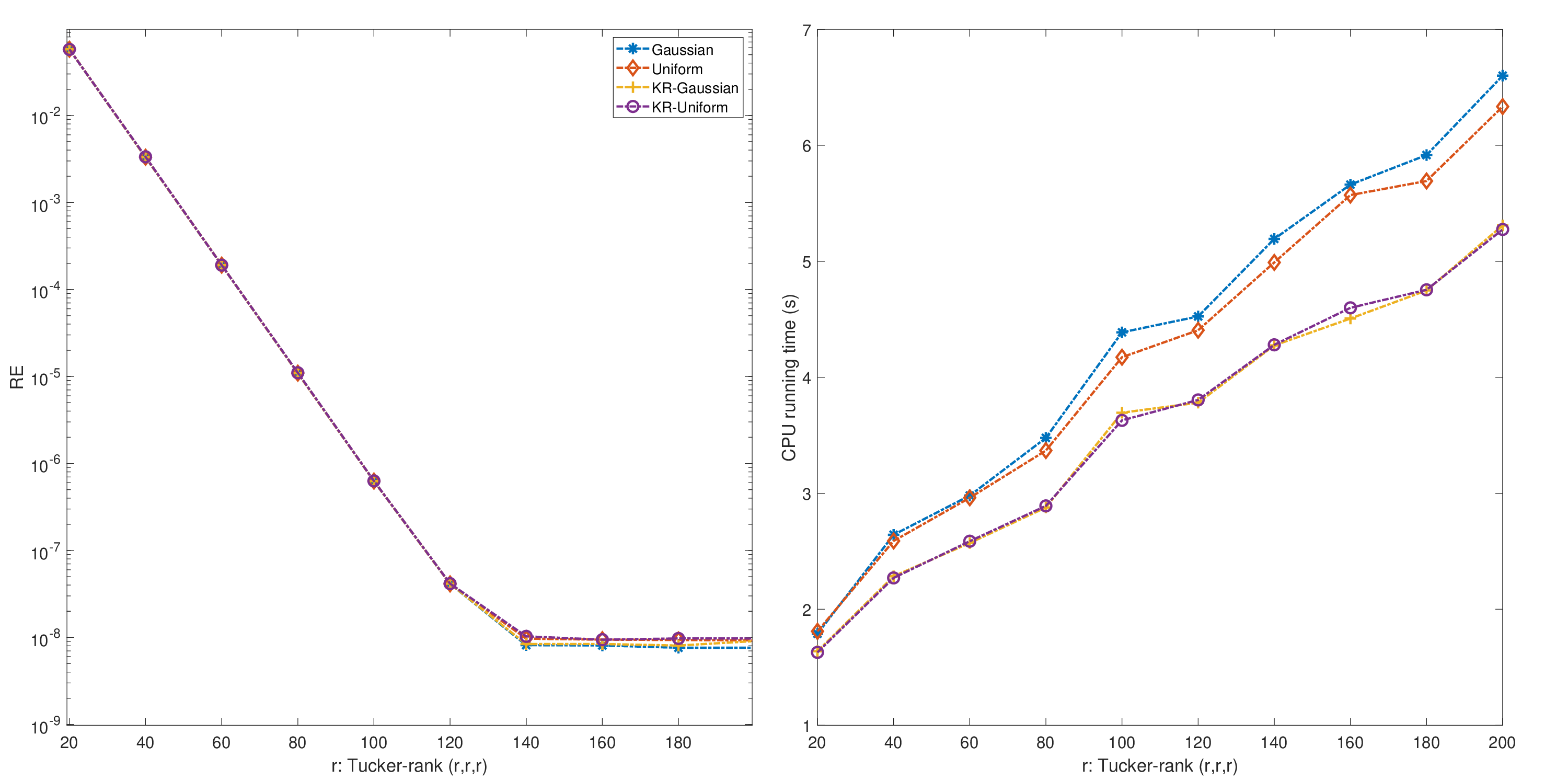}}\\
    \subfigure[The tensor $\mathcal{C}$ with S-shape decay]{\includegraphics[width=4.5in,height=1.6in]{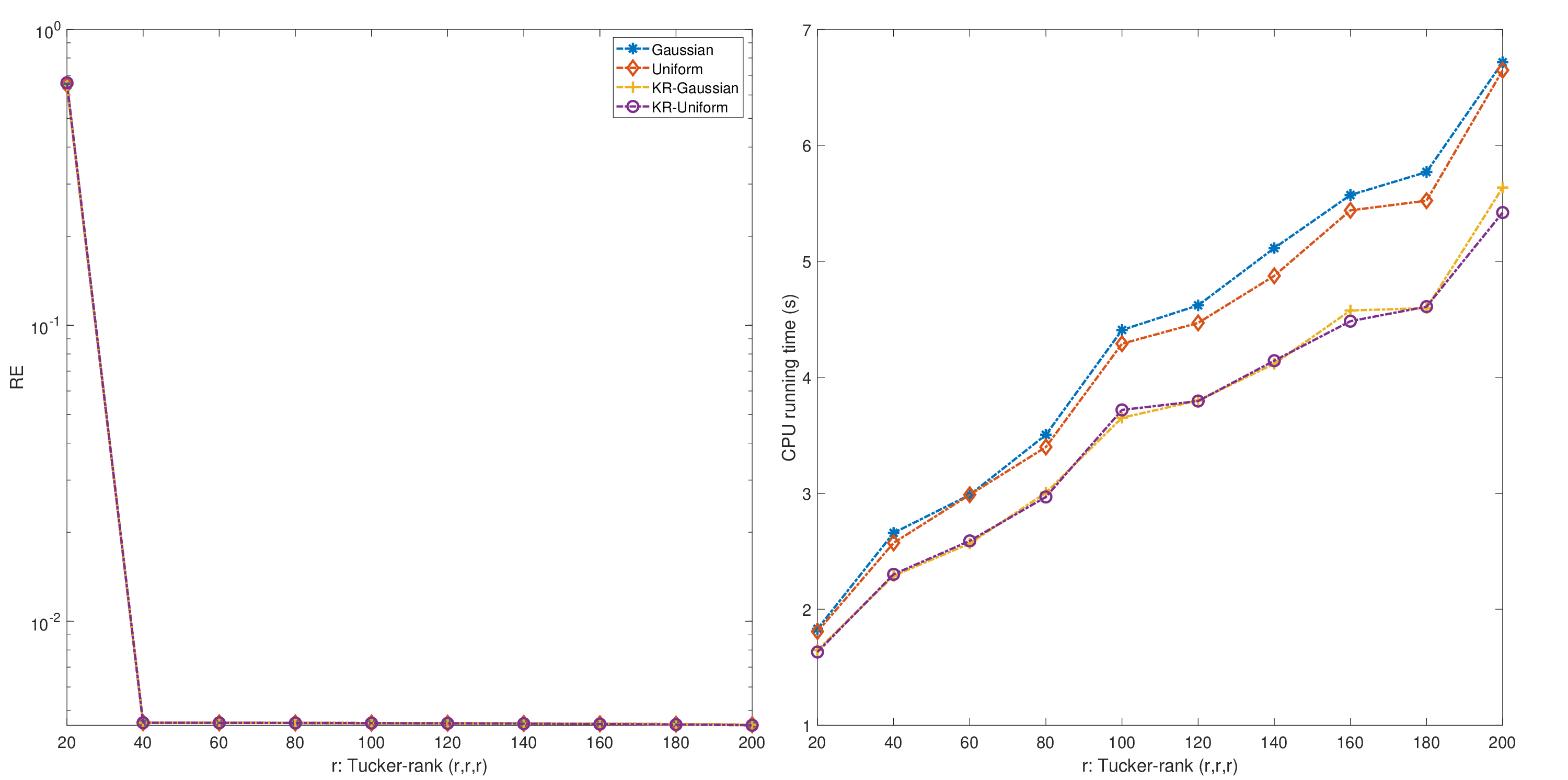}}\\
	\caption{
        For fixed $s_1=s_2=s_3=10$ and $q=1$, the comparison of RE and CPU run time obtained by applying Algorithm \ref{dash-rthosvd:alg2} using four types of random matrices and different Tucker-ranks $(r,r,r)$ on two test tensors in Section \ref{dash-rthosvd:sec5:subsec1}.}\label{dash-rthosvd:figure1:appsec4}
\end{figure}

\begin{figure}
	\setlength{\tabcolsep}{4pt}
	\renewcommand\arraystretch{1}
	\centering
	\subfigure[The tensor $\mathcal{B}$]{\includegraphics[width=4.5in,height=1.6in]{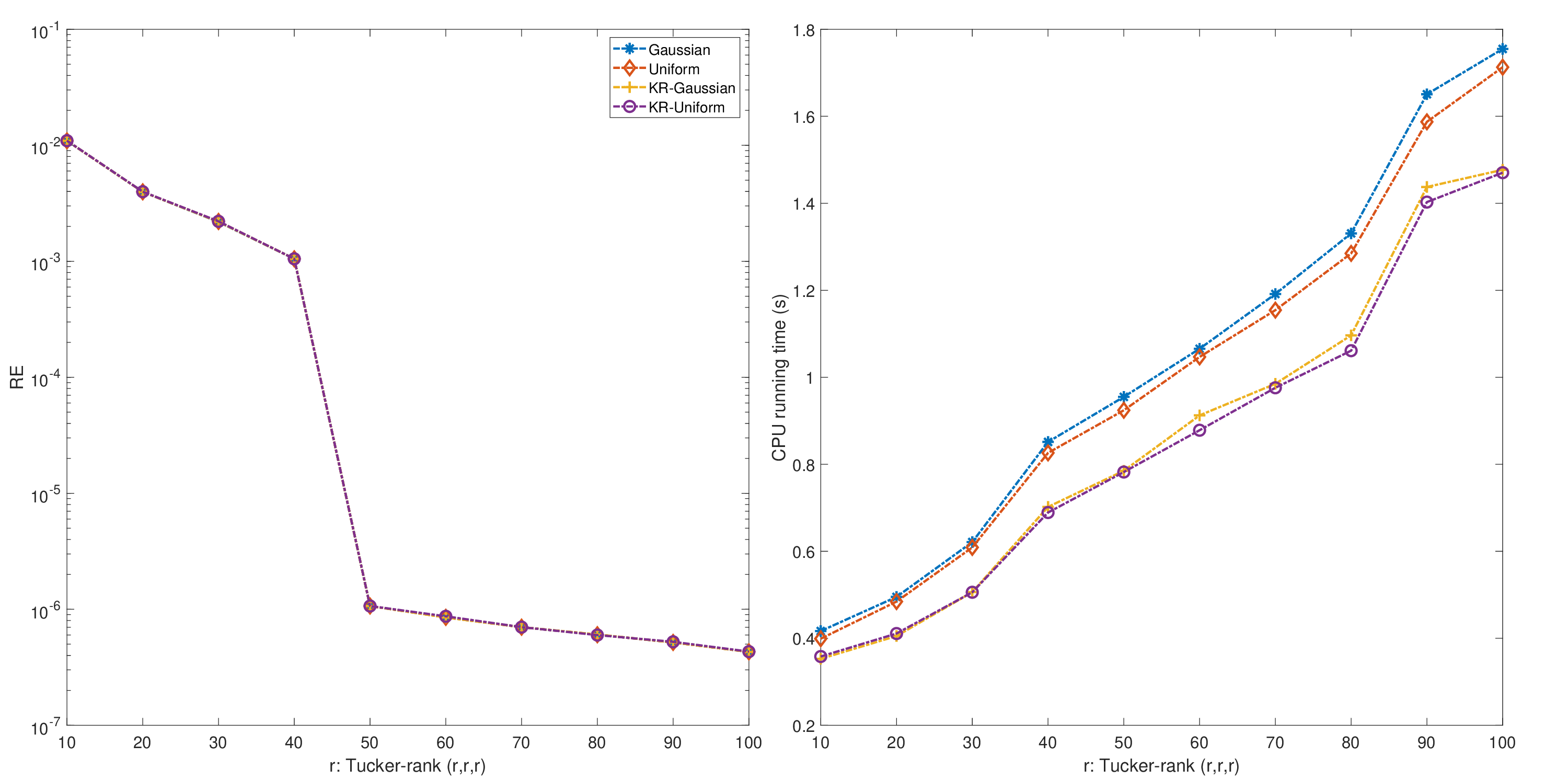}}\\
	\subfigure[The tensor $\mathcal{C}$ with Slow decay]{\includegraphics[width=4.5in,height=1.6in]{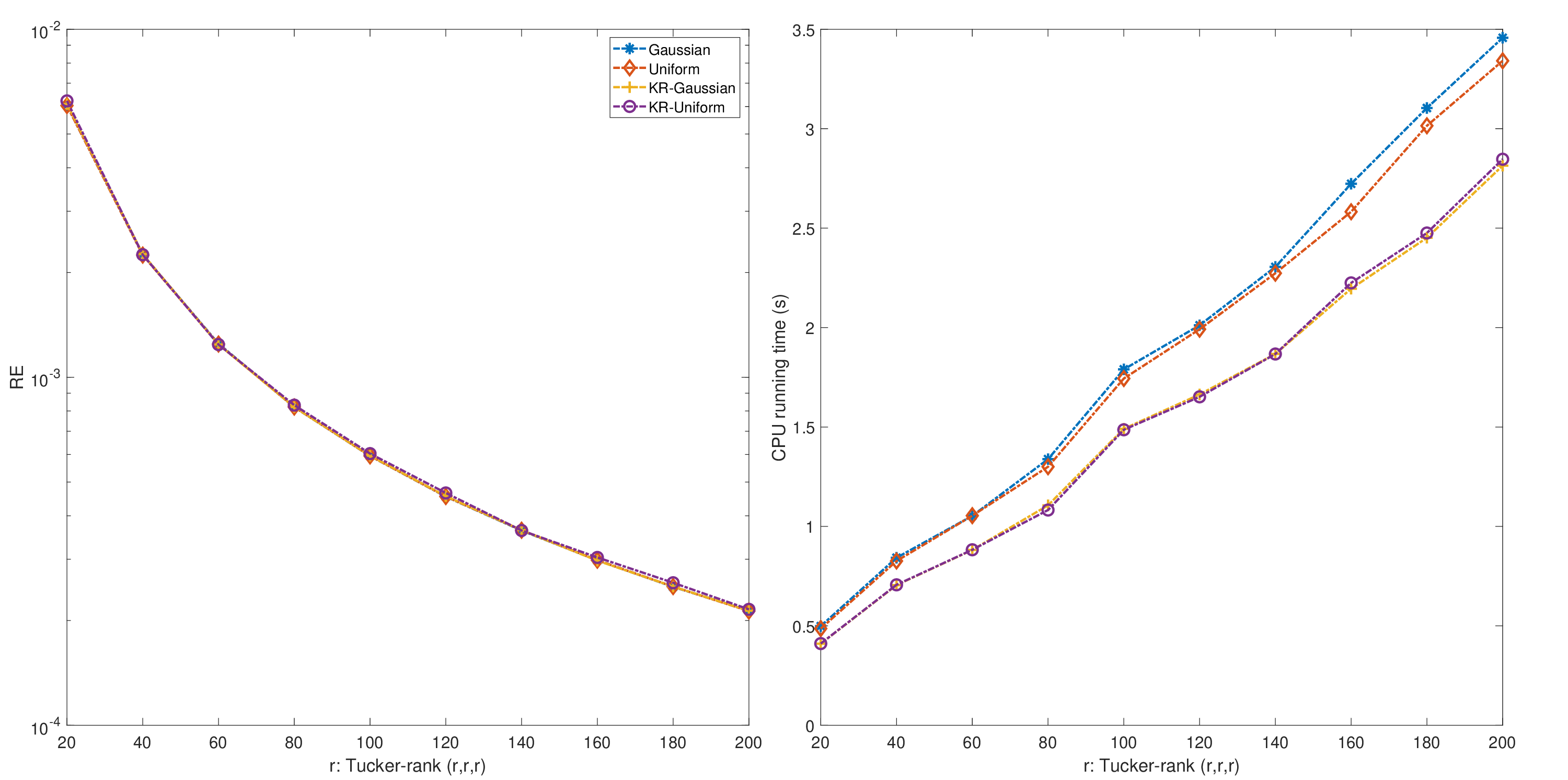}}\\
	\subfigure[The tensor $\mathcal{C}$ with Fast decay]{\includegraphics[width=4.5in,height=1.6in]{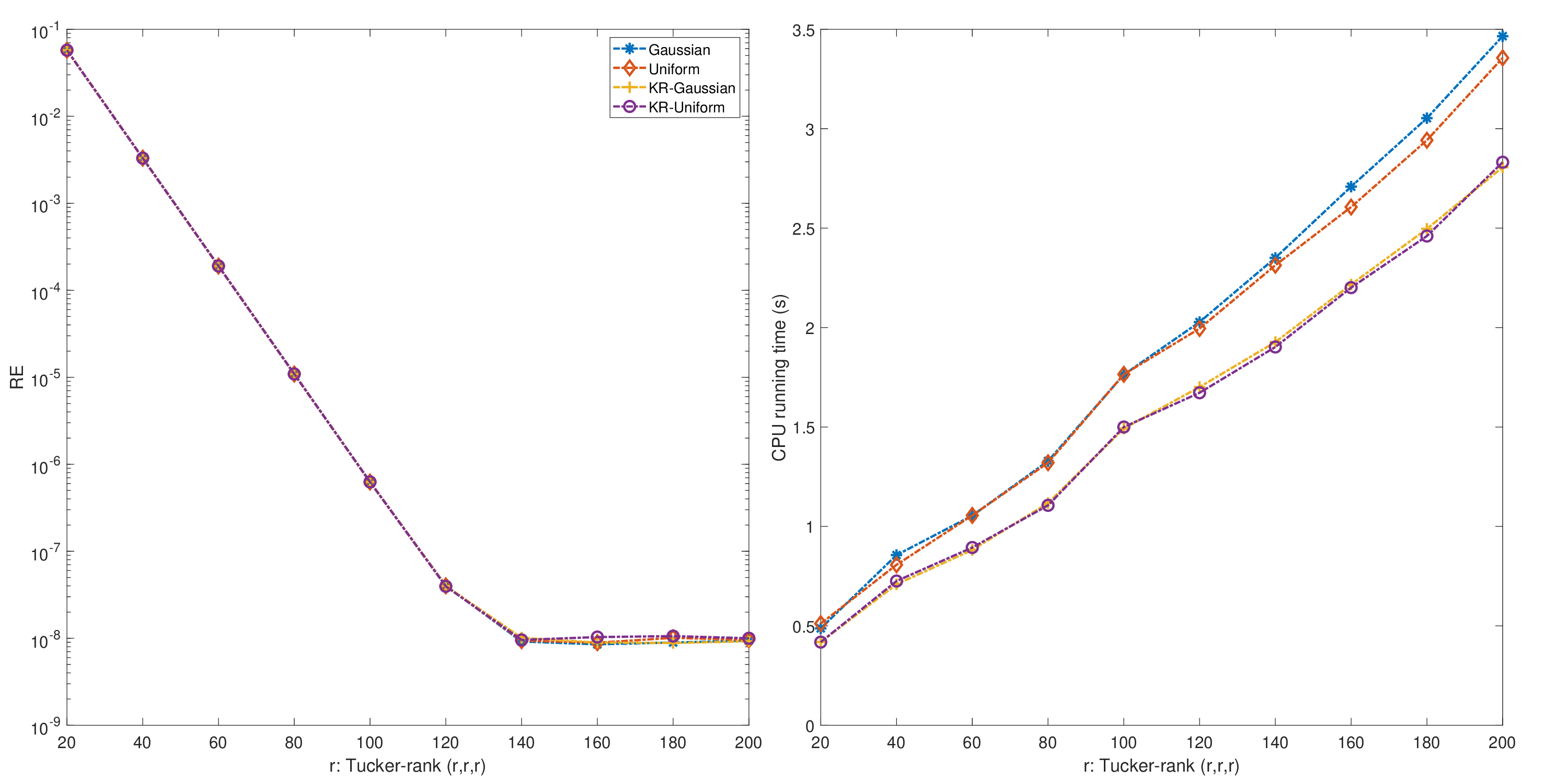}}\\
    \subfigure[The tensor $\mathcal{C}$ with S-shape decay]{\includegraphics[width=4.5in,height=1.6in]{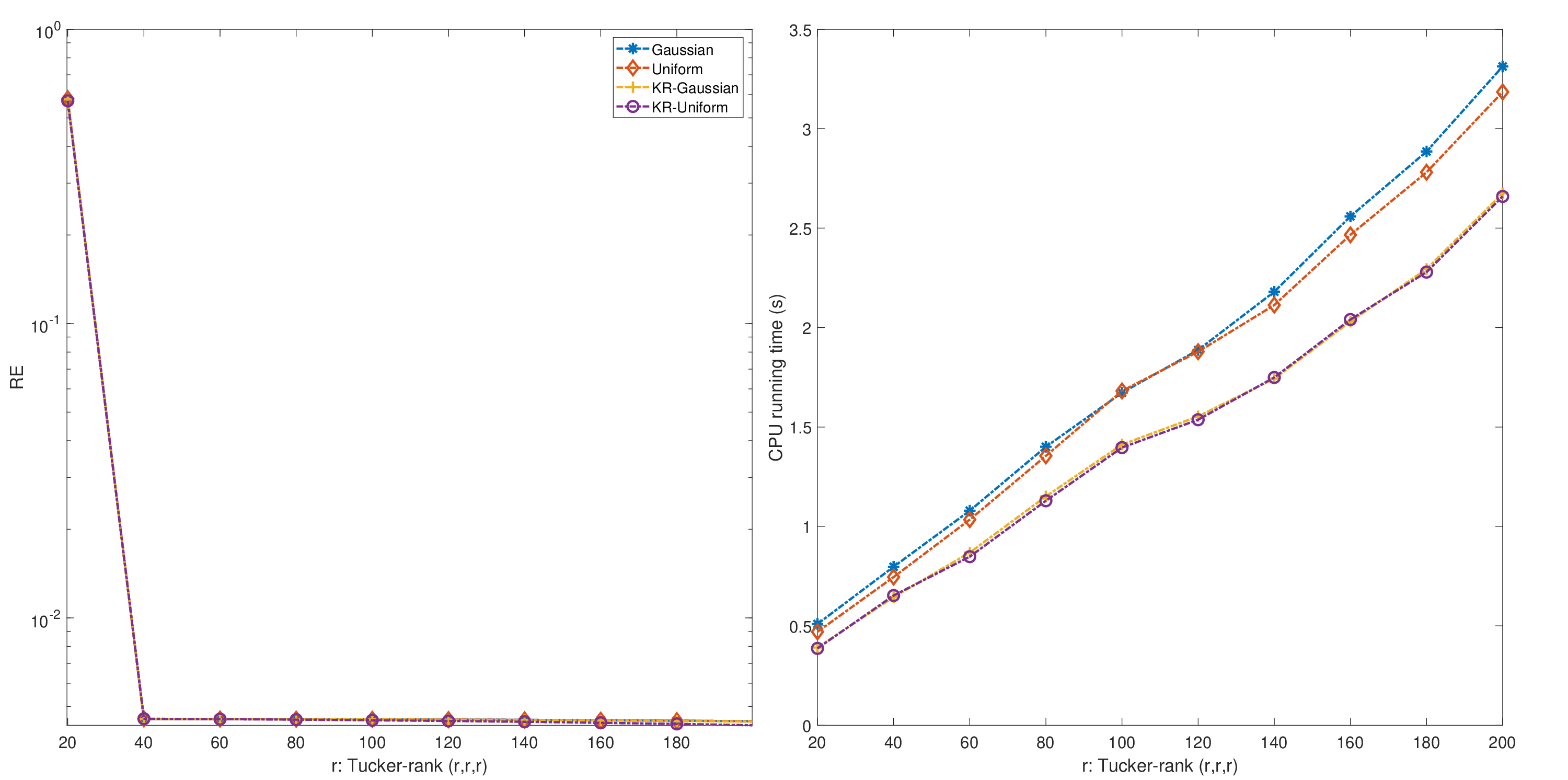}}\\
	\caption{
        For fixed $s_1=s_2=s_3=10$ and $q=1$, the comparison of RE and CPU run time obtained by applying Algorithm \ref{dash-rthosvd:alg2:v2} using four types of random matrices and different Tucker-ranks $(r,r,r)$ on two test tensors in Section \ref{dash-rthosvd:sec5:subsec1}.}\label{dash-rthosvd:figure2:appsec4}
\end{figure}

When fixing $s_1=s_2=s_3=10$ and $q=1$, we compare the efficiency of Algorithms \ref{dash-rthosvd:alg2} and \ref{dash-rthosvd:alg2:v2} using different random matrices via two test tensors $\mathcal{B}$ and $\mathcal{C}$ (see Section \ref{dash-rthosvd:sec5:subsec1}). The related results, obtained by Algorithms \ref{dash-rthosvd:alg2} and \ref{dash-rthosvd:alg2:v2}, are, respectively, illustrated in Figures \ref{dash-rthosvd:figure1:appsec4} and \ref{dash-rthosvd:figure2:appsec4}. From these two figures, we can see that Algorithms \ref{dash-rthosvd:alg2} and \ref{dash-rthosvd:alg2:v2} using Gaussian, Uniform, KR-Gaussian and KR-Uniform are comparable in terms of RE, and Algorithms \ref{dash-rthosvd:alg2} and \ref{dash-rthosvd:alg2:v2} using KR-Gaussian and KR-Uniform are faster than Algorithms \ref{dash-rthosvd:alg2} and \ref{dash-rthosvd:alg2:v2} using Gaussian and Uniform.
\section{Conclusions}
\label{dash-rthosvd:sec6}
In this paper, we discussed efficient randomized variants of T-HOSVD and ST-HOSVD for the fixed-rank approximation in Tucker format based on the integration of adaptive shifted power iterations. We also obtained an error bound for the term $\|\mathcal{G}\times_1\mathbf{U}_1\times_2\mathbf{U}_2\dots\times_d\mathbf{U}_d-\mathcal{A}\|_F^2$, where $\{\mathcal{G},\mathbf{U}_1,\dots,\mathbf{U}_d\}$ was obtained by applying the proposed algorithms to $\mathcal{A}$. For Algorithm \ref{dash-rthosvd:alg2}, when the power parameter tends to infinity, the error bound is consistent with that obtained from T-HOSVD and ST-HOSVD (see Theorem \ref{dash-rthosvd:general:errror}). Numerical examples illustrated that Algorithms \ref{dash-rthosvd:alg2} and \ref{dash-rthosvd:alg2:v2} with $q=1$ were comparable to the existing algorithms in terms of RE and in terms of CPU run time, Algorithm \ref{dash-rthosvd:alg2} was faster than Tucker-ALS, T-HOSVD and their randomized variants, and Algorithm \ref{dash-rthosvd:alg2:v2} was even faster than Algorithm \ref{dash-rthosvd:alg2}, Tucker-ALS, T-HOSVD, ST-HOSVD and their randomized variants.

\section*{Acknowledgements}
The authors would like to thank the associate editor and one anonymous reviewer for their careful and very detailed comments on our paper.


\appendix
\renewcommand\thealgorithm{\Alph{section}.\arabic{algorithm}}
\setcounter{algorithm}{0}
\renewcommand\thefigure{\Alph{section}.\arabic{figure}}
\setcounter{figure}{0}
\renewcommand\thetable{\Alph{section}.\arabic{table}}
\setcounter{table}{0}

\section{The existing randomized variants for T-HOSVD and ST-HOSVD with a fixed Tucker-rank}
\label{dash-rthosvd:appsec1}
We present one model of randomized variants for T-HOSVD and ST-HOSVD with a fixed Tucker-rank as in Algorithm \ref{dash-rthosvd:alg1}.
\begin{algorithm}[htb]
     \caption{One model for randomized variants of T-HOSVD/ST-HOSVD with a fixed Tucker-rank}
     \begin{algorithmic}[1]
        \STATEx {\bf Input}: A tensor $\mathcal{A}\in\mathbb{R}^{n_1\times n_2\times\dots\times n_d}$, the desired $d$-tuple $\mathbf{r}=(r_1,r_2,\dots,r_d)$ with $r_k<n_k$ and $k=1,2,\dots,d$, the $d$-tuple of oversampling parameters $\mathbf{s}=(s_1,s_2,\dots,s_d)$, and the power parameter $q\geq 1$.
        \STATEx {\bf Output}: The core tensor $\mathcal{G}\in\mathbb{R}^{r_1\times r_2\times\dots\times r_d}$ and the mode-$k$ factor matrix $\mathbf{U}_k\in\mathbb{R}^{n_k\times r_k}$ such that $\mathcal{A}\approx\widehat{\mathcal{A}}:=\mathcal{G}\times_1\mathbf{U}_1\times_2\mathbf{U}_2\dots\times_d\mathbf{U}_d$.
        \STATE Let $l_k=r_k+s_k$ with $k=1,2,\dots,d$ and set a temporary tensor $\mathcal{G}:=\mathcal{A}$.
        \FOR{$k=1,2,\dots,d$}
            \IF{in T-HOSVD case}
               \STATE Draw a standard Gaussian matrix $\mathbf{\Omega}_k\in\mathbb{R}^{n_1\dots n_{k-1}n_{k+1}\dots n_{d}\times l_k}$.
               \STATE Form the mode-$k$ unfolding $\mathbf{A}_{(k)}$ from $\mathcal{A}$.
               \STATE Compute $[\mathbf{Q}_k,\sim,\sim]={\rm svd}(\mathbf{A}_{(k)}\mathbf{\Omega}_k,'{\rm econ}')$.
               \FOR{$t=1,2,\dots,q$}
                   \STATE Compute $[\mathbf{Q}_k,\sim,\sim]={\rm svd}(\mathbf{A}_{(k)}(\mathbf{A}_{(k)}^\top\mathbf{Q}_k),'{\rm econ}')$.
               \ENDFOR
            \ELSIF{in ST-HOSVD case}
               \STATE Draw a random matrix $\mathbf{\Omega}_k\in\mathbb{R}^{r_1\dots r_{k-1}n_{k+1}\dots n_{d}\times l_k}$.
               \STATE Form the mode-$k$ unfolding $\mathbf{G}_{(k)}$ from $\mathcal{G}$.
               \STATE Compute $[\mathbf{Q}_k,\sim,\sim]={\rm svd}(\mathbf{G}_{(k)}\mathbf{\Omega}_k,'{\rm econ}')$.
               \FOR{$t=1,2,\dots,q$}
                   \STATE Compute $[\mathbf{Q}_k,\sim,\sim]={\rm svd}(\mathbf{G}_{(k)}(\mathbf{G}_{(k)}^\top\mathbf{Q}_k),'{\rm econ}')$.
               \ENDFOR
            \ENDIF
            \STATE Set $\mathbf{U}_k=\mathbf{Q}_k(:,1:r_k)$ and update $\mathcal{G}=\mathcal{G}\times_k\mathbf{U}_{k}^\top$.
        \ENDFOR
    \end{algorithmic}
    \label{dash-rthosvd:alg1}
\end{algorithm}
\begin{remark}
    Algorithm \ref{dash-rthosvd:alg1} is similar to the work in \cite{ahmadi2021randomized,che2020computation,che2021efficient,che2021randomized,che2023randomized,minster2020randomized,qiu2024towards}. Their differences relate to:
    \begin{enumerate}
        \item[(a)] for each $k$, the matrix $\mathbf{\Omega}_k$ in \cite{che2020computation,che2021randomized} is the Kronecker product of standard Gaussian matrices, the matrix $\mathbf{\Omega}_k$ in \cite{che2021efficient} is the Kronecker product of SRFT (subsampled randomized Fourier transform) matrices, and the matrix $\mathbf{\Omega}_k$ in \cite{che2023randomized} is the Kronecker product of SpEmb (sparse subspace embedding) matrices;
        \item[(b)] the power parameter $q$ in \cite{che2021randomized} is set to 0;
        \item[(c)] the way to obtain $\mathbf{U}_k$ in Algorithm \ref{dash-rthosvd:alg1} is the same as that in \cite{che2020computation,che2021randomized,che2023randomized}, the matrix $\mathbf{U}_k$ in \cite{che2021efficient} is obtained as $\mathbf{U}_k=\mathbf{Q}_k(:,1:r_k)$, and according to \cite{ahmadi2021randomized,minster2020randomized,qiu2024towards}, the matrix $\mathbf{U}_k$ is given by $\mathbf{U}_k=\mathbf{Q}_k\mathbf{V}_k$, where $\mathbf{Q}_k={\rm orth}((\mathbf{G}_{(k)}\mathbf{G}_{(k)}^\top)^q\mathbf{G}_{(k)}\mathbf{\Omega}_k)$ and all the columns of $\mathbf{V}_k$ are the first $r_k$ left singular vectors of $\mathbf{Q}_k^\top\mathbf{A}_{(k)}$.
    \end{enumerate}
\end{remark}

Similar to the randomized projection algorithms for the fixed-rank approximation in Tucker format, except for the truncation approach, Algorithm \ref{dash-rthosvd:alg1:app} is another basic form of randomized variants of T-HOSVD and ST-HOSVD with a fixed Tucker-rank, which is modified from Minster {\it et al.} \cite{minster2024parallel}. Instead of directly applying the randomized SVD algorithm (see \cite{halko2011finding}) to each mode unfolding matrix, as in Algorithm \ref{dash-rthosvd:alg1}, a holistic approach was used.

\begin{algorithm}[htb]
     \caption{Another model for randomized variants of T-HOSVD and ST-HOSVD with power scheme}
     \begin{algorithmic}[1]
        \STATEx {\bf Input}: A tensor $\mathcal{A}\in\mathbb{R}^{n_1\times n_2\times\dots\times n_d}$, the desired $d$-tuple $\mathbf{r}=(r_1,r_2,\dots,r_d)$ with $r_k< n_k$ and $k=1,2,\dots,d$, the $d$-tuple of oversampling parameters $\mathbf{s}=(s_1,s_2,\dots,s_d)$, and the power parameter $q\geq 1$.
        \STATEx {\bf Output}: The core tensor $\mathcal{G}\in\mathbb{R}^{r_1\times r_2\times\dots\times r_d}$ and the mode-$k$ factor matrix $\mathbf{U}_k\in\mathbb{R}^{n_k\times r_k}$ such that $\mathcal{A}\approx\widehat{\mathcal{A}}:=\mathcal{G}\times_1\mathbf{U}_1\dots\times_d\mathbf{U}_d$.
        \STATE Let $l_k=r_k+s_k$ with $k=1,2,\dots,d$ and set a temporary tensor $\mathcal{B}:=\mathcal{A}$.
        \FOR{$k=1,2,\dots,d$}
           \IF{in T-HOSVD case}
               \STATE Draw a standard Gaussian matrix $\mathbf{\Omega}_k\in\mathbb{R}^{n_1\dots n_{k-1}n_{k+1}\dots n_{d}\times l_k}$.
               \STATE Form the mode-$k$ unfolding $\mathbf{A}_{(k)}$ from $\mathcal{A}$.
               \STATE Compute $\mathbf{Q}_k={\rm orth}(\mathbf{A}_{(k)}\mathbf{\Omega}_k)$.
               \FOR{$t=1,2,\dots,q$}
                  \STATE Compute $\mathbf{Q}_k={\rm orth}(\mathbf{A}_{(k)}(\mathbf{A}_{(k)}^\top\mathbf{Q}_k))$.
               \ENDFOR
            \ELSIF{in ST-HOSVD case}
               \STATE Draw a standard Gaussian matrix $\mathbf{\Omega}_k\in\mathbb{R}^{l_1\dots l_{k-1}n_{k+1}\dots n_{d}\times l_k}$.
               \STATE Form the mode-$k$ unfolding $\mathbf{B}_{(k)}$ from $\mathcal{B}$.
            \STATE Compute $\mathbf{Q}_k={\rm orth}(\mathbf{B}_{(k)}\mathbf{\Omega}_k)$.
            \FOR{$t=1,2,\dots,q$}
                \STATE Compute $\mathbf{Q}_k={\rm orth}(\mathbf{B}_{(k)}(\mathbf{B}_{(k)}^\top\mathbf{Q}_k))$.
            \ENDFOR
            \ENDIF
            \STATE Update $\mathcal{B}=\mathcal{B}\times_k\mathbf{Q}_{k}^\top$.
        \ENDFOR
        \STATE Compute $\{\mathcal{G};\mathbf{V}_1,\dots,\mathbf{V}_d\}=\text{ST-HOSVD}(\mathcal{B},\mathbf{r})$.\label{dash-rthosvd:alg1:app:step1}
        \STATE Form $\mathbf{U}_k=\mathbf{Q}_{k}\mathbf{V}_k$ with $k=1,2,\dots,d$.
    \end{algorithmic}
    \label{dash-rthosvd:alg1:app}
\end{algorithm}
\begin{remark}
Step \ref{dash-rthosvd:alg1:app:step1} in Algorithm \ref{dash-rthosvd:alg1:app} can be revised as $$\{\mathcal{G};\mathbf{V}_1,\dots,\mathbf{V}_d\}=\text{T-HOSVD}(\mathcal{B},\mathbf{r}).$$
\end{remark}

\section{Another type for randomized variants of ST-HOSVD with adaptive shift}
\label{dash-rthosvd:appsec2}
Another type for randomized variants of ST-HOSVD with adaptive shift (see Algorithm \ref{dash-rthosvd:alg2:app}), similar to Algorithm \ref{dash-rthosvd:alg2:v2},  is proposed to find a quasi-optimal solution for Problem \ref{dash-rthosvd:prob3} (see Algorithm \ref{dash-rthosvd:alg2:app}).
\begin{algorithm}
     \caption{Another version for the randomized ST-HOSVD with adaptive shift}
     \begin{algorithmic}[1]
        \STATEx {\bf Input}: A tensor $\mathcal{A}\in\mathbb{R}^{n_1\times n_2\times\dots\times n_d}$, the desired $d$-tuple $\mathbf{r}=(r_1,r_2,\dots,r_d)$ with $r_k\leq n_k$ and $k=1,2,\dots,d$, the $d$-tuple of oversampling parameters $\mathbf{s}=(s_1,s_2,\dots,s_d)$, and the power parameter $q\geq 1$.
        \STATEx {\bf Output}: The core tensor $\mathcal{G}\in\mathbb{R}^{r_1\times r_2\times\dots\times r_d}$ and the mode-$k$ factor matrix $\mathbf{U}_k\in\mathbb{R}^{n_k\times r_k}$ such that $\mathcal{A}\approx\widehat{\mathcal{A}}:=\mathcal{G}\times_1\mathbf{U}_1\dots\times_d\mathbf{U}_d$.
        \STATE Let $l_k=r_k+s_k$ with $k=1,2,\dots,d$ and set a temporary tensor $\mathcal{B}:=\mathcal{A}$.
        \FOR{$k=1,2,\dots,d$}
            \STATE Draw a standard Gaussian matrix $\mathbf{\Omega}_k\in\mathbb{R}^{l_1\dots l_{k-1}n_{k+1}\dots n_{d}\times l_k}$.
            \STATE Form the mode-$k$ unfolding $\mathbf{B}_{(k)}$ from $\mathcal{B}$.
            \STATE Compute $\mathbf{Q}_k={\rm orth}(\mathbf{B}_{(k)}\mathbf{\Omega}_k)$ and set {\color{red}$\alpha_{0}^{(k)}=0$}.
            \FOR{$t=1,2,\dots,q$}
                \STATE Compute $[\mathbf{Q}_k,\bm{\Sigma}_k,\sim]={\rm svd}(\mathbf{B}_{(k)}(\mathbf{B}_{(k)}^\top\mathbf{Q}_k)-{\color{red}\alpha_{t-1}^{(k)}} \mathbf{Q}_k,'{\rm econ}')$.
                \IF{$\bm{\Sigma}_k(r_k+p_k,r_k+p_k)>{\color{red}\alpha_{t-1}^{(k)}}$}
                    \STATE Update ${\color{red}\alpha_{t}^{(k)}}=(\mathbf{\Sigma}_k(r_k+p_k,r_k+p_k)+{\color{red}\alpha_{t-1}^{(k)}})/2$.
                \ENDIF
            \ENDFOR
            \STATE Update $\mathcal{B}=\mathcal{B}\times_k\mathbf{Q}_{k}^\top$.
        \ENDFOR
        \STATE Compute $\{\mathcal{G};\mathbf{V}_1,\dots,\mathbf{V}_d\}=\text{ST-HOSVD}(\mathcal{B},\mathbf{r})$.
        \STATE Form $\mathbf{U}_k=\mathbf{Q}_{k}\mathbf{V}_k$ with $k=1,2,\dots,d$.
    \end{algorithmic}
    \label{dash-rthosvd:alg2:app}
\end{algorithm}

\begin{remark}
We can obtain another version of randomized T-HOSVD with adaptive shift based on Algorithm 4.1 in \cite{minster2024parallel}, that is similar to Algorithm \ref{dash-rthosvd:alg2:app} (details omitted). For this case, the corresponding algorithms are denoted as Algorithm \ref{dash-rthosvd:alg2}-v1 and Algorithm \ref{dash-rthosvd:alg2}-v2, and Algorithms \ref{dash-rthosvd:alg2:v2} and \ref{dash-rthosvd:alg2:app} are denoted as Algorithm \ref{dash-rthosvd:alg2:v2}-v1 and Algorithm \ref{dash-rthosvd:alg2:v2}-v2, respectively.
\end{remark}

\begin{figure}
	\setlength{\tabcolsep}{4pt}
	\renewcommand\arraystretch{1}
	\centering
    \subfigure[The tensor $\mathcal{B}$]{\includegraphics[width=4.5in,height=1.6in]{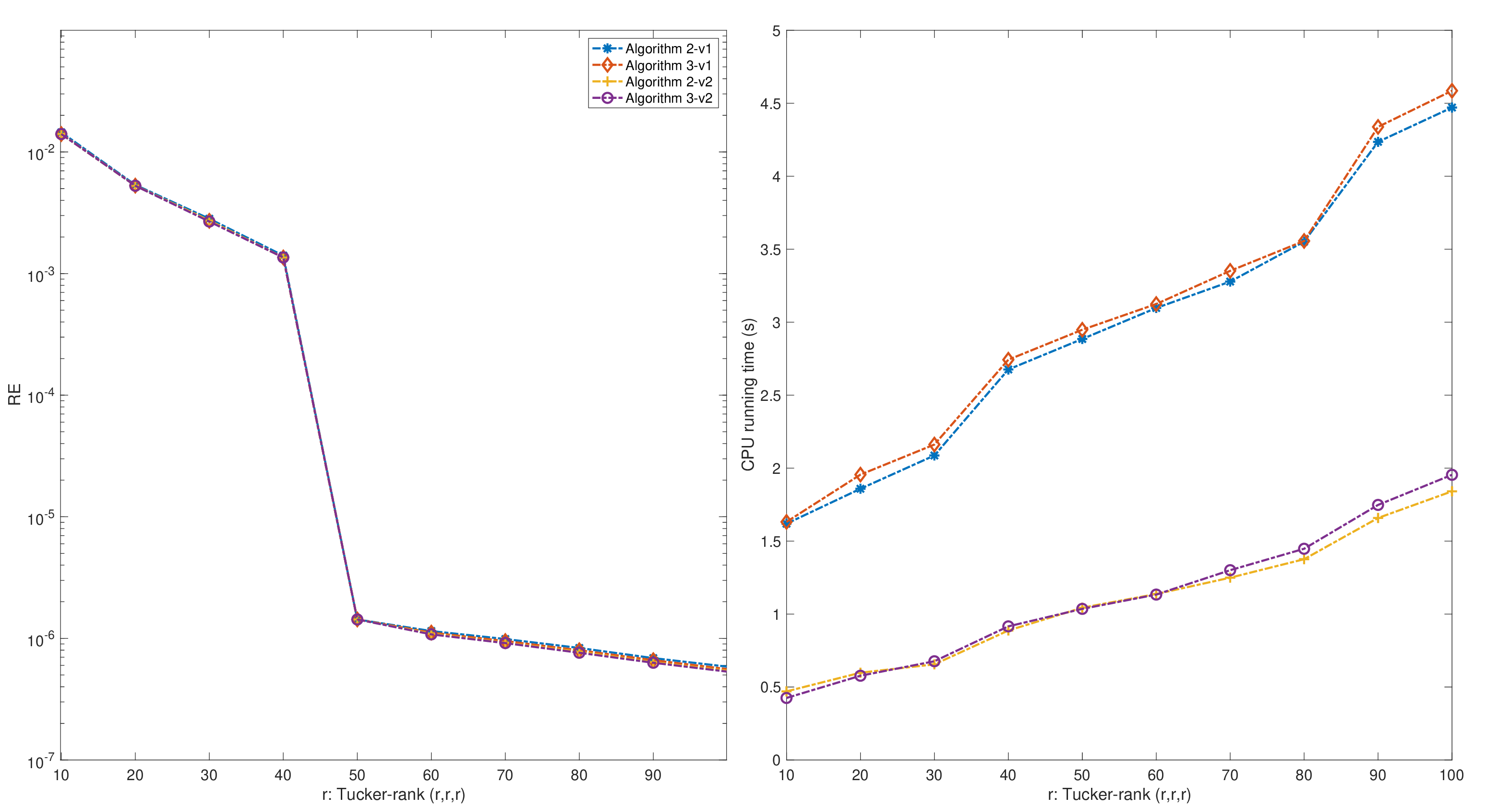}}\\
	\subfigure[The tensor $\mathcal{C}$ with Slow decay]{\includegraphics[width=4.5in,height=1.6in]{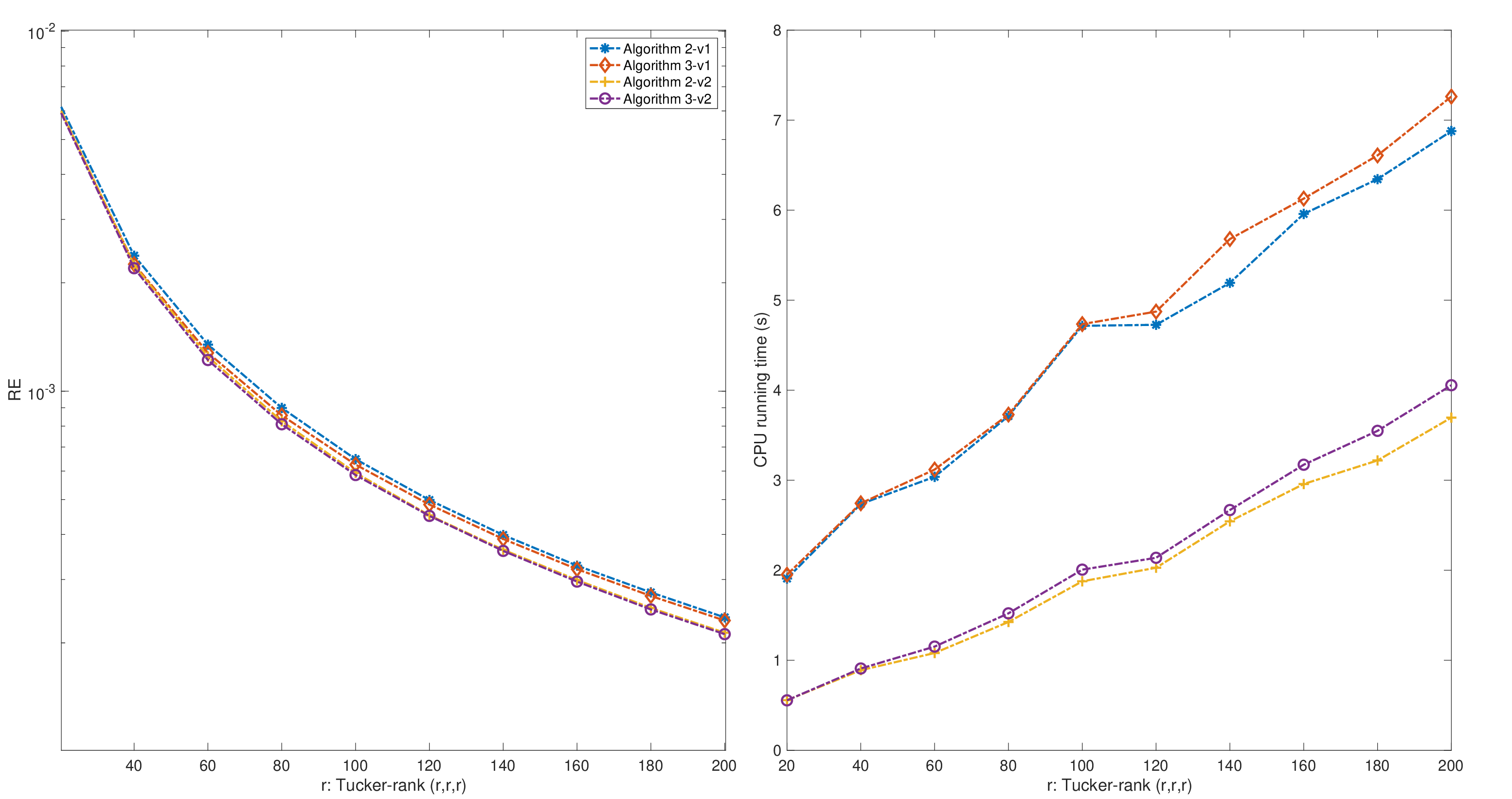}}\\
	\subfigure[The tensor $\mathcal{C}$ with Fast decay]{\includegraphics[width=4.5in,height=1.6in]{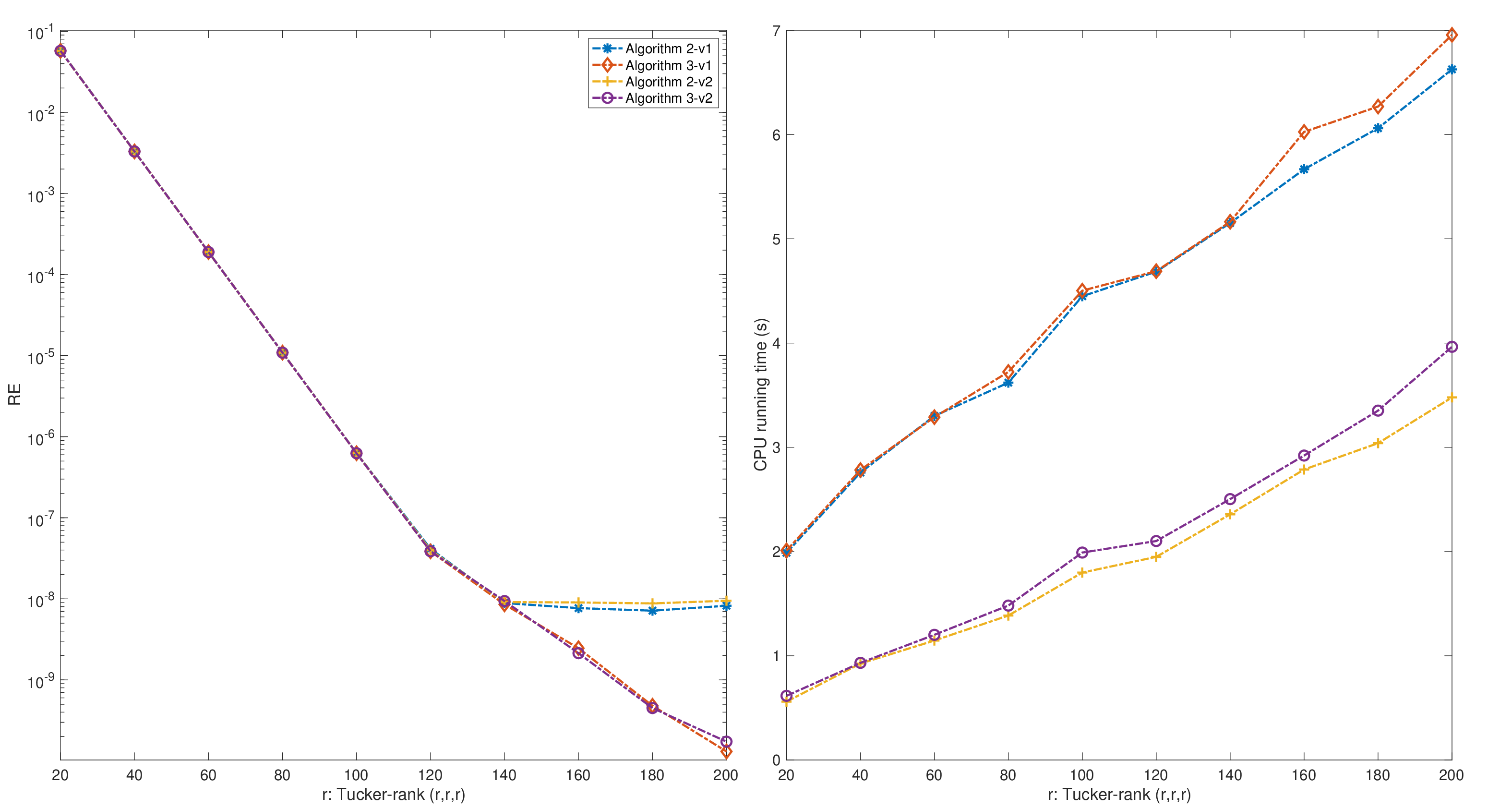}}\\
    \subfigure[The tensor $\mathcal{C}$ with S-shape decay]{\includegraphics[width=4.5in,height=1.6in]{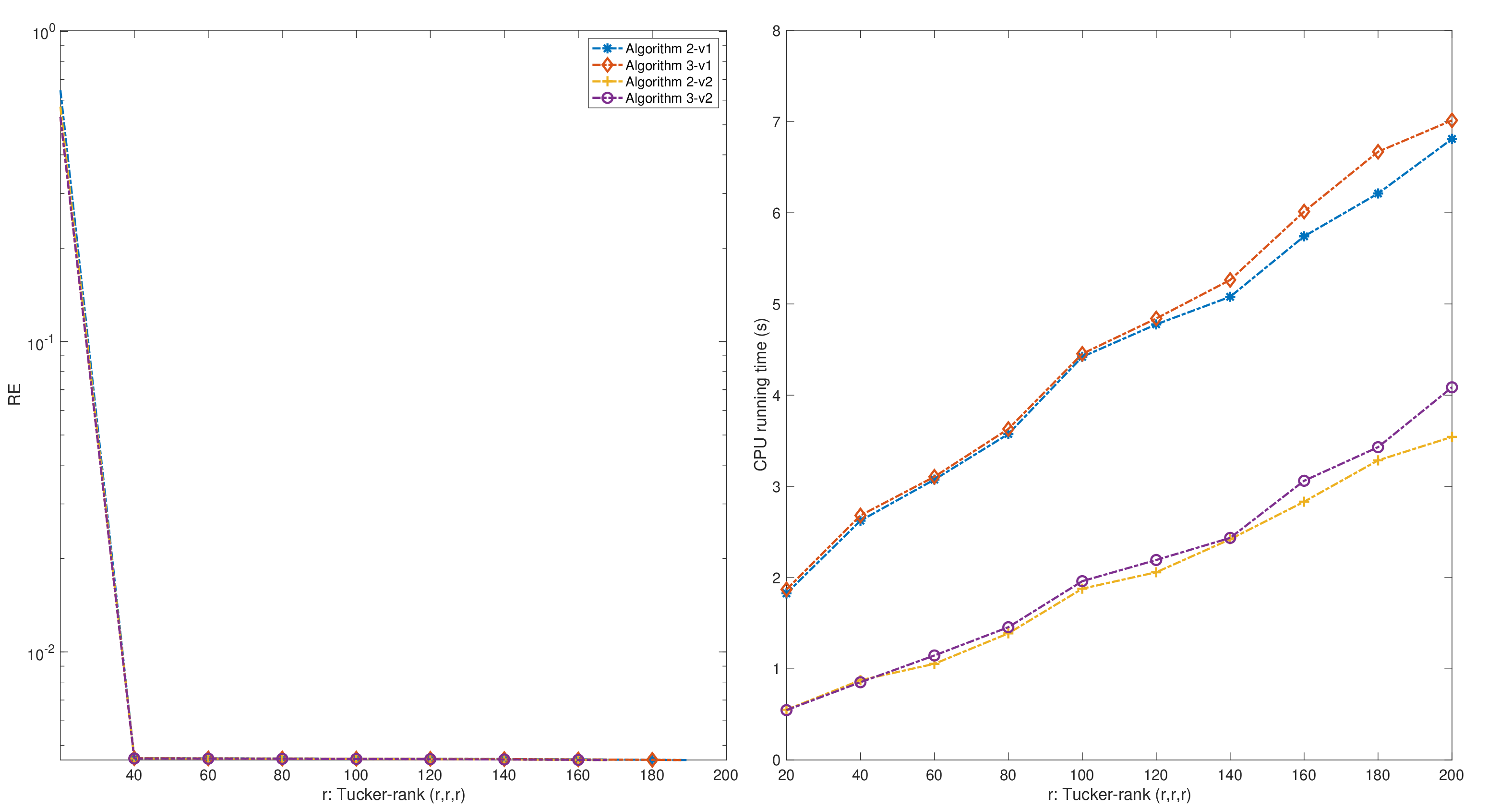}}\\
	\caption{
        For fixed $s_1=s_2=s_3=10$ and $q=1$, for different 3-tuples $(r,r,r)$, the values of RE and CPU run time by applying Algorithms \ref{dash-rthosvd:alg2}-v1/v2 and Algorithms \ref{dash-rthosvd:alg2:v2}-v1/v2 to two test tensors in Section \ref{dash-rthosvd:sec5:subsec1}.}\label{dash-rthosvd:figure:appsec2}
\end{figure}

We now compare the efficiency of Algorithms \ref{dash-rthosvd:alg2}-v1/v2 and Algorithms \ref{dash-rthosvd:alg2:v2}-v1/v2 via two test tensors $\mathcal{B}$ and $\mathcal{C}$ (see Section \ref{dash-rthosvd:sec5:subsec1}). In these algorithms, we set $s_1=s_2=s_3=10$ and $q=1$. The related results for RE and CPU run time are shown in Figure \ref{dash-rthosvd:figure:appsec2}These four algorithms are comparable in terms of RE, and in terms of CPU run time, Algorithms \ref{dash-rthosvd:alg2:v2}-v1/v2 were faster than Algorithms \ref{dash-rthosvd:alg2}-v1/v2, and Algorithm \ref{dash-rthosvd:alg2}-v1 and Algorithm \ref{dash-rthosvd:alg2:v2}-v2 were slightly faster than Algorithm \ref{dash-rthosvd:alg2}-v2 and Algorithm \ref{dash-rthosvd:alg2:v2}-v1, respectively.

\bibliographystyle{amsplain}
\bibliography{reference}

\end{document}